\documentclass{amsart}
\usepackage{amsrefs}
\usepackage{color}
\usepackage{epsfig}
\usepackage{amsmath}
\usepackage{amssymb}
\usepackage{amsfonts}
\usepackage{euscript}
\usepackage[all]{xy}
\usepackage{mathrsfs}
\usepackage{setspace}
\usepackage{filecontents}

\numberwithin{equation}{section}
\newcommand{\bR}{{\mathbb R}}

\newcommand{\cF}{{\mathcal F}}

\newcommand{\smooth}{C^\infty}
\newcommand{\R}{\mathbb{R}}
\newcommand{\C}{\mathbb C}
\newcommand{\Z}{\mathbb Z}
\newcommand{\K}{\mathbb K}

\renewcommand{\K}{\mathbb{K}}

\newcommand{\scrF}{\EuScript F}
\newcommand{\scrR}{\EuScript R}
\newcommand{\scrS}{\EuScript S}
\newcommand{\scrM}{\EuScript M}

\newcommand{\scrC}{\EuScript C}
\newcommand{\scrU}{\EuScript U}
\newcommand{\scrX}{\EuScript X}
\newcommand{\scrY}{\EuScript Y}
\newcommand{\scrQ}{\EuScript Q}
\newcommand{\scrK}{\EuScript K}
\newcommand{\scrJ}{\EuScript J}
\newcommand{\scrP}{\EuScript P}
\newcommand{\scrW}{\EuScript{W}}

\newcommand{\bfepsilon}{\boldsymbol \epsilon}
\newcommand{\bfphi}{\boldsymbol \phi}
\newcommand{\bfzeta}{\boldsymbol \zeta}

\newcommand{\bfp}{\mathbf p}
\newcommand{\bfF}{\mathbf F}
\newcommand{\bfI}{\mathbf I}
\newcommand{\bfQ}{\mathbf Q}
\newcommand{\bfq}{\mathbf q}
\newcommand{\bfalpha}{\boldsymbol\alpha}
\newcommand{\bfbeta}{\boldsymbol\beta}
\newcommand{\bfJ}{\mathbf J}
\newcommand{\bfK}{\mathbf K}
\newcommand{\bfP}{\mathbf P}

\newcommand{\bfx}{\mathbf{x}}
\newcommand{\bfY}{\mathbf{Y}}
\newcommand{\bfw}{\mathbf{w}}
\newcommand{\bfsigma}{\boldsymbol\sigma}
\newcommand{\bfpsi}{\boldsymbol\psi}
\newcommand{\bfL}{\mathbf{L}}

\newcommand{\Sym}{\mathit{Sym}}

\renewcommand{\bigsqcup}{\coprod}
\newcommand{\leftsc}{\langle}
\newcommand{\rightsc}{\rangle}

\newcommand{\half}{\textstyle\frac{1}{2}}
\newcommand{\im}{\mathrm{im}}
\newcommand{\iso}{\cong}
\renewcommand{\deg}{\mathrm{deg}}
\newcommand{\bardeg}{\overline{\deg}}

\newtheorem{thm}{Theorem}[section]
\newtheorem{theorem}[thm]{Theorem}
\newtheorem{addendum}[thm]{Addendum}
\newtheorem{cor}[thm]{Corollary}

\newtheorem{prop}[thm]{Proposition}
\newtheorem{definition}[thm]{Definition}
\newtheorem{remark}[thm]{Remark}
\newtheorem{convention}[thm]{Convention}
\newtheorem{example}[thm]{Example}
\newtheorem{examples}[thm]{Examples}
\newtheorem{lemma}[thm]{Lemma}

\title[Viterbo functoriality]{An open string analogue\\ of Viterbo functoriality}

\author[M.~Abouzaid, P.~Seidel]{Mohammed Abouzaid, Paul Seidel} \date{December 15, 2007, revised August 25 2009}
\thanks{This research was conducted during the period the first author served as a Clay Research Fellow. The second author was partially supported by NSF grants DMS-0405516 and DMS-0652620.}
\begin{document}

\maketitle
\tableofcontents

\section{Introduction}

Symplectic (co)homology is an invariant of symplectic manifolds with boundary, introduced by Cieliebak-Floer-Hofer-Wysocki \cite{floer-hofer94, cieliebak-floer-hofer95, floer-hofer-wysocki95, cieliebak-floer-hofer-wysocki96} and Viterbo \cite{viterbo97a,viterbo97b}. There are actually several flavours of the theory, which appear in different contexts; even when the domains of applicability overlap, their mutual relation is not always clear. In this paper, the symplectic manifolds we consider are {\em Liouville domains}, which are particularly simple to work with, since they combine two advantageous features (convex contact type boundary, and exactness of the symplectic form) and we use Viterbo's definition (except that we call cohomology what he calls homology). The notation will be $SH^*(M)$ for the symplectic cohomology of a Liouville domain $M$. Crucially, for every Liouville subdomain $M^{in} \subset M$, Viterbo defined a restriction (in his terminology, transfer) homomorphism
\begin{equation} \label{eq:viterbo-closed}
SH^*(M) \longrightarrow SH^*(M^{in}).
\end{equation}
The original motivation came from the Weinstein conjecture;  this {\em Viterbo functoriality} property of symplectic cohomology also has many applications to Lagrangian embedding questions, especially for cotangent bundles. There are several surveys covering aspects of this material, for instance \cite{oancea04b,seidel07}.

Symplectic cohomology has an open string analogue, where one considers a Lagrangian submanifold $L  \subset M$ (whose boundary is Legendrian, and with suitable exactness conditions). The special case of cotangent fibres appears in \cite{abbondandolo-schwarz06}, and the general construction is known to specialists, even though it has not received wide attention. Following \cite{fukaya-seidel-smith07b} we call the resulting invariant {\em wrapped Floer cohomology}, and denote it by $HW^*(L)$. In analogy with \eqref{eq:viterbo-closed} one expects, under suitable assumptions, to have a map
\begin{equation} \label{eq:viterbo-open}
HW^*(L) \longrightarrow HW^*(L^{in})
\end{equation}
where $L^{in} = L \cap M^{in}$; indeed, a fairly straightforward translation of the argument in \cite{viterbo97a} shows that this can be done. The aim of the present paper is refine that idea, by first constructing an $A_\infty$-algebra structure on the chain complex underlying $HW^*(L)$, and then proving that \eqref{eq:viterbo-open} is induced by an $A_\infty$-homomorphism. In fact, we also provide an extension to several Lagrangian submanifolds, where wrapped Floer cohomology turns into the {\em wrapped Fukaya $A_\infty$-category} $\scrW(M)$. In that case, Viterbo functoriality takes on the form of an $A_\infty$-functor defined on a suitable full subcategory of $\scrW(M)$, and mapping that to $\scrW(M^{in})$.

We now enter a little more into the details. There are several equivalent versions of Viterbo's definition, but we have found one to be particularly effective. Choose a Hamiltonian function $H$ whose vector field $X$ equals the Reeb vector field along $\partial M$. Then set
\begin{equation}
SH^*(M) = \underrightarrow{\text{\it lim}}_w\, HF^*(M;wH),
\end{equation}
where the right hand side is ordinary (Hamiltonian) Floer cohomology, and the direct limit, over all $w = 1,2,\dots$, is formed with respect to suitable continuation maps. In parallel, one can define $HW^*(L)$ to be the direct limit of Lagrangian Floer cohomologies $HF^*(L;wH)$. Recall that for each single $w$, the Floer cochain complex $CF^*(L;wH)$ is generated by chords
\begin{equation} \label{eq:intro-chords}
x: [0,1] \rightarrow M \text{ satisfying $dx/dt = wX$ and $x(0), x(1) \in L$.}
\end{equation}
If one chooses $H$ carefully, these generators will be in bijection with critical points of the Morse function $H|L$, together with Reeb chords of length $\leq w$ for the Legendrian submanifold $\partial L \subset \partial M$. Hence, wrapped Floer cohomology is a combination of ordinary Floer cohomology (which is isomorphic to the cohomology of $L$ in the present context) and linearized relative SFT (which counts Reeb chords, and sometimes goes under the name of linear Legendrian contact homology). Since direct limits are exact, one could in principle define the chain complex underlying $HW^*(L)$ to be the direct limit of the $CF^*(L;wH)$. However, a direct limit is in a sense a quotient construction, hence is unlikely to be strictly compatible with the $A_\infty$-structure. Instead, we use a larger quasi-isomorphic chain complex $CW^*(L;H)$, which is the homotopy direct limit (the algebraic version of the telescope construction in homotopy theory). By definition, the differential on this complex includes both the ordinary Floer differentials (for all $w$), and continuation maps from $w$ to $w+1$. When constructing the higher order components $\mu^d$ of the $A_{\infty}$-structure, one must similarly include terms arising from solutions to continuation map equations on $(d+1)$-punctured discs. In order to write down such an equation, one needs a one-form $\gamma$ on the disc, with suitable behaviour near the punctures, and which satisfies $d\gamma \leq 0$ everywhere. To organize the choices of one-forms, we introduce the notion of {\em popsicle}, which is a punctured disc equipped with additional geometric data, given by marked points which can move along special lines on the surface (intuitively, one imagines that $d\gamma$ is a smeared out version of the $\delta$-distribution located at those points). Of course, there is a certain amount of necessary preliminary work, which concerns the moduli spaces of popsicles and their compactifications.

The construction of $A_\infty$-homomorphisms underlying \eqref{eq:viterbo-open} has two main ingredients. The first idea, taken from \cite{viterbo97a}, is to compress the inner domain $M^{in} \subset M$ by using the Liouville flow. This shrinks all symplectic areas  by some factor $\rho$, and thereby also the values of the action functional on any chord \eqref{eq:intro-chords} lying inside $M^{in}$. Viterbo exploits this by making a careful choice of $\rho$ and of the Hamiltonian (both depending on $w$), which allows one to define the map $HF^*(M;wH) \rightarrow HF^*(M^{in};wH^{in})$ by projecting to a quotient of the Floer cochain complex. This can be carried out in the Lagrangian case as well, but it is not really suitable for our purpose since we need to consider all $w$ simultaneously. The second idea, used to overcome this difficulty, is adapted from work of Fukaya, Oh, Ohta and Ono \cite{fooo}*{Section 19}. In order to prove the independence of Fukaya-type $A_\infty$-structures from various choices, they considered parametrized moduli spaces arranged in trees, where the parameter is required to increase monotonically from the root to the leaves.  We adapt their construction to the setting of popsicles, and call the resulting objects \emph{cascades} (as far as the authors know, this is not related to the use of the same term in \cite{frauenfelder}). The parameter is then precisely Viterbo's rescaling factor $\rho$, which we allow to become arbitrarily small, meaning that it varies freely in $(0,1]$. We should point out that the resulting $A_\infty$-homomorphism has (in general) infinitely many higher order terms, hence is not a projection in any sense.

The structure of the paper is as follows. Sections \ref{sec:popsicle}--\ref{sec:many-objects} contain the main thread of the argument. The first three of these are devoted to, respectively, the discussion of popsicles and other Riemann surface matters; the definition of the $A_\infty$-algebra structure; and that of the $A_\infty$-homomorphism. Section \ref{sec:many-objects} contains some less central additional material, among it the extension of the previous constructions to wrapped Fukaya categories. The second part of the paper (Sections \ref{sec:lollipop}--\ref{sec:signs}) is a sequence of appendices, covering details which were omitted from the main exposition. These concern the construction of popsicle moduli spaces; a priori estimates that enter into various compactness arguments; transversality issues; and coherent orientations (with the resulting signs).

As should be clear from this description, the present paper is entirely foundational. One could certainly think of applications to the chord conjecture, in line with what has been done for symplectic cohomology, and also to Lagrangian embeddings, in particular fleshing out the last part of \cite{fukaya-seidel-smith07b}. However, for the authors of this paper, the main interest lies in potential applications to homological mirror symmetry, using Mikhalkin's generalized pair-of-pants decompositions \cite{mikhalkin04,abouzaid05}. This relies crucially on having the entire $A_\infty$-structure, rather than only the cohomology level theory $HW^*(L)$ (the same is true for \cite{fukaya-seidel-smith07b} as well). Unfortunately, explaining these ideas is beyond our scope here.

{\bf Acknowledgments.} We would like to thank Denis Auroux, Ivan Smith, and Katrin Wehrheim for valuable discussions. 

\section{Popsicles\label{sec:popsicle}}

We start this section by recalling one of the fundamental facts underlying the construction of Fukaya categories, namely the interpretation of Stasheff associahedra as moduli spaces of Riemann surfaces with boundary \cite{fukaya93}. Our main topic is an enriched version of this setup, where the Riemann surfaces carry holomorphic functions with suitable boundary conditions. We will explain the geometry of the resulting moduli spaces, relegating the more technical details to Section \ref{sec:lollipop}. As an important consequence, there is a natural class of one-forms on our Riemann surfaces, which will be used later to write down inhomogeneous $\bar\partial$-equations. As a final topic, we introduce yet another version of decorated Riemann surfaces, which come with additional real parameters (this time, the parameters have no geometric meaning, and adding them should be thought of as an abstract operadic-type construction).

\subsection{Discs\label{subsec:pointed}}
Fix $d \geq 0$. A {\em disc with $d+1$ ends} is a Riemann surface $S$ isomorphic to the closed unit disc $D$ minus $d+1$ boundary points (called points at infinity), together with a distinguished choice of one of those points (the distinguished point is called negative, and the others positive). We usually denote the points at infinity by $\bfzeta = \{\zeta^0,\dots,\zeta^d\}$. Here, the numbering convention is that $\zeta^0$ is the negative point, and that the positive ones appear in the order determined by the orientation of $\partial S$. The obvious compactification of $S$ is then written as $\bar{S} = S \cup \bfzeta \iso D$.

\begin{example}
A simple example is the infinite strip $Z = \R \times [0,1] \subset \C$, with points at infinity $\zeta^0 = -\infty$, $\zeta^1 = +\infty$. We usually write the coordinates on $Z$ as $z = s+it$. Of course, any other disc with two ends is isomorphic to $Z$.
\end{example}

The infinite strip is unstable, because it admits an infinite group of translational symmetries $(s,t) \mapsto (s+\sigma,t)$. Discs with $d+1 \geq 3$ ends are automatically stable, and in fact have no nontrivial automorphisms (preserving the distinguished point at infinity). In this stable range, it is easy to construct a universal family of discs, which we denote by $\scrS^{d+1} \rightarrow \scrR^{d+1}$. The base (the moduli space of stable discs) is a smooth manifold diffeomorphic to $\R^{d-2}$.

Consider a ribbon tree $T$ with $d+1$ semi-infinite edges, together with a preferred choice of one such edge (the preferred one is called the root, and the others the leaves). A {\em broken disc} modelled on $T$ is a
disjoint union
\begin{equation} \label{eq:broken}
S = \bigsqcup_v \; S_v
\end{equation}
indexed by vertices $v$ of $T$, where each component $S_v$ is a disc whose number of ends is given by the valency $|v|$ of the vertex. We assume from now on that stability holds, which means that $|v| \geq 3$ for all vertices of $T$. Let $\scrR^{d+1,T}$ be the moduli space of broken discs modelled on $T$. If the tree has a single vertex, this just agrees with $\scrR^{d+1}$, while in all other cases, it is a product of lower-dimensional moduli spaces $\scrR^{|v|}$. Take the disjoint union over all stable trees $T$ with $d+1$ semi-infinite edges:
\begin{equation} \label{eq:stasheff}
\bar\scrR^{d+1} = \bigsqcup_T \, \scrR^{d+1,T}.
\end{equation}
This can be equipped with a natural topology, and indeed the structure of a manifold with corners, which is such that the stratification \eqref{eq:stasheff} is the one into boundary and corner strata of various codimensions. In particular, the interior of $\bar\scrR^{d+1}$ is precisely $\scrR^{d+1}$. We call these compactifications the \emph{moduli spaces of stable broken discs} (they are homeomorphic to the classical Stasheff associahedra).

At this point, it is convenient to introduce some more combinatorial terminology. Take a tree $T$ as before. A flag (sometimes also called a half-edge) in $T$ is a pair consisting of a vertex and an adjacent edge. In our case, there is a preferred way of numbering the flags adjacent to any given vertex by $\{0,\dots,|v|\}$. Namely, one starts with the unique flag (called negative) which points towards the root, and continues with the other ones (called positive) according to the ribbon structure. Hence, if $S$ is a broken disc modelled on $T$, there is a natural correspondence between flags (adjacent to $v$) and points at infinity (belonging to the component $S_v$). Note also that the semi-infinite edges of $T$ can be numbered in a preferred way by $\{0,\dots,d\}$, by starting with the root and proceeding in accordance with the cyclic ordering given by any plane embedding $T \hookrightarrow \R^2$ compatible with the ribbon structure. Since each semi-infinite edge gives rise to a unique flag, we have $d+1$ distinguished points at infinity in $S$.

With this in mind, one defines the compactification $\bar{S}$ of a broken disc by first taking the compactifications $\bar{S}_v$ of all components, and then identifying those pairs of points at infinity which correspond to flags belonging to the same interior edge of $T$. Equivalently, one takes the disjoint union of the discs $S_v$, and adds to that $d+1$ points at infinity (corresponding to the semi-infinite edges of $T$) as well as singular points (corresponding to the finite edges). We call the added points of both type {\em special points} of $\bar{S}$. There is a natural compactification $\bar\scrS^{d+1}$ of the universal family $\scrS^{d+1}$, coming with a map
\begin{equation} \label{eq:compactified-family}
\bar\scrS^{d+1} \longrightarrow \bar\scrR^{d+1},
\end{equation}
such that the fibre over any point of $\bar\scrR^{d+1}$ is precisely the compactification of the broken disc corresponding to that point. Everywhere except at the singular points of compactified broken discs, $\bar\scrS^{d+1}$ is a manifold with corners. At those singular points, we have a different local model, namely
\begin{equation} \label{eq:singularity}
 \R^k \times [0,\infty)^l \times \{(\delta_+,\delta_-) \in \C^2 \;:\;
 \im(\delta_+) \geq 0,\im(\delta_-) \geq 0, \delta_+\delta_- \in [0,\infty) \}
\end{equation}
where of course, $k+l=d-3$ for dimension reasons. Still, the transition functions between these local models are differentiable (more precisely, restrictions of differentiable maps, in the coordinates given above). Using that, one can make sense of the notion of smooth map from a manifold into $\bar\scrS^{d+1}$, and also in the inverse direction. For instance, the projection map in \eqref{eq:compactified-family} is smooth; and the distinguished points at infinity give rise to a set of $d+1$ smooth, pairwise disjoint sections of that map. We will encounter the same situation again later on, see Section \ref{subsec:define-popsicle}, and return to the details in Section \ref{subsec:singularity}.
\begin{figure}[ht]
\begin{centering}
\epsfig{file=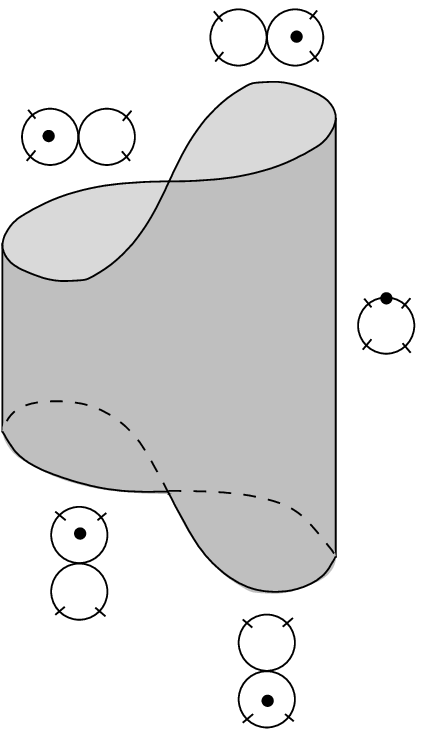}
\caption{\label{fig:vase-space}}
\end{centering}
\end{figure}
\begin{example}
Figure \ref{fig:vase-space} shows the three-dimensional space $\bar\scrS^4$. By definition, points of this space are pairs $(S,z)$, where $z \in \bar{S}$ is arbitrary. The picture shows what such pairs look like on each of the five boundary faces of $\bar\scrS^4$. In the case of the ``biggest'' face, $z$ can lie anywhere on $\partial \bar{S}$ (over the main stratum $\scrR^4 \iso \R$, where $\bar{S}$ is smooth, one gets a circle bundle; and as one gets to the limits, where the disc degenerates to a broken one, $\partial\bar{S}$ turns into a figure-eight).
%Note that there are no actual codimension three corners, only boundary faces, edges, and the two singular %points (this absence of top codimension corners is a general feature of $\bar\scrS^{d+1}$).
\end{example}

\subsection{Parametrizing the ends\label{subsec:strip-like-ends}}
Let $Z_- = \{s \leq 0\}, Z_+ = \{s \geq 0\} \subset Z$ be the half-infinite strips. Take a disc $S$ with $d+1$ ends. A {\em strip-like end} around a point at infinity $\zeta^k$ is a holomorphic embedding
\begin{equation}
\begin{cases} \epsilon: Z_- \longrightarrow S & \text{if } k = 0, \text{ or } \\
\epsilon: Z_+ \longrightarrow S & \text{if } k>0,
\end{cases}
\end{equation}
satisfying $\epsilon^{-1}(\partial S) = Z_{\pm} \cap \partial Z$ and $\lim_{s \rightarrow \pm \infty} \epsilon(s,\cdot) = \zeta^k$. As usual, the first kind of end is called negative, and the other positive. A set of strip-like ends for $S$ is a collection $\bfepsilon = \{\epsilon^0,\dots,\epsilon^d\}$, where each $\epsilon^k$ is an end around $\zeta^k$.

Strip-like ends are an auxiliary technical device, which enters into the standard gluing (or end connected sum) construction. To be specific, suppose that we have two discs $S_{\pm}$, a positive end $\epsilon_+: Z_+ \rightarrow S_+$, a negative end $\epsilon_-: Z_- \rightarrow S_-$, and a gluing parameter $\delta \in (-1,0)$. The associated gluing length is set to be $\sigma = -\log(-\delta)$. One then defines another disc
\begin{equation} \label{eq:gluing}
S = S_+ \# S_-
\end{equation}
by first removing $\epsilon_+([\sigma, +\infty) \times [0,1])$ as well as $\epsilon_-((-\infty,-\sigma] \times [0,1])$ from $S_+ \sqcup S_-$, and then identifying the remaining finite pieces of the ends through the formula $\epsilon_+(s,t) \sim \epsilon_-(s-\sigma,t)$. The distinguished negative point at infinity in $S$ is the one inherited from $S_+$.

This gluing process generalizes to families, and can be used to obtain collar neighbourhoods around the boundary strata of the compactified moduli space \eqref{eq:stasheff}. Fix some tree $T$ subject to the stability condition. Next, for each vertex $v$ of that tree, choose a set of strip-like ends for the universal family $\scrS^{|v|} \rightarrow \scrR^{|v|}$ (this means a set of strip-like ends on each fibre of that family, varying smoothly over the base). Every finite edge $e$ of the tree designates two vertices $v^\pm$ and the associated flags $(v^\pm,e)$, which in view of the considerations above, correspond to points at infinity in $\scrS^{|v_\pm|}$. One can then glue together the associated ends with some gluing parameter (say $\delta_e < 0$), or else leave them alone (which is the degenerate case $\delta_e = 0$). This gives rise to a gluing map on moduli spaces, of the form
\begin{equation} \label{eq:gluing-map}
\scrR^T \times \Big( \prod_e (-1,0] \Big) \supset \scrU^T \longrightarrow \bar\scrR^{d+1}.
\end{equation}
Here $\scrU^T$ is a subset where the gluing parameters are small, which means a neighbourhood
of $\scrR^T \times \{0,\dots,0\}$ (smallness of the parameter ensures that the various gluing processes don't interfere with each other). The main property of \eqref{eq:gluing-map} is that it is smooth and, after shrinking $\scrU^T$ if necessary, a diffeomorphism onto a neighbourhood of $\scrR^T \subset \bar\scrR^{d+1}$. This is proved for instance in \cite[Section 9]{seidel04}.

A universal choice of strip-like ends (for discs) is the choice of a set $\bfepsilon_S = \{\epsilon_S^k\}$ of ends for every stable disc $S$, which is compatible with isomorphisms of such discs, and smooth with respect to deformations. In practice, this means that one chooses ends for the universal family $\scrS^{d+1}$, varying smoothly over $\scrR^{d+1}$, and then obtains ends for any $S$ by identifying it with a fibre of the universal family. It is useful to impose an additional {\em consistency} requirement, which makes such a choice compatible with gluing processes. Namely, assume that a universal choice has been made, and take the situation from \eqref{eq:gluing}, where both $S_\pm$ are assumed stable. As part of our global choice, these two surfaces come with sets of strip-like ends; and if the gluing parameter is small, $S$ inherits a set of strip-like ends from that (basically, those ends which are not used up in the gluing process). The consistency condition says that for small values of the gluing parameter, these induced strip-like ends should be equivalent to those which are assigned to $S$ as part of the universal choice. The existence of consistent choices can be proved by an elementary inductive process, using coordinates \eqref{eq:gluing-map} near the boundary strata (see \cite[Section 9]{seidel04}; we will be using the same kind of construction repeatedly later on).

\begin{addendum} \label{th:rational-ends}
The space of all possible strip-like ends around any point at infinity is contractible (in an appropriate topology) but infinite-dimensional, and this can be occasionally awkward to work with. To cut down on the amount of choice, one can restrict attention to {\em rational} strip-like ends $Z_{\pm} \rightarrow S$, which are the ones that extend to an isomorphism $\bar{Z} \iso \bar{S}$. The rational strip-like ends around a given point at infinity form a principal homogeneous space over the contractible group $Aut(D,1)$ of holomorphic automorphisms fixing a point. Since the gluing process preserves rationality, consistency still makes sense in this more restrictive context. The same argument as before can then be used to produce consistent choices of rational strip-like ends.
\end{addendum}

\begin{addendum} \label{th:equivalent-ends}
The unstable case is problematic because strip-like ends cannot be made invariant under automorphisms of the Riemann surface. However, we can weaken the definition slightly, and then include at least the infinite strip. Call two strip-like ends {\em equivalent} if they differ by a translation of the domain $Z_\pm$. Equip $Z$ itself with the tautological ends, which are the inclusions $Z_\pm \hookrightarrow Z$. From that, any $S \iso Z$ inherits a set of strip-like ends which are unique up to equivalence. Moreover, if we take a given consistent universal choice and enlarge it to $d = 1$ in this way, the compatibility with the gluing process is maintained (up to equivalence, of course).
\end{addendum}

\subsection{The definition\label{subsec:define-popsicle}}
Fix $d \geq 1$ and in addition, labels $p_f \in \{1,\dots,d\}$ indexed by some finite set $F$. We write $\bfp = \{p_f\}$ for a collection of such labels. A {\em $\bfp$-flavoured popsicle} is a disc $S$ with $d+1$ ends, together with a collection $\bfphi = \{\phi_f\}$ of holomorphic maps
\begin{equation} \label{eq:popsicle-map}
\phi_f: S \longrightarrow Z,
\end{equation}
each of which extends to an isomorphism $\bar{S} \iso \bar{Z}$, taking $\zeta^0$ to $-\infty$ and $\zeta^{p_f}$ to $+\infty$. If we fix $S$, then the group $\R^F$ of translations acts simply transitively on the set of all possible popsicle structures $\bfphi$. Additionally, one can permute the $\phi_f$ within each subset where $p_f$ is constant, and that gives rise to an action of
\begin{equation} \label{eq:sym}
\Sym^{\bfp} = \{\text{permutations of $F$ stabilizing $\bfp$}\}.
\end{equation}

\begin{remark} \label{th:pre-sticks}
One way to represent the additional data geometrically is as follows. Take some $k \in \{1,\dots,d\}$, and any holomorphic map $\phi: S \rightarrow Z$ which extends to an isomorphism of compactifications, sending $\zeta^0$ to $-\infty$ and $\zeta^k$ to $+\infty$. Consider the infinite line in $S$ which is the preimage of $\{t = \half\} \subset Z$. This line is obviously the same for all maps, but a specific choice of $\phi$ singles out a point on it, namely the preimage of $(0,\half) \in Z$. In this way, a $\bfp$-flavoured popsicle structure can be encoded in a collection of points on $S$ indexed by $F$ (some of the points may coincide, and of course their position is not free, each being constrained to the infinite line determined by $k = p_f$). The group \eqref{eq:sym} acts by permuting those points.
\end{remark}

Popsicles are stable, and in fact have no nontrivial automorphisms (compatible with the given structure), if and only if $d+|F| \geq 2$. In that range we have a universal family of popsicles, denoted by $\scrS^{d+1,\bfp} \rightarrow \scrR^{d+1,\bfp}$, whose base has dimension $d-2+|F|$. The structure of these moduli space is actually quite easy to understand. Assume first that $d \geq 2$. In that case, each fibre of the forgetful map
\begin{equation} \label{eq:forgetful}
\scrR^{d+1,\bfp} \longrightarrow \scrR^{d+1}
\end{equation}
can be identified with $\R^F$. This identification depends on a choice of origin in the fibre, but that can be chosen to vary smoothly with $\scrR^{d+1}$, giving rise to a non-canonical diffeomorphism
$\scrR^{d+1,\bfp} \iso \scrR^{d+1} \times \R^F \iso \R^{d-2+F}$. The action of $\R^F$ is just translation on the fibres of \eqref{eq:forgetful}, and (if one picks the origin to be invariant, which is always possible) $\Sym^{\bfp}$ acts by permuting coordinates. In the remaining case $d = 1$, one can identify $\scrR^{d+1,\bfp}$ with the quotient $\R^F/\R$ (where $\R$ acts diagonally by translation on all factors). Note that in either case, any transposition in $\Sym^{\bfp}$ acts orientation-reversingly.

\begin{addendum} \label{th:isotropy}
Suppose that we have a partition $\bfP$ of $F$ into subsets $\{P_m\}$, such that the map $f \mapsto p_f$ is constant on each subset. The permutations of $F$ which preserve $\bfP$ form a subgroup $\Sym^{\bfP} \subset \Sym^{\bfp}$. These are precisely the subgroups appearing as isotropy groups for the action of \eqref{eq:sym}. The corresponding fixed point set in $\scrR^{d+1,\bfp}$ is the set of those $(S,\bfphi)$ such that for any two $f$ belonging to the same $P_m$, the associated maps $\phi_f$ coincide.
\end{addendum}

Take a tree $T$ of the same kind as before, and a decomposition $\bfF = \{F_v\}$ of $F$ into subsets labeled by vertices of that tree. Moreover, the decomposition needs to be compatible with $\bfp$ in the following sense: for each $f \in F_v$, the vertex $v$ must lie on the path from the root to the $p_f$-th leaf. Note that $\bfp$ and the decomposition $\bfF$ determine numbers $\bfp_v = \{p_{v,f}\}$ for each vertex, where $f$ runs through $F_v$ and $p_{v,f} \in \{1,\dots,|v|\}$. This is because the path mentioned above leaves $v$ in the direction of a specific flag, whose number is taken to be $p_{v,f}$. For instance, if $T$ is the tree with just one vertex $v$, the compatibility condition implies that we necessarily have $\bfp_v = \bfp$. Then, a {\em broken popsicle} is a disjoint union \eqref{eq:broken}, together with a $\bfp_v$-flavoured popsicle structure $\bfphi_v$ on each $S_v$. This time, the stability condition is that two-valent vertices $v$ are allowed as long as $F_v \neq \emptyset$. The moduli space of stable broken popsicles modelled on $(T,\bfF)$ is a product
\begin{equation} \label{eq:product}
\scrR^{T,\bfF} \iso \prod_v \scrR^{|v|,\bfp_v}.
\end{equation}
Take the disjoint union of those spaces,
\begin{equation} \label{eq:compactification}
\bar\scrR^{d+1,\bfp} = \bigsqcup_{(T,\bfF)} \scrR^{T,\bfF}.
\end{equation}
We will see later (Corollary \ref{th:real-and-negative}) that $\bar\scrR^{d+1,\bfp}$ carries a canonical structure of a smooth manifold with corners, such that \eqref{eq:compactification} is the associated standard stratification, generalizing the corresponding statement for \eqref{eq:stasheff} (which would be the special case $F = \emptyset$). The forgetful map \eqref{eq:forgetful} (defined if $d \geq 2$), as well as the actions of $\R^F$ and $\Sym^{\bfp}$, extend smoothly to the compactification. However, their geometry becomes a little more complicated: the compactified forgetful map is no longer a submersion, and the action of $\R^F$ on the compactification is not free for any value of $d$ (both phenomena are due to the appearance of unstable components $S_v \iso Z$ in the boundary strata). Finally, we should point out that on the closure of a stratum \eqref{eq:product}, the action of $\Sym^{\bfp}$ embeds into one of a potentially larger group, namely $\prod_v \Sym^{{\bfp}_v}$.
\begin{figure}[hb]
\begin{centering}
\epsfig{file=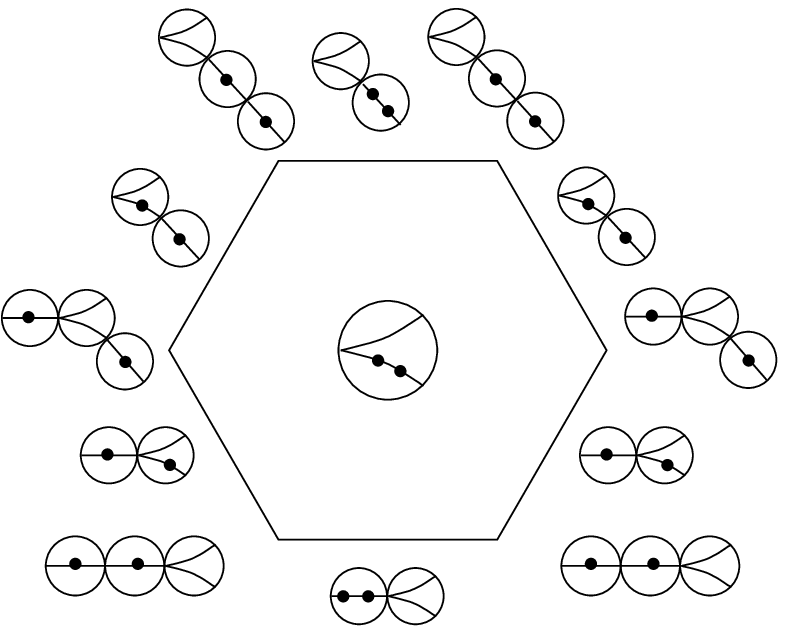}
\caption{\label{fig:hexagon-space}}
\end{centering}
\end{figure}
\begin{examples} (i) \label{th:ex-popsicle}
Take $d = 2$, $F = \{1,2\}$, and $p_1 = p_2 = 1$. Then the compactified moduli space is a hexagon (Figure \ref{fig:hexagon-space}). The group of symmetries $\Sym^{\bfp} = \Z/2$ acts by reflection along the vertical axis.

\begin{figure}[hb]
\begin{centering}
\epsfig{file=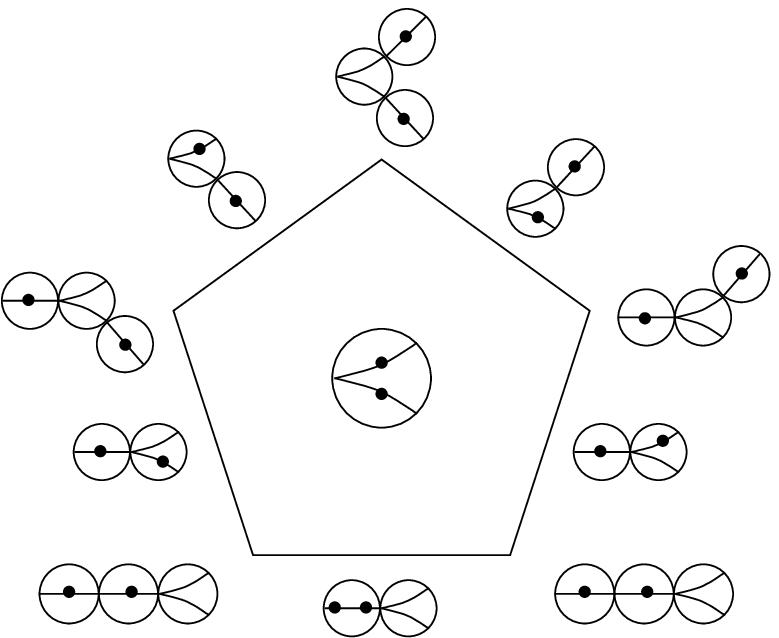}
\caption{\label{fig:pentagon-space}}
\end{centering}
\end{figure}

(ii) Take the same $d$ and $F$, but now with $p_1 = 1$, $p_2 = 2$. Then the compactified moduli space is a pentagon (Figure \ref{fig:pentagon-space}). Note that while $\Sym^{\bfp}$ is trivial in this case, one of the boundary edges (the bottom one in our picture) carries a nontrivial symmetry group $\prod_v \Sym^{{\bfp}_v} = \Z/2$. In particular, the two corners adjacent to this edge are isomorphic moduli spaces. By this, we mean that the product expressions \eqref{eq:product} for these two corner strata consist of the same factors, even though the strata themselves are distinguished by how the given $F$ is partitioned into $F_v$.

It is also instructive to consider the $\R^2$-action on this moduli space. The orbits of the two vector fields generating the action are (schematically) drawn in Figure \ref{fig:translate-pentagon}. Note that where orbits appear to be transverse to the boundary, their speed actually slows down to zero, and the limiting boundary edge consists of stationary points.
\begin{figure}[h]
\begin{centering}
\epsfig{file=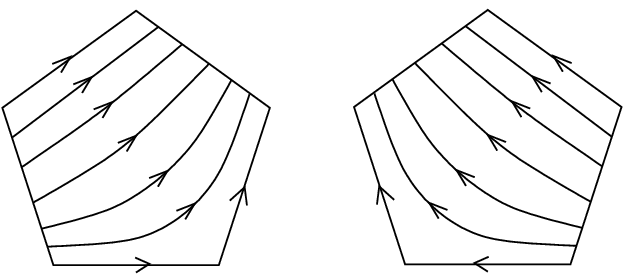}
\caption{\label{fig:translate-pentagon}}
\end{centering}
\end{figure}
\end{examples}

In exactly the same way as in \eqref{eq:compactified-family}, one can define compactifications of the universal families of popsicles,
\begin{equation} \label{eq:compactified-total-space}
\bar\scrS^{d+1,\bfp} \longrightarrow \bar\scrR^{d+1,\bfp}.
\end{equation}
The total space has the same geometry as before, meaning that the local models are either the standard corners or \eqref{eq:singularity}. In addition to the $d+1$ canonical sections, the compactified family now comes with a collections of smooth maps to $\bar{Z}$. If we consider the fibre over a boundary point, which is the compactification of a broken popsicle, then each of these maps $\bar{\phi}_f: \bar{S} \rightarrow \bar{Z}$ is constant and equal to $\pm\infty$ on all components except one (the non-constant component $\bar{S}_v$ being precisely the one such that $f \in F_v$).

\subsection{Popsicle sticks\label{subsec:sticks}}
The gluing process for popsicles is a little less straightforward than that for pointed discs. It seems difficult to work directly on the level of maps \eqref{eq:popsicle-map}. We will therefore adapt the idea from Remark \ref{th:pre-sticks}, replacing these maps by additional marked points lying on certain infinite lines. However, instead of taking those lines to be preimages of the line $\{t = 1/2\} \subset Z$, we will bend them slightly to achieve better compatibility with the strip-like ends, and hence with the gluing process.

Concretely, let $S$ be a disc with $d+1$ ends, which comes equipped with a set $\bfepsilon$ of strip-like ends. Fix some $k \in \{1,\dots,d\}$. A {\em popsicle stick} connecting $\zeta^0$ to $\zeta^k$ is an infinite line $\R \iso Q \subset S$ satisfying
\begin{equation} \label{eq:stick}
\begin{cases}
\!\!\! & \phi(Q) = \{t = \tau(s)\} \text{ for some } \tau: \R \longrightarrow (0,1), \\
\!\!\! & Q = \epsilon^0(\{t = \half\}) \text{ in a neighbourhood of $\zeta^0$}, \\
\!\!\! & Q = \epsilon^k(\{t = \half\}) \text{ in a neighbourhood of $\zeta^k$}.
\end{cases}
\end{equation}
Here, $\phi: S \rightarrow Z$ is any map as in Remark \ref{th:pre-sticks}. To see that the three conditions in \eqref{eq:stick} are not mutually contradictory, one looks at the compositions $\rho_- = \phi \circ \epsilon^0: Z_- \rightarrow Z$ and $\rho_+ = \phi \circ \epsilon^k: Z_+ \rightarrow Z$. Either of these can be written as the sum of a translation and a remainder term which decays exponentially as $s \rightarrow \pm\infty$. Therefore, the image of $\{t = 1/2, \pm s \gg 0\}$ under $\rho_\pm$ is of the form $\{t = \tau_\pm(s)\}$. To build a popsicle stick, it is then sufficient to choose a function $\tau$ which at infinity agrees with the given $\tau_\pm$. The particular choice of $\phi$ is irrelevant for the definition of popsicle stick (and, incidentally, only the equivalence class of the strip-like ends matters). However, given $Q$, any choice of $\phi$ singles out a point $q \in Q$, namely the one which satisfies
\begin{equation} \label{eq:marker}
 \phi(q) \in \{0\} \times (0,1),
\end{equation}
or equivalently in terms of \eqref{eq:stick}, $\phi(q) = (0,\tau(0))$. This sets up a bijection between choices of $\phi$ and points $q$. We call those points {\em sprinkles}.

Define a set of popsicle sticks for $S$ to be a collection $\bfQ = \{Q^1,\dots,Q^d\}$ where each $Q^k$ connects $\zeta^0$ to $\zeta^k$. Clearly, given a set of popsicle sticks, there is a bijection between $\bfp$-flavoured popsicle structures and collections $\bfq = \{q_f\}$ of sprinkles, such that $q_f \in Q^{p_f}$ for all $f \in F$.

With this at hand, we return to the gluing construction. Take discs $S_\pm$, which come with sets of strip-like ends and popsicle sticks. For small values of the gluing parameter, the surface $S$ defined as in \eqref{eq:gluing} inherits the same additional structures. Moreover, if $S_\pm$ have $\bfp_\pm$-flavoured popsicle structures, then $S$ inherits a $\bfp$-flavoured one, provided that $\bfp$ is related to $\bfp_\pm$ in a suitable way (which is straightforward to infer from our prescription for gluing popsicle sticks). As before, this process can be applied to families, and in particular to the boundary strata in \eqref{eq:compactification}. Pick a pair $(T,\bfF)$ describing such a stratum. Assume that for each vertex $v$, we have chosen strip-like ends and popsicle sticks for the family $\scrS^{|v|,\bfp_v} \rightarrow \scrR^{|v|,\bfp_v}$. The resulting gluing map has the form
\begin{equation} \label{eq:gluing-with-sprinkles}
\scrR^{T,\bfF} \times \Big( \prod_e (-1,0] \Big) \supset \scrU^{T,\bfF} \longrightarrow \bar\scrR^{d+1,p}.
\end{equation}
As in the previous case \eqref{eq:gluing}, it will turn out that this is smooth, and a diffeomorphism onto a neighbourhood of $\scrR^{T,\bfF}$. Strictly speaking, we will prove this only under the additional assumption that the strip-like ends used are rational ones (Corollary \ref{th:glue-well}). However, the general case is only slightly more involved, and interested readers (if any) should be easily able to extend the given argument.

Still proceeding in parallel with the previous discussion, we now introduce the relevant notions of universal choices. Suppose first that we have made a universal choice of strip-like ends for pointed discs, which should be consistent. Then, a universal choice of popsicle sticks is a choice of a set $\bfQ_S$ for every stable pointed disc $S$, varying smoothly over the moduli space. There is an appropriate notion of consistency, which as before means that the popsicle sticks are compatible with the gluing construction, for small values of the gluing parameter. Consistent choices can be constructed inductively. Moreover, there is an obvious way to include the case $d = 1$ as well: for $S \iso Z$, take the popsicle stick to be $\{t = \half\}$. The final step is to pick, for every stable popsicle $(S,\bfphi)$ of some flavour $\bfp$, a set of strip-like ends
$\bfepsilon_{S,\bfphi}$, which are equivalent to those previously chosen for the underlying pointed disc $S$ if $d \geq 2$, and which in the remaining case $d = 1$ lie in the equivalence class from Addendum \ref{th:equivalent-ends}. As usual, the new ends should vary smoothly over the relevant moduli spaces. Additionally, they should be invariant under the action of \eqref{eq:sym}. Given choices of strip-like ends and popsicle sticks, one then has a gluing construction for popsicle structures, and can impose the relevant consistency condition on the $\bfepsilon_{S,\bfphi}$. As before, the structure of the maps \eqref{eq:gluing-with-sprinkles} is crucial in proving that this condition can be met. We omit the details of this, since the argument is essentially the same as that for pointed discs.

\begin{remark}
The choice of $\bfQ_S$ itself only involves the underlying surface $S$. In other words, it is pulled back through the forgetful map \eqref{eq:forgetful}. In principle, it would be possible to lift that restriction, making the popsicle sticks depend on $\bfphi$ as well, but there seems little point in doing so. The situation for strip-like ends is slightly different. While the equivalence class of $\bfepsilon_{S,\bfphi}$ depends only on $S$, the specific choice of end varies with the popsicle structure. This is necessary in order to achieve consistency, because there are stable popsicles $(S,\bfphi)$ with unstable underlying disc $S \iso Z$ (and moreover, even if one restricts to $d \geq 2$, such popsicles occur as components in the compactification of $\scrR^{d+1,\bfp}$).
\end{remark}

We want to fix these auxiliary structures once and for all. To reiterate, this involves: first, a universal choice of strip-like ends for pointed discs (supposed to be consistent, and for technical reasons, also rational); second, a consistent universal choice of popsicle sticks; and third, a consistent universal choice of strip-like ends for popsicles. All these choices will be kept constant throughout the rest of the paper.

\subsection{Closed one-forms\label{subsec:closed-forms}}
We now introduce a class of one-forms on discs, which can be seen as generalizations of $dt \in \Omega^1(Z)$. Namely, let $S$ be a disc equipped with a set of strip-like ends. By a set of {\em basic closed one-forms} for $S$ we mean a collection $\bfalpha = \{\alpha^1,\dots,\alpha^d\}$, where each $\alpha^k \in \Omega^1(S)$ is closed, with $\alpha^k|\partial S = 0$, and satisfies
\begin{equation} \label{eq:basic-closed}
(\epsilon^j)^*\alpha^k = \begin{cases} dt & \text{if $j = 0$ or $k$}, \\
0 & \text{otherwise}
\end{cases}
\end{equation}
over the ends (more precisely, on a subset of $Z_\pm$ where $\pm s \gg 0$).

Basic closed one-forms behave well with respect to the gluing process \eqref{eq:gluing}. To spell this out, let $S_\pm$ be pointed discs with $d_\pm + 1$ ends, respectively, equipped with strip-like ends and basic closed one-forms. Suppose that we glue together the ends corresponding to the points at infinity $\zeta^0_-$ and $\zeta^i_+$, for some $i \in \{1,\dots,d_+\}$, using a small gluing parameter. The result is a disc $S$ with $d+1 = d_- + d_+$ ends. This will come with induced basic closed one-forms, which are constructed as follows:
\begin{equation}
\begin{cases}
& \!\! \text{in the case $k<i$, take $\alpha^k_+$ and extend it trivially (by zero);} \\
& \!\! \text{for $k \geq i+d_-$ do the same, but starting with $\alpha^{k-d_- + 1}_+$;} \\
& \!\! \text{for $i \leq k < i+d_-$, take $\alpha^i_+$ and glue it together with $\alpha^{k-i+1}_-$.}
\end{cases}
\end{equation}
For a fixed $S$, it is clear that basic closed one-forms exist, and in fact form a contractible space. In the stable range $d+1 \geq 3$, we can choose such one-forms $\bfalpha_S = \{\alpha_S^k\}$ for each pointed disc $S$, varying smoothly over the moduli space (and, of course, compatibly with our fixed universal choice of strip-like ends). Moreover, this can be done in such a way that the obvious consistency condition, for all gluing processes with small gluing parameter, is satisfied. Finally, returning to the original motivation, we extend this to $d = 1$ by setting $\alpha^1_Z = dt$.
%(this is the unique possible translation-invariant choice, which makes sense if one thinks of ``$\scrR^2$'' as %the stack $point/\R$).

\subsection{Sub-closed one-forms\label{subsec:sub-closed-forms}}
Let $(S,\bfphi)$ be a $\bfp$-flavoured popsicle, equipped with strip-like ends. A set $\bfbeta = \{\beta_f\}$ of {\em basic sub-closed one-forms} consists of a one-form $\beta_f$ for each $f \in F$, such that $\beta_f|\partial S = 0$, $d\beta_f = 0$ in a neighbourhood of $\partial S$, and $d\beta_f \leq 0$ everywhere (this is sub-closedness; the sign is with respect to the complex orientation of $S$). The analogue of \eqref{eq:basic-closed} is the requirement that at infinity,
\begin{equation} \label{eq:basic-sub-closed}
(\epsilon^j)^*\beta_f = \begin{cases} dt & \text{if $j = 0$}, \\
0 & \text{otherwise.}
\end{cases}
\end{equation}
This may seem puzzling at first sight, because $f$ does not appear on the right hand side of \eqref{eq:basic-sub-closed}; however, looking at the gluing process provides some clarification. Take two popsicles $S_\pm$ with flavours $\bfp_\pm$. These should be equipped with strip-like ends, popsicle sticks, and basic closed as well as sub-closed one-forms. Glue $S_\pm$ together as before, using the process from
Section \ref{subsec:sticks} to produce a popsicle structure on the resulting disc $S$. Recall that the flavour $\bfp$ of the new structure is indexed by the disjoint union $F = F_+ \sqcup F_-$. We then define sub-closed one-forms $\beta_f$ on the surface obtained by gluing $S_+$ to $S_-$ at the $i$ incoming marked point of $S_+$ as follows:
\begin{equation}
\begin{cases}
& \!\! \text{if $f \in F_+$, take the given $\beta_{+,f}$ and extend it trivially;} \\
& \!\! \text{if $f \in F_-$, take $\beta_{-,f}$ and glue it together with $\alpha_+^i$.}
\end{cases}
\end{equation}

Take a universal choice of sub-closed one-forms $\{\bfbeta_{S,\bfphi}\}$. We will always assume that this choice is invariant under the action of $\Sym^\bfp$ (which is unproblematic, since sub-closed one-forms satisfying \eqref{eq:basic-sub-closed} form a convex subspace, where averages can be formed). Moreover, we can impose a consistency condition with respect to the gluing process introduced above (and our previous choice of strip-like ends as well as popsicle sticks); choices satisfying this condition can be constructed by applying the usual process, based on the recursive structure of the compactification \eqref{eq:compactification}.

\begin{example}
For ease of visualization, it's better to think of the two-forms $d\beta_f$, rather than the $\beta_f$ themselves. For instance, consistency (compatibility with gluing) has a straightforward meaning, because each $d\beta_f$ is compactly supported. Let's take a concrete look at the simplest situations. The trivial situation is when $d = 1$ and $|F| = 1$. The (unique) point in that moduli space is represented by a surface $S \iso Z$, and we can fix the isomorphism in such a way that the unique sprinkle gets mapped to $(0,1/2)$. We then have a single two-form $d\beta \leq 0$ supported on some compact subset of $int(Z)$. In principle, this can be chosen arbitrarily, subject only to the condition that $\int_Z d\beta = -1$, but it is convenient to imagine its support lying in a small neighbourhood of $(0,1/2)$. Now take $d = 1$ and $F = \{1,2\}$. A $\bfp$-flavoured popsicle is a surface $S \iso Z$ equipped with two sprinkles, whose images in $Z$ are $(s_1,1/2)$ and $(s_2,1/2)$. The modular parameter is $\sigma = s_1-s_2$, giving an identification $\scrR^{2,\bfp} \iso \R$. On the subset where $|s_1-s_2| \gg 0$, one starts with $d\beta_1$ and $d\beta_2$ both being copies of the previously introduced $d\beta$, translated away from each other exactly by the amount $\sigma$, so that each $d\beta_f$ is supported near the corresponding sprinkle. Over the rest of the moduli space, one interpolates between those given choices, while preserving the symmetry which exchanges both the sprinkles and the $d\beta_f$.
\end{example}

As before we end our discussion by fixing, once and for all, a consistent universal choice of closed basic forms $\bfalpha_S$ and of sub-closed basic forms $\bfbeta_{S,\bfphi}$.

\subsection{Weighted popsicles\label{subsec:weights}}
Take some $d \geq 1$ and $\bfp = \{p_f\}$, $f \in F$. A set of {\em weights} is a collection of positive integers $\bfw = \{w^0,\dots,w^d\}$ satisfying
\begin{equation} \label{eq:weights}
w^0 = w^1 + \cdots + w^d + |F|.
\end{equation}
The primary use of this concept is as follows. Suppose that $(S,\bfphi)$ is a $\bfp$-flavoured popsicle, which also comes with strip-like ends $\bfepsilon$, basic closed one-forms $\bfalpha$, and sub-closed one-forms $\bfbeta$. Given weights, one then defines the {\em total sub-closed one-form} $\gamma \in \Omega^1(S)$ to be
\begin{equation} \label{eq:total-one-form}
\gamma = \sum_{k=1}^d w^k\alpha^k + \sum_{f \in F} \beta_f.
\end{equation}
By construction, see \eqref{eq:basic-closed}, \eqref{eq:basic-sub-closed} and \eqref{eq:weights}, this has the following properties:
\begin{equation} \label{eq:gamma-properties}
\parbox{30em}{
$\gamma|\partial S = 0$. Moreover, $d\gamma$ is non-positive everywhere, and vanishes in a neighbourhood of $\partial S$. Finally, on the strip-like ends we have $(\epsilon^k)^*\gamma = w^k dt$, at all points $(s,t)$ with $\pm s \gg 0$.
}
\end{equation}
Recall that in Section \ref{subsec:sub-closed-forms} we made a universal choice of basic closed and sub-closed one-forms. This determines, for each stable popsicle with weights, a total sub-closed one-form \eqref{eq:total-one-form}, which we denote by $\gamma_{S,\bfphi,\bfw}$. Note that by construction, our choice of total sub-closed one-form is compatible with the action of \eqref{eq:sym} on the moduli space. Concretely, this means that if we have two weighted popsicles which differ only by a reordering of the maps $\phi_f$, then they carry the same total sub-closed one-form. Finally, in the unstable case where $S \iso Z$ and $F = \emptyset$, hence $\bfw = \{w^0 = w^1 = w\}$, we set
\begin{equation} \label{eq:unstable-gamma}
\gamma_{Z,\bfw} = w\, dt.
\end{equation}
Because of the previously imposed consistency conditions on the one-forms $\alpha$ and $\beta$, their sums $\gamma$ are compatible with gluing processes (assuming that any two ends which are glued together carry the same weight; and that the gluing parameters remain small).

We write $\scrR^{d+1,\bfp,\bfw}$ for the moduli space of stable popsicles with weight $\bfw$. Of course, this is just a copy of $\scrR^{d+1,\bfp}$, but the separate notation is still useful. For instance, if we take the compactification $\bar\scrR^{d+1,\bfp,\bfw}$ defined as in \eqref{eq:compactification}, then the boundary strata are canonically products of lower-dimensional moduli spaces of weighted popsicles. To see this explicitly, fix a set of weights $\bfw$, and consider a pair $(T,\bfF)$ which labels one of the boundary strata. Recall from Section \ref{subsec:pointed} that pairs $(v,k)$, where $v$ is a vertex and $k$ is a number in $\{0,\dots,|v|-1\}$, correspond canonically to flags of $T$. We can then define induced weights $\bfw_v = \{w_v^k\}$ for every vertex $v$, by the following conditions:  first, if $(v,k)$ corresponds to the $j$-th semi-infinite edge of $T$, then $w_v^k = w^j$; secondly, two flags belonging to the same finite edge should carry the same weight; and finally, the analogue of \eqref{eq:weights} should be satisfied at each vertex, meaning that $w_v^0 = w_v^1 + \cdots + w_v^{|v|-1} + |F_v|$. In these terms, the boundary stratum of $\bar\scrR^{d+1,\bfp,\bfw}$ indexed by $(T,\bfF)$ is isomorphic to the product of the spaces $\bar\scrR^{|v|,\bfp_v,\bfw_v}$.

\subsection{Cascades\label{subsec:cascade}}
What follows now is a bad case of mixed metaphor. Fix $(d,\bfp)$. Take a pair $(T,\bfF)$ as in the definition of broken popsicle, and consider the induced flavours $\{\bfp_v\}$. A {\em cascade} modelled on $(T,\bfF)$ is a collection $\{(\rho_v,S_v,\bfphi_v)\}$, where each $(S_v,\bfphi_v)$ is a $\bfp_v$-flavoured popsicle, equipped with an additional parameter value $\rho_v \in (0,1]$. These parameters are subject to the following {\em causal ordering} condition:
\begin{equation} \label{eq:time-ordering}
\parbox{30em}{
If a finite edge of $T$ has endpoints $v_\pm$, with $v_+$ being closer to the root than $v_-$, then $\rho_{v_+} < \rho_{v_-}$.}
\end{equation}
We can impose a stability condition, which is the same as for broken popsicles, and take the disjoint union of the resulting moduli spaces for all stable choices $(T,\bfF)$. We call the resulting space \emph{moduli space of stable cascades} and denote it by $\scrQ^{d+1,\bfp}$. Note that even though the notation is analogous to \eqref{eq:compactification}, the structure of the moduli space is entirely different, since each $(T,\bfF)$ defines a connected component of the same dimension $d-1+|F|$.

There is also a natural partial compactification, which is a manifold with corners denoted by $\bar\scrQ^{d+1,\bfp}$. The strata of this space are indexed by combinatorial data
\begin{equation} \label{eq:5-label}
(T,\bfF,\bar{E},E)
\end{equation}
where $(T,\bfF)$ is as before, and $\bar{E} \subset E$ are subsets of the set of finite edges of $T$. A point of such a stratum is again given by a collection $\{(\rho_v, S_v, \bfphi_v)\}$, with the difference that we now require the parameter values to satisfy \eqref{eq:time-ordering} only if the edge connecting $v_\pm$ is not in $E$, whereas for the other edges the equality $\rho_{v_+} = \rho_{v_-}$ should hold. The codimension of the stratum is then precisely $|E|$. To understand the topology of the partial compactification, note that one can associate to each datum \eqref{eq:5-label} a collapsed tree $\bar{T}$, obtained by contracting the edges in $\bar{E}$ to points. For any vertex $\bar{v}$ of $\bar{T}$, let $\bar{F}_{\bar{v}}$ be the union of $F_v$ over all vertices $v$ mapped to $\bar{v}$ under the collapsing map. Denote the resulting decomposition of $F$ by $\bar\bfF$. Reversing the direction of the argument, take a fixed $(\bar{T},\bar{\bfF})$, with its associated collection of labels $\{\bar\bfp_{\bar{v}}\}$, and consider all choices of \eqref{eq:5-label} which yield the given pair after collapsing. The union of all the associated strata in $\bar\scrQ^{d+1,\bfp}$ is a copy of
\begin{equation} \label{eq:delta-simplex}
 \Delta^{\bar{T}} \times \prod_{\bar{v}} \bar\scrR^{|\bar{v}|,\bar\bfp_{\bar{v}}},
\end{equation}
where $\Delta^{\bar{T}}$ is the set of all parameter values $\{\bar\rho_{\bar{v}}\}$ indexed by vertices, such that in the situation analogous to \eqref{eq:time-ordering}, the weak inequality $\bar\rho_{\bar{v}_+} \leq \bar\rho_{\bar{v}_-}$ holds. Then \eqref{eq:delta-simplex}, with its natural topology, is a connected component of $\bar\scrQ^{d+1,\bfp}$. As a direct consequence, the map
\begin{equation} \label{eq:properity}
\begin{aligned}
& \bar\scrQ^{d+1,\bfp} \longrightarrow (0,1], \\
& \{(\rho_v,S_v,\bfphi_v)\} \longmapsto \rho_{v_0},
\end{aligned}
\end{equation}
where $v_0$ is the vertex closest to the root, is proper. This means that the non-compactness of $\bar\scrQ^{d+1,\bfp}$ only reflects that of the parameter interval $(0,1]$.

It is important to note that a given collection $\{(\rho_v,S_v,\bfphi_v)\}$ can appear more than once in $\bar\scrQ^{d+1,\bfp}$. The reason is that, while $E$ is determined by the relation between parameters $\rho_v$, the choice of $\bar{E} \subset E$ is free. In particular, codimension one faces defined by the equality of exactly two parameters occur in isomorphic pairs (the two faces in such a pair belong to different connected components of $\bar\scrQ^{d+1,\bfp}$).

\begin{example} \label{th:two-boundary}

In the simplest example, $T$ is a tree with two vertices $v_\pm$, and we have corresponding triples $\{(\rho_\pm,S_\pm,\bfphi_\pm)\}$ with $\rho_+ = \rho_- = \rho$, so that $E$ consists of the unique finite edge of the tree. If one then sets $\bar{E} = \emptyset$, the resulting point in $\bar\scrQ^{d+1,\bfp}$ is the limit of a sequence $\{(\rho_{l,\pm},S_{l,\pm},\bfphi_{l,\pm})\}$ with the same combinatorial data $(\bar{T},\bar{\bfF}) = (T,\bfF)$, and where $\rho_{l,\pm}$ become equal in the limit $l \rightarrow \infty$. The other possibility is $\bar{E} = E$. In that case $\bar{T}$ has a single vertex, so the resulting point in $\bar\scrQ^{d+1,\bfp}$ is the limit of a sequence $(\rho_l,S_l,\bfphi_l)$ in which $\rho_l \rightarrow \rho$, and the popsicle $(S_l,\bfphi_l)$ degenerates into the broken popsicle given by $\{(S_\pm,\bfphi_\pm)\}$
(see Figure \ref{fig:double} for a schematic picture).
\end{example}
\begin{figure}[ht]
\begin{centering}
\begin{picture}(0,0)%
\includegraphics{double.pstex}%
\end{picture}%
\setlength{\unitlength}{3355sp}%
\begingroup\makeatletter\ifx\SetFigFont\undefined%
\gdef\SetFigFont#1#2#3#4#5{%
  \reset@font\fontsize{#1}{#2pt}%
  \fontfamily{#3}\fontseries{#4}\fontshape{#5}%
  \selectfont}%
\fi\endgroup%
\begin{picture}(5034,3414)(-72,-1444)
%\put(2746,-436){\makebox(0,0)[lb]{\smash{{\SetFigFont{10}{12.0}{\rmdefault}{\mddefault}{\updefault}{\color[rgb]{0,0,0}$\rho$}%
%}}}}
%\put(2026,-1036){\makebox(0,0)[lb]{\smash{{\SetFigFont{10}{12.0}{\rmdefault}{\mddefault}{\updefault}{\color[rgb]{0,0,0}$\rho$}%
%}}}}
\put(296,614){\makebox(0,0)[lb]{\smash{{\SetFigFont{10}{12.0}{\rmdefault}{\mddefault}{\updefault}{\color[rgb]{0,0,0}$\rho_{l,+}$}%
}}}}
\put(1066,1214){\makebox(0,0)[lb]{\smash{{\SetFigFont{10}{12.0}{\rmdefault}{\mddefault}{\updefault}{\color[rgb]{0,0,0}$\rho_{l,-}$}%
}}}}
\put(3976,1814){\makebox(0,0)[lb]{\smash{{\SetFigFont{10}{12.0}{\rmdefault}{\mddefault}{\updefault}{\color[rgb]{0,0,0}$\zeta^3_l$}%
}}}}
\put(4501,1664){\makebox(0,0)[lb]{\smash{{\SetFigFont{10}{12.0}{\rmdefault}{\mddefault}{\updefault}{\color[rgb]{0,0,0}$\zeta^2_l$}%
}}}}
\put(4051,-136){\makebox(0,0)[lb]{\smash{{\SetFigFont{10}{12.0}{\rmdefault}{\mddefault}{\updefault}{\color[rgb]{0,0,0}$\zeta^2_l, \zeta^3_l \rightarrow \zeta^2_+$}%
}}}}
\put(-340,-286){\makebox(0,0)[lb]{\smash{{\SetFigFont{10}{12.0}{\rmdefault}{\mddefault}{\updefault}{\color[rgb]{0,0,0}$ \rho_{l,+}, \rho_{l,-} \rightarrow \rho$}%
}}}}
\end{picture}%
\caption{\label{fig:double}}
\end{centering}
\end{figure}%

Finally, $\scrQ^{d+1,\bfp}$ carries an action of $\Sym^{\bfp}$. If we fix a connected component, and restrict to the subgroup of permutations which preserve each subset $F_v$, then this is just the product of the previously defined actions on the factors $\scrR^{|v|,\bfp_v}$, plus trivial actions on $(0,1]$. Outside that subgroup, the action exchanges different connected components. The same applies to the partial compactification.

\section{Wrapped Floer cohomology\label{sec:wrapped}}

In this section, we construct the wrapped Floer cochain complex and its $A_\infty$-structure. As mentioned before, this should be thought of as a kind of direct limit, although the definition itself is not explicitly formulated in such terms. We emphasize again that, even though the cochain level construction itself is new, the underlying cohomology level theory is well-known to specialists, as one will see by looking at \cite{abbondandolo-schwarz06}. Our focus will be on the geometry of the relevant moduli spaces of (deformed) pseudo-holomorphic maps. The most interesting aspect is that the contribution of most of these spaces is zero, because of cancellations induced by orientation-reversing symmetries; in fact, this cancellation phenomenon is crucial in order to obtain the $A_\infty$-associativity relations. Analytic issues, which do not go beyond the standard technological level, are relegated to Sections \ref{sec:compactness}--\ref{sec:signs}. To conclude, we should mention that, in principle, there is an alternative and more direct definition of wrapped Floer cohomology using a single Hamiltonian function with unboundedly increasing slope instead of a sequence of functions with constant slope. However, the alternative approach has its own technical quirks and is not as suitable for discussing Viterbo functoriality.

\subsection{Basic geometry}
Let $M$ be a Liouville domain. By definition, $M$ is a compact $2n$-dimensional manifold with boundary, together with a one-form $\theta$ such that $\omega = d\theta$ is symplectic, and the dual Liouville vector field $Z$, uniquely determined by the requirement that $i_Z\omega = \theta$, points strictly outwards along $\partial M$. This implies that $\theta|\partial M$ is a contact one-form. By flowing inward  from the boundary along $Z$, one obtains a collar embedding $(0,1] \times \partial M \hookrightarrow M$ modelled on the small half of the symplectization of $\partial M$. One can then complete $M$ by attaching the big half:
\begin{equation} \label{eq:completion}
\hat{M} = M \cup_{\partial M} ([1,\infty) \times \partial M).
\end{equation}
The piece $[1,\infty) \times \partial M$ of \eqref{eq:completion} is called the {\em infinite cone}. The completion carries a natural one-form $\hat\theta$ such that $d\hat\theta = \hat\omega$ is symplectic, and an associated Liouville vector field $\hat{Z}$. On $M \subset \hat{M}$ these restrict to the previously given data, while on the overlapping piece $(0,\infty) \times \partial M \subset \hat{M}$ we have $\hat\theta = r(\theta|\partial M)$ and $\hat{Z} = r\partial_r$; here $r$ is the variable in $(0,\infty$). We will use a particular class of Hamiltonian functions $H \in \smooth(M,\R)$, namely those for which:
\begin{align} \label{eq:class-of-h}
& \parbox{30em}{$H > 0$ everywhere. Moreover, $H$ should admit a smooth extension $\hat{H}$ to the completion $\hat{M}$ such that $\hat{H}(r,y) = r$ on the infinite cone.}
\end{align}
Let $\hat{X}$ be the Hamiltonian vector field of $\hat{H}$. On the infinite cone we have $\hat{X} = (0,R)$, where $R$ is the Reeb vector field associated to $\theta|\partial M$. By restricting, one sees that the Hamiltonian vector field $X$ of $H$ satisfies $X|\partial M = R$.

Let $L \subset M$ be a Lagrangian submanifold intersecting $\partial M$ transversally, and which has the following property:
\begin{equation} \label{eq:asymptotically-parallel}
\parbox{30em}{
$\theta|L \in \Omega^1(L)$ is exact, $\theta|L = df$. Moreover, $\theta|L \in \Omega^1(L)$ vanishes to infinite order along the boundary $\partial L = L \cap \partial M$.}
\end{equation}
The second part implies that $\partial L$ is a Legendrian submanifold of $\partial M$, and also that, by attaching an infinite cone to that boundary, one can extend $L$ to a smooth noncompact Lagrangian submanifold
\begin{equation} \hat{L} = L \cup_{\partial L} ([1,\infty) \times \partial L) \subset \hat{M}. \end{equation} Recall that a Reeb chord of length $w>0$ is a trajectory $x: [0,1] \rightarrow \partial M$ of $wR$, such that $x(0), x(1) \in \partial L$. An integer Reeb chord is one which has integer length. We will make the following technical assumption:
\begin{align} \label{eq:no-integer-chords}
& \parbox{30em}{There are no integer Reeb chords.}
\end{align}
In parallel with the previous terminology, an {\em integer $X$-chord} is a trajectory $x: [0,1] \rightarrow M$ of $wX$, for some positive integer $w$, such that $x(0), x(1) \in L$. In these terms, \eqref{eq:no-integer-chords} says that all such trajectories must lie in the interior of $M$. Additionally, we assume:
\begin{align}
\label{eq:nondegenerate-x}
& \parbox{30em}{All integer $X$-chords are nondegenerate.} \\
\label{eq:nonconcatenate-x}
& \parbox{30em}{No point of $L$ is at once a starting point of an integer $X$-chord, and an endpoint of one (either the same chord or a different one).}
\end{align}
Note that as a consequence of either of these two assumptions, no critical point of $H$ may lie on $L$. Besides that, we will also impose certain topological conditions which are traditional in Floer theory:
\begin{equation} \label{eq:standard-floer}
\parbox{30em}{
The relative Chern class $2 c_1(M,L) \in H^2(M,L;\Z)$ and the second Stiefel-Whitney class $w_2(L) \in H^2(L;\Z/2)$ both vanish.
}
\end{equation}

It seems appropriate to briefly discuss the significance of these conditions. There are at least two alternative versions of \eqref{eq:asymptotically-parallel}. One is stronger:
\begin{equation} \label{eq:stronger}
\parbox{30em}{$\theta|L$ is exact and vanishes on a neighbourhood of $\partial L$.}
\end{equation}
The second part implies that the Liouville flow is parallel to $L$ near its boundary, making it locally into a cone over $\partial L$. The other version is weaker:
\begin{equation} \label{eq:weaker}
\parbox{30em}{$\theta|L$ is exact and the restriction $\theta|\partial L \in \Omega^1(\partial L)$ is zero.}
\end{equation}
The second part of this is equivalent to saying that $\partial L$ is Legendrian. While \eqref{eq:weaker} may be the most natural assumption, it is not directly suitable for our construction. On the other hand, simply by allowing  isotopies which are constant on the boundary, the difference between any of these versions disappears.

\begin{lemma} \label{th:bend-boundary}
Let $L_0$ be a Lagrangian submanifold satisfying \eqref{eq:weaker}. Then there is a Hamiltonian isotopy rel $\partial M$, which deforms $L_0$ into another submanifold $L_1$ satisfying \eqref{eq:stronger}.
\end{lemma}

\proof We start with a construction which is local near $\partial L_0$. By assumption, $\theta|L_0$ vanishes along $\partial L_0$, hence is exact nearby. This means that we can find a function $h \in \smooth(M,\R)$ vanishing along $\partial M$, such that $(\theta + dh)|L_0 \in \Omega^1(L_0)$ is zero in a neighbourhood of $\partial L_0$. The Hamiltonian vector field $Y$ of $h$ is parallel to $\partial M$. Hence, the Liouville vector field of $\theta_t = \theta + (1-t) dh$, which is $Z_t = Z + (1-t) Y$, still points outwards along $\partial M$. Take $\partial L_0$ and push it inwards using $-Z_t$, for times in some interval $[0,\infty)$. The orbits of each such flow yield a local Lagrangian submanifold $L_t^{\mathrm{loc}}$ with boundary $\partial L_t^{\mathrm{loc}} = \partial L_0$. More precisely, we consider $L_t^{\mathrm{loc}}$ as a germ defined near the boundary, hence reserve the right to shrink it whenever necessary. For $t = 0$ this germ agrees with that of $L_0$, since $Z_0$ is parallel to that submanifold; on the other hand, for $t = 1$ we have $\theta|L_1^{\mathrm{loc}} = 0$ by construction.

This local Lagrangian isotopy is described by a family of closed one-forms $\beta_t \in \Omega^1(L_t^{\mathrm{loc}})$, each of which is zero along the boundary. Take functions $g_t \in \smooth(M,\R)$ vanishing along the boundary, such that $dg_t|L_t^{\mathrm{loc}} = \beta_t$. By construction, the Hamiltonian isotopy generated by the family $g_t$ is constant along the boundary, and moves $L_0$ to a Lagrangian submanifold $L_1$ which extends the previously given germ $L_1^{\mathrm{loc}}$. \qed

Next, among the more technical conditions, \eqref{eq:no-integer-chords} can be satisfied by a small rescaling of $\omega$ and $\theta$, simply because the periods of all Reeb chords form a measure zero subset of $\R$ (by a suitable application of Sard's theorem). Suppose that this has been done, and additionally that $2n = \mathrm{dim}\,M \geq 4$. Then, a generic choice of $H$ ensures that \eqref{eq:nondegenerate-x} and \eqref{eq:nonconcatenate-x} hold (Lemma \ref{th:nondegenerate-g} and \ref{th:nonconcatenate-g}, respectively). Unfortunately, in the two-dimensional case the second of these conditions is not generic, and a workaround is needed (Section \ref{subsec:stabilization}). Finally, one can drop \eqref{eq:standard-floer} altogether, at the cost of the usual loss of Floer-theoretic structure. This would mean working with coefficients in $\Z/2$ and obtaining only an ungraded $A_\infty$-structure, since degrees and signs can no longer be assigned consistently.

\subsection{Almost complex structures\label{subsec:j}}
We begin by recalling a classical definition from \cite{hofer93}. Suppose that we are given an almost complex structure on the contact hyperplane field $\xi = ker(\theta|\partial M) \subset T(\partial M)$, compatible with its symplectic structure. From that, one constructs an almost complex structure on the symplectization $(0,\infty) \times \partial M$ as follows: write
\begin{equation} \label{eq:standard-splitting}
T((0,\infty) \times \partial M) = \R^2 \oplus \xi,
\end{equation}
where the first factor is spanned by $r\partial_r$ and $(0,R)$. Take the direct sum of the standard complex structure on the first summand, meaning that $r\partial_r \mapsto (0,R)$, and of the given one on the second summand. The almost complex structures on the symplectization constructed in this way are called {\em of contact type}. Any such structure is invariant under the flow of $r\partial_r$.

Take a $\hat\omega$-compatible almost complex structure $\hat{I}$ on $\hat{M}$, whose restriction to the infinite cone is of contact type. The restriction $I = \hat{I}|M$ is called an almost complex structure {\em of contact type at the boundary}. Note that $I$ determines $\hat{I}$ uniquely. Denote the space of all such structures by $\scrJ(M)$. This is an infinite-dimensional manifold in a loose sense. Its tangent space $T\scrJ(M)$ at any point $I$ consists of endomorphisms $K$ of $TM$ which anti-commute with $I$, which have the property that $\omega(\cdot, K \cdot)$ is symmetric, and which moreover are restrictions of suitable endomorphisms of $T\hat{M}$ (we omit the details of the last condition, which are straightforward but lengthy). We call these endomorphisms {\em infinitesimal deformations} of the almost complex structure. One can get an actual deformation by exponentiating $K$, which means setting
\begin{equation} \label{eq:perturbed}
J = I \, \exp(-IK).
\end{equation}

Let $S$ be a disc equipped with a set of strip-like ends $\bfepsilon$. Let $\bfI = \{I_z\}$ be a family of almost complex structures in $\scrJ(M)$ parametrized by points $z \in S$. We say that this family is {\em strictly compatible with the strip-like ends} if for each $k \in \{0,\dots,d\}$ there is a family of such structures $\bfI^k = \{I^k_t\}$ parametrized by $t \in [0,1]$, such that
\begin{equation} \label{eq:j-on-the-end}
I_{\epsilon^k(s,t)} = I^k_t
\end{equation}
when $\pm s \gg 0$. In other words, if we go sufficiently far along any end, $\bfI$ becomes independent of $s$, hence compatible with translations in the $s$ direction.

A universal choice of almost complex structures consists of the following data. For each stable weighted popsicle $(S,\bfphi,\bfw)$, choose a family $\bfI_{S,\bfphi,\bfw}$ parametrized by $S$, and strictly compatible with (the previous universal choice of) strip-like ends. We require that these structures should vary smoothly over the moduli spaces $\scrR^{d+1,\bfp,\bfw}$, and be invariant with respect to the action of $\Sym^{\bfp}$. Additionally, for each $w$ choose a family $\bfI_w$ parametrized by $[0,1]$. This can actually be viewed as an extension of the previous choice to the unstable case, invariant under automorphisms:
\begin{equation} \label{eq:unstable-j}
I_{Z,\bfw,s+it} = I_{w,t},
\end{equation}
where $\bfw = \{w^0 = w^1 = w\}$. A universal choice is called consistent if the following two conditions are satisfied. First of all, over the $k$-th strip-like end the behaviour of $\bfI_{S,\bfphi,\bfw}$, in the sense of \eqref{eq:j-on-the-end}, is actually governed by $\bfI_{w^k}$. Secondly, the choice is compatible with gluing together weighted popsicles, for small values of the gluing parameter. Here, the gluing process is defined by our previous universal choices of strip-like ends and popsicle sticks, as in Section \ref{subsec:sticks} (and carrying over the weights in the obvious way).

\begin{remark}
The consistency requirements ensure the existence of well-behaved compactifications for the moduli spaces of holomorphic maps to be defined later on. However, they are not convenient from the point of view of transversality results, because the behaviour of the almost complex structure on parts of the Riemann surface is restricted (in a non-local way, even). To see an example of this, suppose that $(S,\bfphi,\bfw)$ is obtained from $(S_+,\bfphi_+,\bfw_+)$ and two copies of $(S_-,\bfphi_-,\bfw_-)$ (glued into different positive ends of $S_+$, which necessarily should have the same weight). Then, consistency dictates that the family of almost complex structures associated to $S$ should be equal on the two pieces coming from $S_-$. To avoid these problems, we will introduce perturbations of the almost complex structures which allow the conditions to be relaxed, while retaining enough asymptotic control for the purposes of compactification.
\end{remark}

Take a disc with a set of strip-like ends, and suppose that we have chosen $\bfI$ on it satisfying \eqref{eq:j-on-the-end}. An infinitesimal deformation of $\bfI$ is a family $\bfK = \{K_z\}$, where each $K_z$ lies in the tangent space to $T\scrJ(M)$ at $I_z$, and with the property that as $z = \epsilon^k(s,t)$ goes to infinity along a strip-like end, $K_z$ and all its derivatives go to zero faster than any exponential function $\exp(-C|s|)$. Equivalently, let $\overline{\bfK}$ be the trivial extension of $\bfK$ to $\bar{S}$; the condition is that this extension is smooth, and vanishes to infinite order at the points at infinity. Then, the perturbed family $\bfJ = \{J_z\}$ obtained by applying \eqref{eq:perturbed} at every point of $S$ agrees asymptotically with the original $\bfI$, with good control over the rate of convergence.

We now consider the corresponding notion for families. Assume we have already made a universal choice of almost complex structures, assumed to be consistent. A universal infinitesimal deformation
consists of an infinitesimal deformation $\bfK_{S,\bfphi,\bfw}$ of each $\bfI_{S,\bfphi,\bfw}$, which varies smoothly with respect to the moduli, and is compatible with the action of $\Sym^{\bfp}$. In slightly different terms, the universal infinitesimal deformation is given by a $\Sym^{\bfp}$-invariant section of the pullback bundle
\begin{equation} \label{eq:pullback-end}
\mathit{End}(\mathit{TM}) \longrightarrow \scrS^{d+1,\bfp,\bfw} \times M
\end{equation}
for each $(d,\bfp,\bfw)$. We say that the universal infinitesimal deformation is {\em asymptotically consistent} if the sections extend smoothly to the compactifications $\bar\scrS^{d+1,\bfp,\bfw} \times M$, and the following two additional properties hold.

For the first condition, fix a pair $(T,\bfF)$ which represents a boundary stratum $\scrR^{T,\bfF,\bfw} \subset \bar\scrR^{d+1,\bfp,\bfw}$. For every vertex $v$ of $T$, with its associated $\bfp_v$ and weights $\bfw_v$, we take the projection $\pi_v : \scrR^{T,\bfF,\bfw} \longrightarrow \scrR^{|v|,\bfp_v,\bfw_v}$, use that to pull back the universal family, then take the disjoint union over all $v$. By definition of the compactification, there is a natural embedding
\begin{equation} \label{eq:compare-c}
\coprod_v \pi_v^* \scrS^{|v|,\bfp_v,\bfw_v} \longrightarrow \bar\scrS^{d+1,\bfp,\bfw}.
\end{equation}
Given a universal infinitesimal deformation, there are two ways of making a section of $\mathit{End}(\mathit{TM}) \rightarrow \coprod_v \pi_v^*\scrS^{|v|,\bfp_v,\bfw_v} \times M$. We could take the sections attached to the smaller infinitesimal families $\scrS^{|v|,\bfp_v,\bfw_v}$ and combine them; on the other hand, we could take the section attached to $\scrS^{d+1,\bfp,\bfw}$, extend it to the compactification, and then pull it back via \eqref{eq:compare-c}. The condition then says that these two constructions yield the same result, for all $(T,\bfF)$.

For the second condition, take a point of the compactified moduli space $\bar\scrS^{d+1,\bfp,\bfw}$ which is represented by a pair $(\bar{S},z)$, where $z$ is either a point at infinity of $S$ or (in the case of a broken popsicle) a singular point. Then, our section should vanish when restricted to that point (times $M$), and so should all its derivatives.

From now on, we will work with a fixed consistent universal choice of almost complex structures, as well as a fixed asymptotically consistent infinitesimal deformation. By applying \eqref{eq:perturbed} at each point of $S$, one can then form perturbed families $\bfJ_{S,\bfphi,\bfw}$. By definition, this kind of perturbation does not affect the limits of our almost complex structures over the ends of $S$, so we set
\begin{equation} \label{eq:no-perturbation}
\bfJ_w = \bfI_w, \;\; \bfJ_{Z,\bfw} = \bfI_{Z,\bfw} \text{ for $\bfw = (w^0 = w^1 = w)$}.
\end{equation}

\begin{remark} \label{th:concrete-consistency}
To understand the concrete implications of asymptotic consistency, one has to look at the geometry of the compactified universal family, and more specifically at the behaviour of functions with similar vanishing conditions at special points. This is done in Section \ref{subsec:singularity}, and the outcome of that discussion can be summarized as follows. As one approaches a point in the boundary of $\bar\scrR^{d+1,\bfp,\bfw}$, the fibre of the universal family acquires ``necks'' which are finite strips of increasing length. On each of these finite strips, as well as on the strip-like ends, the infinitesimal deformation $\bfK_{S,\bfphi,\bfw}$ decays faster than exponentially, in the sense of \eqref{eq:superexponential-1}, \eqref{eq:superexponential-2}. Hence, $\bfJ_{S,\bfphi,\bfw}$ converges everywhere to the structure which is inherited from the components of the limiting broken popsicle, and moreover, the rate of convergence is precisely controlled on the finite strips and strip-like ends.
\end{remark}

\subsection{Pseudo-holomorphic maps\label{subsec:define}}
For each $w$, let $\scrX_w$ be the set of integer $X$-chords of length $w$. From \eqref{eq:nondegenerate-x} we know that each set $\scrX_w$ is  finite.

Take a weighted popsicle $(S,\bfphi,\bfw)$, equipped with strip-like ends, as well as sets of basic closed and sub-closed one-forms. Let $\gamma$ be the induced total sub-closed one-form \eqref{eq:total-one-form}. In addition, we assume that we have a family $\bfI$ which is compatible with the strip-like ends, as well as an infinitesimal deformation $\bfK$, which together define the perturbed family $\bfJ$. Pick $\bfx = \{x^0,\dots,x^d\}$, where each $x^k$ is an element of $\scrX_{w^k}$. Consider the inhomogeneous $\bar\partial$-equation
\begin{equation} \label{eq:dbar}
\left\{
\begin{aligned}
 & u: S \longrightarrow M, \\
 & u(\partial S) \subset L, \\
 & \textstyle\lim_{s \rightarrow \pm\infty} u(\epsilon^k(s,\cdot)) = x^k(\cdot), \\
 & (du_z - X_{u(z)} \otimes \gamma_z) \circ j + J_{z,u(z)} \circ (du_z - X_{u(z)} \otimes \gamma_z) = 0.
\end{aligned}
\right.
\end{equation}
Here $j$ is the complex structure on $S$; $X \otimes \gamma$ is the section of $\mathit{Hom}(\mathit{TS},u^*\mathit{TM})$ obtained by composing $\gamma \in \smooth(\mathit{TS}^*)$ with $u^*X \in \smooth(u^*\mathit{TM})$; and the overall condition says that when we subtract that inhomogeneous term from $du$, the outcome is of type $(1,0)$ with respect to $j$ and $J_z$. We will sometimes abbreviate this equation by
\begin{equation}
 (du - X \otimes \gamma)^{0,1} = 0.
\end{equation}

Fix $d \geq 1$ and $\bfp = \{p_f\}_{f \in F}$, assuming first that these satisfy the stability condition $d+|F| \geq 2$. Choose weights $\bfw$ satisfying \eqref{eq:weights}, and limits $\bfx$ as before. The {\em space of stable popsicle maps}, denoted by
\begin{equation} \label{eq:moduli-space}
\scrR^{d+1,\bfp,\bfw}(\bfx) \quad \text{(or $\scrR^{d+1,\bfw}(\bfx)$ if $F = \emptyset$)},
\end{equation}
is the set of triples $(S,\bfphi,u)$ of the following kind. $(S,\bfphi)$ is a $\bfp$-flavoured popsicle, representing a point of $\scrR^{d+1,\bfp} = \scrR^{d+1,\bfp,\bfw}$. As such, it carries additional data specified by our previous choices, in particular: ends $\bfepsilon_{S,\bfphi}$, a one-form $\gamma_{S,\bfphi,\bfw}$, and a family of almost complex structures $\bfJ = \bfJ_{S,\bfphi,\bfw}$. Using all that, one can write down \eqref{eq:dbar} on $S$, and $u$ is then a solution of that equation, with limits given by $\bfx$. Note that \eqref{eq:moduli-space} carries an induced action of $\Sym^{\bfp}$, which consists of permuting the sprinkles $\phi_f$ while keeping everything else the same.

It remains to extend the definition of the moduli spaces to the case where the underlying popsicle is unstable, namely $d = 1$ and $F = \emptyset$. Equip $S = Z$ with \eqref{eq:unstable-gamma} as well as the family $\bfJ = \bfJ_{Z,\bfw}$ of almost complex structures taken from \eqref{eq:no-perturbation}, \eqref{eq:unstable-j}. Then \eqref{eq:dbar} reduces to the classical equation for Floer trajectories:
\begin{equation} \label{eq:old-floer}
\left\{
\begin{aligned}
 & u: Z \longrightarrow M, \\
 & u(\partial Z) \subset L, \\
 & \textstyle\lim_{s \rightarrow -\infty} u(s,\cdot) = x^0(\cdot), \\
 & \textstyle\lim_{s \rightarrow +\infty} u(s,\cdot) = x^1(\cdot), \\
 & \partial_s u + J_{t,u(t)}(\partial_t u - w X) = 0.
\end{aligned}
\right.
\end{equation}
As usual, one defines $\scrR^{2,\bfw}(\bfx)$ to be the space of non-stationary solutions of that equation (solutions for which $\partial_s u$ does not vanish identically), divided by the translational $\R$-action.

\begin{convention} \label{th:waffle}
From now on, we will lump both cases together. This means that, even though we will mostly use the notation $(S,\bfphi,u)$ which is appropriate for stable popsicles, it is tacitly understood that the discussion includes the case of Floer moduli spaces, with the necessary minor modifications.
\end{convention}

\subsection{Smoothness\label{subsec:smoothness}}
Let $(S,\bfphi,u)$ be a point in \eqref{eq:moduli-space}. The virtual dimension of the moduli space at that point is
\begin{equation} \label{eq:vdim}
 \mathrm{vdim}\,\scrR^{d+1,\bfp,\bfw}(\bfx) = \mathrm{dim}\,\scrR^{d+1,\bfp} +
 \deg(x^0) - \deg(x^1) - \cdots - \deg(x^d).
\end{equation}
Here, $\deg(x)$ is the Maslov index of $x \in \scrX_w$ (in general, this depends on a choice of grading). In the case of Floer's equation, the dimension is set to ``$\dim\,\scrR^2 = -1$'', recovering the standard formula.

\begin{theorem} \label{th:smoothness-1}
Suppose that the almost complex structures $\bfI_w$ have been chosen generically. Then, for a generic choice of infinitesimal deformations $\bfK_{S,\bfphi,\bfw}$, subject to asymptotic compatibility, all the moduli spaces $\scrR^{d+1,\bfp,\bfw}(\bfx)$ are regular, hence smooth manifolds of the expected dimension \eqref{eq:vdim}.
\end{theorem}

The action of $\Sym^{\bfp}$ on $\scrR^{d+1,\bfp,\bfw}(\bfx)$ has the same structure as the underlying action on the moduli space of popsicles. By this, we mean that the isotropy groups which appear are the subgroups $\Sym^{\bfP}$ from Addendum \ref{th:isotropy}. The virtual codimension of the corresponding fixed point set is
\begin{equation} \label{eq:codim}
\mathrm{vcodim}\,\scrR^{d+1,\bfp,\bfw}(\bfx)^{\Sym^{\bfP}} = \sum_m (|P_m| - 1).
\end{equation}

\begin{theorem} \label{th:smoothness-2}
Generically, in the same sense as in Theorem \ref{th:smoothness-1}, the subspaces of $\scrR^{d+1,\bfp,\bfw}(\bfx)$ fixed by any subgroup $\Sym^{\bfP}$ are themselves regular, hence smooth submanifolds of codimension \eqref{eq:codim}.
\end{theorem}

The proofs of Theorem \ref{th:smoothness-1} and \ref{th:smoothness-2} are pretty much standard transversality theory (see Section \ref{subsec:linear} for details). From now on, we will assume throughout our discussion that all the almost complex structures involved have been chosen generically, so that the conclusions of these theorems hold.

\subsection{Compactness\label{subsec:gromov}}
The spaces \eqref{eq:moduli-space} admit compactifications whose structure is modelled on \eqref{eq:compactification}. Suppose that we have a pair $(T,\bfF)$ as in the definition of broken popsicle. Define flavours $\bfp_v$ as in Section \ref{subsec:define-popsicle}, and weights $\bfw_v$ as in Section \ref{subsec:weights}. Moreover, for each $v$ we want to have a collection $\bfx_v = \{x_v^k\}$, where each $x_v^k \in \scrX_{w_v^k}$, with the following properties:
\begin{equation} \label{eq:x-properties}
\parbox{30em}{
If the pair $(v,k)$ correspond to the $j$-th semi-infinite edge of $T$, then $x_v^k = x^j$ is one of the given $X$-chords. Otherwise, given a finite edge of $T$, and the two associated flags $(v_+,k)$ and $(v_-,0)$, we require that $x_{v_+}^k = x_{v_-}^0$.
}
\end{equation}
Finally, for every $v$ we want to have an element
\begin{equation} \label{eq:broken-map}
(S_v,\bfphi_v,u_v) \in \scrR^{|v|,\bfp_v,\bfw_v}(\bfx_v).
\end{equation}
Then, the collection $\{(S_v,\bfphi_v,u_v)\}$, which we call a broken popsicle map, defines a point in the compactified moduli space. As the terminology indicates, this space is compact in an appropriately defined Gromov topology. A basic role in the compactness argument is played by an priori estimate for the energy $||du-X \otimes \gamma||_{L^2}$, derived from an inequality between that norm and the action functional; see Section \ref{subsec:energy}, or more specifically \eqref{eq:energy-inequality} and \eqref{eq:a-priori}. Other important players are the various consistency conditions; these ensure that, when the domain degenerates into a broken popsicle, the limiting map satisfies \eqref{eq:dbar} on each component of that popsicle, with the correct almost complex structure and inhomogeneous term. The proof of compactness starts with these observations, and then proceeds exactly as in the case of ordinary Fukaya categories. We will not comment on it further.

The simplest application is when the virtual dimension is zero. Then the boundary of the compactification is empty, hence $\scrR^{d+1,\bfp,\bfw}(\bfx)$ itself is a finite set. Next, suppose that the virtual dimension is one. Then $\bar\scrR^{d+1,\bfp,\bfw}(\bfx)$ is a one-dimensional compact manifold with boundary. Boundary points correspond to broken pseudo-holomorphic maps, where $T$ is a tree with two vertices, and each component of the broken map is itself part of a zero-dimensional moduli space. Part of this statement requires a suitable gluing theorem, which we again take for granted, given its similarity with established results for Fukaya categories. Another point, which is more specific to the geometric situation here, is the following one: in principle, one could think that there might be other boundary points of one-dimensional moduli spaces, related to the fact that our pseudo-holomorphic maps can hit the boundary of $M$. It turns out that these points are actually not boundary points, essentially due the convexity of $\partial M$ (see Section \ref{subsec:escape} for a proof; alternatively one could use the maximum principle, as in \cite{oancea04b} or \cite{khovanov-seidel98}).

\subsection{Algebraic relations\label{subsec:algebraic-relation}}
At this point, we fix a coefficient field $\K$. We also choose an orientation of each popsicle moduli space $\scrR^{d+1,\bfp}$. Finally, we make a choice of grading and $Pin$ structure for our Lagrangian submanifold. The first choice fixes the Maslov indices $\deg(x) \in \Z$ of all $x \in \scrX_w$, and the second one allows
one to associate to any $x$ a one-dimensional $\K$-vector space, the so-called {\em $\K$-normalized orientation space} $|o_x|_{\K}$ (see Section \ref{subsec:grading-and-stuff} for a brief recapitulation). Moreover, for every point $(S,\bfphi,u)$ in a zero-dimensional moduli space $\scrR^{d+1,\bfp,\bfw}(\bfx)$, we get a preferred isomorphism
\begin{equation} \label{eq:orientation-map}
|o_{S,\bfphi,u}^{red}|_{\K} : |o_{x^d}|_\K \otimes \cdots \otimes |o_{x^1}|_\K \longrightarrow |o_{x^0}|_{\K}.
\end{equation}
This is, in a slightly more abstract formulation than usual, the sign with which $(S,\bfphi,u)$ contributes to the algebraic count of points in the moduli space. Denote by $m^{d,\bfp,\bfw}(\bfx)$ the sum of the homomorphisms \eqref{eq:orientation-map}
\begin{equation*}m^{d,\bfp,\bfw}(\bfx) = \sum_{(S,\bfphi,u)} |o_{x^d}|_\K.   \end{equation*}
 When $F = \emptyset$ we may omit $\bfp$ and just write $m^{d,\bfw}(\bfx)$, in parallel with the notation for the spaces \eqref{eq:moduli-space} themselves.

\begin{lemma} \label{th:vanishing}
If $\Sym^\bfp$ is nontrivial, $m^{d,\bfp,\bfw}(\bfx)$ vanishes.
\end{lemma}

To see why this is the case, recall that Theorem \ref{th:smoothness-2} implies that all strata of $\scrR^{d+1,\bfp,\bfw}(\bfx)$ consisting of points with nontrivial isotropy group have positive codimension. Hence, in the zero-dimensional case, $\Sym^{\bfp}$ necessarily acts freely on the moduli space. Now, any two points related by a transposition in that group actually contribute with opposite sign, hence cancel out (for details, see Lemma \ref{th:cancellation}).

The vanishing argument outlined above applies exactly when the map $\bfp: F \rightarrow \{1,\dots,d\}$ fails to be injective (or to put it more intuitively, if there is a popsicle stick carrying more than one sprinkle). We now turn to the remaining injective case. In that situation, it is more natural think of $F$ as a subset of $\{1,\dots,d\}$, with $\bfp$ being the inclusion, and we correspondingly replace $\bfp$ with $F$ everywhere in the notation. Fix $d+1$, $F$, $\bfw$ and $\bfx$ such that the moduli space $\scrR^{d+1,F,\bfw}(\bfx)$ is one-dimensional.

\begin{definition} \label{th:cut}
An {\em admissible cut} of $F$ consists of $d_+,d_- \geq 1$ such that $d_- + d_+ = d+1$, a number $i \in \{1,\dots,d_+\}$, and subsets $F_\pm \subset \{1,\dots,d_\pm\}$ satisfying $|F_-| + |F_+| = |F|$, with the following property:
\begin{equation}
\parbox{30em}{
$F$ contains all $k \in F_+$ satisfying $k < i$, the numbers $k+d_- - 1$ for all $k \in F_+$ with $k>i$, and the numbers $k+i-1$ for $k \in F_-$. If $i \notin F_+$, this completely describes $F$. Otherwise, $F$ has one element other than the ones we have already given, which lies in the range $\{i,\dots,i+d_- -1\}$.
}
\end{equation}
Note that these subsets automatically come with canonical injective maps $\iota_\pm: F_\pm \rightarrow F$, which together cover the whole of $F$. One of them is simply $\iota_-(k) = k+i-1$. The other one satisfies $\iota_+(k) = k$ for $k<i$, $\iota_+(k) = k+d_--1$ for $k>i$, and maps $i \in F_+$ (if that is the case) to the unique element in $F$ specified above.
\end{definition}

Given such a cut, set $w^{\mathit{new}} = w^i + \cdots + w^{i+d_- - 1} + |F_-|$. Choose any $x^{\mathit{new}} \in \scrX_{w^{\mathit{new}}}$ such that $\deg(x^{\mathit{new}}) = \deg(x^i) + \cdots + \deg(x^{i+d_- - 1}) + 2 - d_- - |F_-|$. Then, define weights $\bfw_\pm = \{w_\pm^0,\dots,w_\pm^{d_\pm}\}$ and collections of chords $\bfx_\pm = \{x_\pm^0,\dots,x_\pm^{d_\pm}\}$ as follows:
\begin{equation} \label{eq:w-x}
\begin{aligned}
 &
 (w_+^k,x_+^k) = \begin{cases} (w^k,x^k) & k<i, \\ (w^{\mathit{new}},x^{\mathit{new}}) & k = i, \\
 (w^{k+d_--1},x^{k+d_--1}) & k>i \end{cases}, \\
 &
 (w_-^k,x_-^k) = \begin{cases} (w^{\mathit{new}},x^{\mathit{new}}) & k = 0, \\
 (w^{k+i-1},x^{k+i-1}) & k > 0.
 \end{cases}
\end{aligned}
\end{equation}
In this terminology, $\scrR^{d_+ + 1,F_+,\bfw_+}(\bfx_+) \times \scrR^{d_- + 1,F_-,\bfw_-}(\bfx_-)$ appears as one of the boundary strata in the compactification of $\scrR^{d,F,\bfw}(\bfx)$. In fact, this construction describes precisely those strata which have a trivial group of symmetries. The contributions from the other boundary strata cancel, for a reason parallel to that in Lemma \ref{th:vanishing} (see Section \ref{subsec:1-dim} for the precise argument). We therefore get a relation
\begin{equation} \label{eq:algebraic-relation}
\sum (-1)^\S \, m^{d_+, F_+,\bfw_+}(\bfx_+) \circ (id^{\otimes d_+ - i} \otimes
m^{d_-,F_-,\bfw_-}(\bfx_-) \otimes id^{\otimes i-1}) = 0,
\end{equation}
where the sum is over all admissible cuts $(d_\pm,i,F_\pm)$ and all possible $x^{\mathit{new}}$. The sign depends on the choice of orientation of the underlying spaces $\scrR^{d+1,F}$ (see Section \ref{subsec:choose-orientation} for a detailed discussion). For a specific such choice, it turns out to be
\begin{equation} \label{eq:orientation-signs}
\begin{aligned}
\S = &\; d_-i + i + 1 + d_+ |F_-| \\
 & + (d_- + |F_-|)(\deg(x^{i+d_-}) + \cdots + \deg(x^d)) \\ & +
|\{(k_+,k_-) \in F_+ \times F_- \;:\; \iota_+(k_+) < \iota_-(k_-)\}|.
\end{aligned}
\end{equation}

\begin{examples} (i) \label{th:algebraic-relation}
Take $d = 1$ and $F = \{1\}$. In that case, the admissible cuts necessarily have $d_\pm = 1$. One of the two sets $F_\pm$ is $\{1\}$ while the other is empty. Omitting weights for the sake of brevity, we find that \eqref{eq:algebraic-relation} specializes to
\begin{equation}
\begin{aligned}
 & \textstyle\sum m^{1,F}(x^0,x^{\mathit{new}}) \circ m^1(x^{\mathit{new}},x^1) \\ & -
 \textstyle\sum m^1(x^0,x^{\mathit{new}}) \circ m^{1,F}(x^{\mathit{new}},x^1) = 0,
\end{aligned}
\end{equation}
where the sum is over all $x^{\mathit{new}}$ with the appropriate Maslov index. Geometrically, this kind of degeneration corresponds to an ordinary Floer trajectory bubbling off at the ends $\pm\infty \in \bar{Z}$.

(ii) Take $d = 2$ and $F = \{1,2\}$, which is the situation from Example \ref{th:ex-popsicle}(ii). Out of the five boundary edges in $\bar\scrR^{3,F}$, four give nontrivial contributions to \eqref{eq:algebraic-relation}, while the remaining one cancels for symmetry reasons. Again omitting weights, these contributions are
\begin{equation}
\begin{aligned}
 & -\textstyle\sum m^{1,\{1\}}(x^0,x^{\mathit{new}}) \circ m^{2,\{1\}}(x^{\mathit{new}},x^1,x^2)
 \\ & + \textstyle\sum m^{1,\{1\}}(x^0,x^{\mathit{new}}) \circ m^{2,\{2\}}(x^{\mathit{new}},x^1,x^2)
 \\ &
 + \textstyle\sum m^{2,\{1\}}(x^0,x^1,x^{\mathit{new}}) \circ (m^{1,\{1\}}(x^{\mathit{new}},x^2) \otimes id) \\ &
 - \textstyle\sum m^{2,\{2\}}(x^0,x^{\mathit{new}},x^2) \circ (id \otimes m^{1,\{1\}}(x^{\mathit{new}},x^1)).
\end{aligned}
\end{equation}
There are three more summands coming from bubbling off of Floer trajectories as in (i):
\begin{equation}
\begin{aligned}
 & \textstyle\sum m^1(x^0,x^{\mathit{new}}) \circ m^{2,F}(x^{\mathit{new}},x^1,x^2) \\ &
 - \textstyle\sum m^{2,F}(x^0,x^1,x^{\mathit{new}}) \circ (m^1(x^{\mathit{new}},x^2) \otimes id) \\ & +
 (-1)^{\deg(x^2)+1} \textstyle\sum m^{2,F}(x^0,x^{\mathit{new}},x^2) \circ (id \otimes m^1(x^{\mathit{new}},x^1)).
\end{aligned}
\end{equation}
\end{examples}

\begin{remark}
The $m^{d,F,\bfw}$ and their relations are all the analytic input underlying our geometric construction. In particular, moduli spaces with nontrivial symmetry are actually relevant only insofar as they appear in the compactification of other spaces. The crucial fact is that whenever $\bar\scrR^{d+1,F,\bfw}(\bfx)$ one-dimensional, those of its boundary strata which do not correspond to cuts admit orientation-reversing free involutions, hence do not appear in \eqref{eq:algebraic-relation}. With that in mind, the various consistency conditions imposed on the additional geometric data living on $\scrR^{d+1,\bfp,\bfw}$ could be somewhat relaxed when $\Sym^{\bfp}$ is nontrivial (but nothing substantial is gained in doing that).
\end{remark}

\subsection{Wrapped Floer cohomology}
For any $w$, one defines the Floer cochain group (actually a finite-dimensional graded $\K$-vector space) of $L$ with respect to $wH$ as a direct sum of $\K$-normalized orientation spaces, placed in the degree indicated by the Maslov index:
\begin{equation} \label{eq:wrapped}
 CF^*(L; wH) = \bigoplus_{x \in \scrX_w} |o_x|_{\K} [-\deg(x)].
\end{equation}
This group carries the usual differential $\delta$ defined by counting Floer gradient trajectories. In our terminology, the $\bfx = \{x^0,x^1\}$ matrix coefficient of $\delta$ is
\begin{equation*}
 \delta_{(x^0,x^1)} = m^{1,\{w,w\}}(x^0,x^1) : |o_{x^1}|_{\K} \longrightarrow |o_{x^0}|_{\K}.
\end{equation*}
There is one other nontrivial linear operation which appears naturally in our framework. Namely, take $d = 1$ with $F = \{1\}$, so $\bfw = \{w^0 = w+1, w^1 = w\}$. Counting points in these moduli spaces yields maps $\kappa: CF^*(L;wH) \rightarrow CF^*(L;(w+1)H)$, whose matrix coefficients are
\begin{equation}
 \kappa_{(x^0,x^1)} = m^{1,\{1\},\{w+1,w\}}(x^0,x^1) : |o_{x^1}|_{\K} \longrightarrow |o_{x^0}|_{\K}.
\end{equation}
Example \ref{th:algebraic-relation}(i) translates into the statement that $\kappa$ is a chain map (in more standard terminology, this is the continuation map associated to the homotopy from $wH$ to $(w+1)H$; this is a monotone homotopy since $H \geq 0$ by assumption).

Let $q$ be a formal variable of degree $-1$, satisfying $q^2 = 0$. Take the space $CF^*(L;wH)[q] = CF^*(L;wH) \oplus q\,CF^*(L;wH)$, whose elements are formal sums $c = a + bq$. We consider this as a module over the one-dimensional exterior algebra $\C[\partial_q]$, where $\partial_q$ has degree $+1$ and acts by formally differentiating in $q$-direction: $\partial_q(a+q b) = b$. The wrapped Floer cochain space is the infinite direct sum
\begin{equation} \label{eq:wrapped-cochain}
 CW^*(L; H) = \bigoplus_{w = 1}^\infty CF^*(L; wH)[q]
\end{equation}
We equip it with a differential which commutes with $\partial_q$:
\begin{equation} \label{eq:wrapped-differential}
 \mu^1(a + q b) = (-1)^{\deg(a)} \delta(a) + (-1)^{\deg(b)}(q \delta(b) + \kappa(b) - b).
\end{equation}
Visually (and with signs omitted), the resulting chain complex looks like this:
\begin{equation}
 \xymatrix{  {CF^*}{(}L{;}H{)}  \ar@(ur,ul)[]_{\delta} & CF^*(L ; 2H)  \ar@(ur,ul)[]_{\delta} & CF^*(L ; 3H)  \ar@(ur,ul)[]_{\delta} & \cdots  \\
q\,CF^*(L ; H) \ar[u]^{id} \ar[ru]^{\kappa} \ar@(dr,dl)[]^{\delta} & q\,CF^*(L ; 2H) \ar[u]^{id} \ar[ru]^{\kappa}  \ar@(dr,dl)[]^{\delta} & q\,CF^*(L ; 3H)  \ar[u]^{id} \ar[ru]^{\kappa}  \ar@(dr,dl)[]^{\delta} }
\end{equation}
The cohomology of \eqref{eq:wrapped-differential} is called {\em wrapped Floer cohomology} and denoted by $HW^*(L)$ (we omit the function $H$ from the notation because the cohomology groups are independent of it, even though we won't give a proof of that fact here).

\begin{lemma} \label{th:partially-forget}
For every $\nu$, the subcomplex $C^\nu$ consisting of those summands in \eqref{eq:wrapped-cochain} with $w \geq \nu$ is quasi-isomorphic to $CW^*(L;H)$.
\end{lemma}

\proof The $C^\nu$ form a decreasing filtration of $CW^*(L;H)$. Each quotient $C^\nu/C^{\nu+1}$ is the mapping cone of the identity map on $CF^*(L;\nu H)$, hence acyclic. This implies that the inclusions $C^{\nu+1} \rightarrow C^{\nu}$ are quasi-isomorphisms. \qed

\begin{lemma} \label{th:homotopy-limit}
There is a canonical isomorphism
\begin{equation} \label{eq:directlim}
HW^*(L) \iso \underrightarrow{\text{\it lim}}_w \, HF^*(L;wH),
\end{equation}
where the maps in the direct system are induced by $\kappa$.
\end{lemma}

\proof Consider the exhausting increasing filtration of $CW^*(L;H)$ by subcomplexes $C_w = CF^*(L;H)[q] \oplus \cdots \oplus CF^*(L;(w-1)H)[q] \oplus CF^*(L;w H)$, where the last summand has no $q$ component. Obviously, this means that
\begin{equation} \label{eq:c-limit}
HW^*(L;H) \iso \underrightarrow{\text{\it lim}}_w \, H^*(C_w).
\end{equation}
In turn, each $C_w$ itself carries a finite decreasing filtration, $C_w^\nu = CF^*(L;\nu H)[q] \oplus \cdots \oplus CF^*(L;(w-1)H)[q] \oplus CF^*(L;w H)$. For $\nu \neq w$, the quotients $C_w^\nu/C_w^{\nu+1}$ are acyclic, for the same reason as in Lemma \ref{th:partially-forget}.
Hence, the inclusion $C_w^w \rightarrow C_w$ is a quasi-isomorphism. Moreover, the diagram
\begin{equation} \label{eq:h-commute}
 \xymatrix{
 C_w^w \ar@{^{(}->}[d] \ar[rr]^-{\kappa} && C_{w+1}^{w+1} \ar@{^{(}->}[d] \\
 C_w \ar@{^{(}->}[rr] && C_{w+1}
}
\end{equation}
commutes up to a chain homotopy, which is $a \mapsto (-1)^{\deg(a)} qa$. Therefore, if in \eqref{eq:c-limit} we replace the cohomology of $C_w$ by that of its quasi-isomorphic subcomplex $C_w^w$, the connecting maps in the resulting direct system are the ones induced by $\kappa$. Now, $C_w^w$ agrees with $CF^*(L;w H)$ up to a sign change in the differential, which does not affect cohomology. \qed

\begin{remark}
$CW^*(L;H)$ should be seen as the homotopy direct limit of the $CF^*(L;wH)$. Indeed, this is just the chain level version of a well-known construction in classical homotopy theory (see \cite[p.\ 457]{hatcher00}, for instance). Note also that, while wrapped Floer cohomology carries an induced $\C[\partial_q]$-module structure, that structure is actually trivial. This is an easy consequence of the proof above: every class in wrapped Floer cohomology can be represented by a cocycle in $C_w^w$, which is obviously killed by $\partial_q$.
\end{remark}

\subsection{The $A_\infty$-structure}
We will now equip $CW^*(L;H)$ with the structure of an $A_\infty$-algebra, generalizing the differential $\mu^1$ introduced above. For each $(d,F,\bfw)$ we will in fact define a multilinear map
\begin{equation} \label{eq:mu-f}
CF^*(L;w^d H)[q] \otimes \cdots \otimes CF^*(L;w^1 H)[q]
\xrightarrow{\mu^{d,F,\bfw}} CF^*(L;w^0 H)[q]
\end{equation}
of degree $2-d$, which encodes all the $m^{d,F,\bfw}(\bfx)$. We will first define the first summand of \eqref{eq:mu-f}, namely the one which takes values in $CF^*(L;w^0 H)$. That summand vanishes unless each entry lies in $q^{i^k}CF^*(L;w^k H)$, where by definition $i^k = 1$ if $k \in F$, and $i^k = 0$ otherwise. The remaining coefficients are
\begin{equation} \label{eq:ad-hoc-sign}
\begin{aligned}
& (-1)^{\ast} m^{d,F,\bfw}: q^{i^d} |o_{x^d}|_{\K} \otimes \cdots \otimes q^{i^1} |o_{x^1}|_{\K} \longrightarrow |o_{x^0}|_{\K}, \\
&
\ast = \textstyle \sum_j j\, \deg(x_j) \;+\;
\sum_{\substack{j \in F \\ k>j}} (\deg(x_k)-1).
\end{aligned}
\end{equation}
Having done that, there is a unique way to define the second summand, which is the one taking values in $q\,CF^*(L;w^0 H)$, such that the outcome commutes with the action of $\partial_q$ when considering reduced degrees $\bardeg(c) = \deg(c)-1$. By this we mean that
\begin{equation} \label{eq:commutes-with-dq}
\begin{aligned}
& \partial_q \mu^{d,F,\bfw}(c^d,\dots,c^1) = \\ & \qquad =
\textstyle\sum_{k=1}^d
(-1)^{\bardeg(c^{k+1}) + \cdots + \bardeg(c^{d})}
\mu^{d,F,\bfw}(c^d,\dots,\partial_q c^k,\dots,c^1).
\end{aligned}
\end{equation}
%(actually, because we are working with a fixed $\bfw$, at most one term on the right hand side can be non-vanishing).
Finally, we extend $\mu^{d,F,\bfw}$ trivially (by zero) to a multilinear map defined on the whole of $CW^*(L;H)$, and let $\mu^d$ be the sum of those maps over all $F$ and $\bfw$. For $d = 1$, this reproduces \eqref{eq:wrapped-differential} except for the last term $qb \mapsto (-1)^{\deg(b)+1} b$, which one needs to add by hand.

\begin{example}
Suppose that we have fixed isomorphisms $|o_x|_{\K} \iso \K$ for each $x$. Then $CF^*(L;wH)$ becomes a graded vector space with distinguished basis indexed by $x \in \scrX_w$, and each $m^{d,F,\bfw}(\bfx)$ is just an element of $\K$. In those terms (and omitting weights for brevity, as in Example \ref{th:algebraic-relation}) we have
\begin{equation}
\begin{aligned}
 & \mu^2(x^2,x^1) = (-1)^{\deg(x^1)}
 \textstyle \sum_{x^0} m^2(x^0,x^1,x^2)\, x^0, \\
 & \mu^2(q x^2,x^1) = (-1)^{\deg(x^1)}
 \textstyle \sum_{x^0} \big( m^2(x^0,x^1,x^2)\,q x^0
 + m^{2,\{2\}}(x^0,x^1,x^2)\,x^0 \big), \\
 & \mu^2(x^2,q x^1) = (-1)^{\deg(x^1)+\deg(x^2)-1}
 \textstyle \sum_{x^0} \big( m^2(x^0,x^1,x^2)\,q x^0 \\ & \qquad \qquad \qquad +
 m^{2,\{1\}}(x^0,x^1,x^2)\,x^0 \big),
 \\
 & \mu^2(q x^2,q x^1) = (-1)^{\deg(x^1)+\deg(x^2)-1}
 \textstyle \sum_{x^0} \big( -m^{2,\{2\}}(x^0,x^1,x^2)\, q x^0
 \\ & \qquad \qquad \qquad + m^{2,\{1\}}(x^0,x^1,x^2) \, q x^0
 + m^{2,\{1,2\}}(x^0,x^1,x^2)\, x^0 \big).
\end{aligned}
\end{equation}
\end{example}

\begin{prop}
The maps $\mu^1,\mu^2,\cdots$ satisfy the $A_\infty$-associativity equations. This means that for every $d \geq 1$ and $c^1,\dots,c^d \in CW^*(L;H)$, we have
\begin{equation} \label{eq:as}
\begin{aligned} &
\textstyle\sum (-1)^{\bardeg(c^1) + \cdots + \bardeg(c^{i-1})}
\mu^{d_+}(c^d,\dots,c^{i+d_-}, \\
& \qquad \qquad \mu^{d_-}(c^{i+d_--1},\dots,c^i),c^{i-1},\dots,c^1) = 0,
\end{aligned}
\end{equation}
where the sum is over all $d_+ + d_- = d+1$ and all $i$. Hence, $CW^*(L;H)$ equipped with these maps is an $A_\infty$-algebra.
\end{prop}

This follows from \eqref{eq:algebraic-relation} and the definition of $\mu^d$ by a direct (if somewhat tedious, because of the signs) computation. The only terms in \eqref{eq:as} which do not have direct counterparts in \eqref{eq:algebraic-relation} come from last term in $\mu^1$, but those cancel out because of \eqref{eq:commutes-with-dq}.

\section{Viterbo functoriality\label{sec:functoriality}}

This section sets up the $A_\infty$-homomorphisms $CW^*(L;H) \longrightarrow CW^*(L^{in};H^{in})$ associated to an embedding $(M^{in},L^{in}) \subset (M,L)$. In the corresponding construction for symplectic cohomology \cite{viterbo97a}, Viterbo's strategy was to shrink the interior domain $M^{in}$ by a conformal symplectic factor $\rho \ll 1$. This makes all energies inside that domain small compared to those outside, and when carried out in a carefully designed way, allows one to identify the Floer complex $CF^*(L^{in};w H^{in})$ with a quotient complex of $CF^*(L;w H)$. One takes the map on cohomology induced by this projection, which under the direct limit $w \rightarrow \infty$ yields the desired homomorphism between symplectic cohomology groups. In our case, the chain complexes for different values of $w$ are tied together by the $A_\infty$-structure. It seems that there is no single value of $\rho$ which achieves the desired properties for all moduli spaces simultaneously, and as a consequence, we can't define $\scrF$ as a projection. Instead, we let $\rho$ vary in $(0,1]$, and look at how the moduli spaces change depending on that parameter, which then yields the correction terms that need to be added to the naive projection map. The relevant formalism of parametrized moduli spaces, and its relation to the $A_\infty$-homomorphism equation, is an adaptation of that in \cite[Section 19]{fooo} (see also \cite[Section 10]{seidel04}). In the interest of simplicity and brevity, we give a full description of $\scrF^1$, but deal in a slightly more abbreviated way with the higher order components $\scrF^d$, $d \geq 2$, in particular omitting sign issues.

\subsection{Geometric setup\label{subsec:geometry-of-viterbo}}
Let $M$ be a Liouville domain, and $M^{in}$ a Liouville subdomain. This means that $M^{in} \subset M \setminus \partial M$ is a compact submanifold with boundary (of equal dimension), such that $Z$ points outwards along $\partial M^{in}$. In particular, $\theta^{in} = \theta|M^{in}$ turns $M^{in}$ itself into a Liouville domain. We call $\partial M^{in}$ the {\em dividing hypersurface}, and write $M^{out} = (M \setminus M^{in}) \cup \partial M^{in}$ for the part outside that hypersurface. By integrating $Z$ starting from the boundary, respectively from the dividing hypersurface, one gets collars, which we denote by
\begin{equation} \label{eq:2-collars}
\begin{aligned}
& \kappa: (0,1] \times \partial M \longrightarrow M, \\
& \kappa^{in}: (0,1+\epsilon) \times \partial M^{in} \longrightarrow M
\end{aligned}
\end{equation}
(for some $\epsilon>0$). In the case of $\kappa^{in}$, points $(r,y)$ with $r \leq 1$ get mapped to $M^{in}$, and those with $r \geq 1$ to $M^{out}$. Let $L \subset M$ be a Lagrangian submanifold, intersecting $\partial M$ as well as $\partial M^{in}$ transversally, and write $L^{in} = L \cap M^{in}$, $L^{out} = L \cap M^{out}$. We make the following assumption:
\begin{equation} \label{eq:relative-exactness}
\parbox{30em}{
$\theta|L^{out}$ vanishes on a neighbourhood of the boundary $\partial L^{out} = \partial L \cup \partial L^{in}$. Moreover, one can write $\theta|L = dh$, where $h|\partial L \cup \partial L^{in}$ is zero.
}
\end{equation}

This requirement merits some discussion. The first part of \eqref{eq:relative-exactness} says that $Z$ is parallel to $L^{out}$ near its boundaries. To put it more explicitly in terms of \eqref{eq:2-collars}, there is some small $\lambda>0$ such that
\begin{equation} \label{eq:piece-of-tubular}
\begin{aligned}
& \kappa^{-1}(L) \cap ([1-\lambda,1] \times \partial M) = [1-\lambda,1] \times \partial L, \\
& (\kappa^{in})^{-1}(L) \cap ([1,1+\lambda] \times \partial M^{in}) = [1,1+\lambda] \times \partial L^{in}.
\end{aligned}
\end{equation}
This is similar to \eqref{eq:stronger}, hence marginally stronger than its counterpart in \eqref{eq:asymptotically-parallel}, but that is just for technical convenience. A more interesting aspect is the condition on $h$, which amounts to a strengthened form of exactness. In more algebro-topological terms, the statement is that the relative class $[\theta|L] \in H^1(L,\partial L^{in} \cup \partial L;\R)$ is zero.

\begin{lemma} \label{th:not-relatively-exact}
Let $L_0 \subset M$ be a Lagrangian submanifold satisfying the following weaker property:
\begin{equation} \label{eq:locally-constant}
\parbox{30em}{
$\theta|L_0 = dh_0|L_0$, where $h_0 \in \smooth(M,\R)$ is a function which is locally constant
on $\partial M \cup \partial M^{in}$.
}
\end{equation}
Then there is a Hamiltonian isotopy (stationary on $\partial M \cup \partial M^{in}$) which deforms $L_0$ into a submanifold $L_1$ such that $[\theta|L_1] \in H^1(L_1,\partial L_1^{in} \cup \partial L_1;\R)$ is zero.
\end{lemma}

\proof We may assume that $\theta|L_0$ vanishes in a neighbourhood of $\partial L_0^{in} \cup \partial L_0$ (since that can be achieved by the same argument as in Lemma \ref{th:bend-boundary}). In that case, we may further assume that $h_0$ is locally constant on a neighbourhood of $\partial M^{in} \cup \partial M \subset M$. Let $Y$ be the Hamiltonian vector field of $h_0$. Then $L_Y\theta$ is a closed one-form vanishing near $\partial M^{in} \cup \partial M$. The relative cohomology class of this form is
\begin{equation}
 [L_Y\theta] = [i_Y \omega] = [-dh_0] \in H^1(M,\partial M \cup \partial M^{in};\R).
\end{equation}
Integrating this out, we find that if $\xi^t$ is the flow of $Y$, then $[(\xi^1)^*\theta] = [\theta - dh_0]$ becomes zero when restricted to $H^1(L_0,\partial L_0 \cup \partial L_0^{in};\R)$, which of course means that $[\theta]$ itself has the same property with respect to $L_1 = \xi^1(L_0)$. This provides the desired isotopy. \qed

Note that unlike Lemma \ref{th:bend-boundary}, the isotopy constructed here is not necessarily local near $\partial M \cup \partial M^{in}$. On the other hand, the following example shows that exactness of $L$, which would be vanishing of $[\theta|L] \in H^1(L;\R)$, is not sufficient to construct meaningful Viterbo functoriality maps.

\begin{example}
Take $M = D^*S^1$ to be the annulus, which is the disc cotangent bundle of $S^1$ with fibres of some radius $r$. Let $L$ be the union of two different cotangent fibres $L_i = D^*_{q_i}S^1$ ($i = 1,2$). Take the zero-section and perturb it locally near $q_1$ in a Hamiltonian way, such that it intersects $L_1$ transversally in three points. A neighbourhood of this perturbed zero-section looks like a smaller cotangent disc bundle $M^{in} = D^*S^1$ with $r^{in} \ll r$. One can arrange that the intersections $L_i^{in} = L_i \cap M^{in}$ are unions of cotangent fibres (three and one, respectively). Moreover, the canonical one-forms on our two cotangent disc bundles will satisfy $\theta|M^{in} = \theta^{in} - dh$ for some function $h$. Take $h$ and (shrinking $r^{in}$ first, if necessary) extend it to a function on the whole of $M$, which vanishes near $\partial M$. Then, if we replace the original $\theta$ by $\theta + dh$, we will actually have $\theta^{in} = \theta|M^{in}$ as in the general framework set out above.  The Lagrangian $L$ satisfies \eqref{eq:piece-of-tubular}, since $\theta|L = dh|L$ and $\theta^{in}|L^{in} = 0$. On the other hand, let $\Delta$ is the shaded region in Figure \ref{fig:annulus}, and write $\partial_l\Delta$ and $\partial_r\Delta$ for its left and right boundary arcs, and $x_u$, $x_d$ for its upper and lower corners. Then
\begin{equation}
\textstyle 0 < \int_\Delta \omega = \int_{\partial_r\Delta} \theta + \int_{\partial_l\Delta} \theta =
\int_{\partial_r\Delta} \theta +  h(x_d) - h(x_u),
\end{equation}
and by construction $\theta|\partial M^{in}>0$, which in view of the orientations occurring means that $\int_{\partial_r\Delta} \theta < 0$. Hence, the values of $h$ at the two corners are different, so that \eqref{eq:relative-exactness} as well as the weaker version \eqref{eq:locally-constant} are violated.

\begin{figure}[ht]
\begin{centering}
\begin{picture}(0,0)%
\includegraphics{annulus.pstex}%
\end{picture}%
\setlength{\unitlength}{3355sp}%
\begingroup\makeatletter\ifx\SetFigFont\undefined%
\gdef\SetFigFont#1#2#3#4#5{%
  \reset@font\fontsize{#1}{#2pt}%
  \fontfamily{#3}\fontseries{#4}\fontshape{#5}%
  \selectfont}%
\fi\endgroup%
\begin{picture}(3316,3618)(593,-2473)
\put(751,-2086){\makebox(0,0)[lb]{\smash{{\SetFigFont{10}{12.0}{\rmdefault}{\mddefault}{\updefault}{\color[rgb]{0,0,0}$M$}%
}}}}
\put(1108,-1632){\makebox(0,0)[lb]{\smash{{\SetFigFont{10}{12.0}{\rmdefault}{\mddefault}{\updefault}{\color[rgb]{0,0,0}$M^{in}$}%
}}}}
\put(2326,-2161){\makebox(0,0)[lb]{\smash{{\SetFigFont{10}{12.0}{\rmdefault}{\mddefault}{\updefault}{\color[rgb]{0,0,0}$L_2$}%
}}}}
\put(2326,989){\makebox(0,0)[lb]{\smash{{\SetFigFont{10}{12.0}{\rmdefault}{\mddefault}{\updefault}{\color[rgb]{0,0,0}$L_1$}%
}}}}
\end{picture}%
\caption{\label{fig:annulus}}
\end{centering}
\end{figure}%

In this situation, the wrapped Floer cohomologies of $L$ and $L^{in}$ are well-defined. However, as we will now show, no Viterbo-type restriction homomorphism (with reasonable properties) can exist. Write $R = \K[t,t^{-1}]$, which is the wrapped Floer cohomology of a single cotangent fibre in $D^*S^1$. The wrapped Floer cohomology of $L$ is a matrix algebra $R^{2 \times 2}$. More geometrically, we write this as
\begin{equation} \label{eq:split-11}
HW^*(L) = \begin{pmatrix} R & R \\ R & R
\end{pmatrix},
\end{equation}
where the $(j,i)$-th piece in the splitting is represented by $X$-chords going from $L_i$ to $L_j$. Similarly, the wrapped Floer cohomology of $L^{in}$ is a matrix algebra $R^{4 \times 4}$ (actually, it is now nontrivially graded, but we will not discuss that since it's irrelevant for our purpose). In terms of the pieces $L^{in}_i$ this has a block matrix decomposition
\begin{equation} \label{eq:split-13}
HW^*(L^{in}) = \begin{pmatrix} R^{3 \times 3} & R^3 \\
R^3  & R \end{pmatrix}.
\end{equation}
Suppose that we had a restriction map $HW^*(L) \rightarrow HW^*(L^{in})$ with the following properties: it is a map of rings preserving the unit elements of the ring structure and is compatible with the splittings \eqref{eq:split-11}, \eqref{eq:split-13}. In particular, multiplication with the images of $x_{12} = \left(\begin{smallmatrix} 0 & 1 \\ 0 & 0  \end{smallmatrix}\right)$ and $x_{21}  = \left(\begin{smallmatrix} 0 & 0 \\ 1 & 0 \end{smallmatrix}\right)$ on the left and right would yield a ring isomorphism between two subalgeras of $HW^*(L^{in})$ which are isomorphic to $R^{3 \times 3}$ and $R$, which is of course impossible (this example finds a more natural interpretation in the framework of wrapped Fukaya categories, where $L$ is the direct sum of two isomorphic objects, while $L^{in}$ is the sum of two non-isomorphic ones).
\end{example}

\subsection{Hamiltonian functions\label{subsec:shrinkable-h}}
We will consider functions $H \in \smooth(M,\R)$ of the following type:
\begin{equation} \label{eq:class-of-h-2}
\parbox{30em}{
Both $H$ and its restriction $H^{in} = H|M^{in}$ satisfy \eqref{eq:class-of-h}; equivalently, $H$ is everywhere positive, and both $H(\kappa(r,y)) - r$ and $H(\kappa^{in}(r,y))-r$ vanish to infinite order along $r = 1$. In addition to that, we assume that the restriction of $H - dH(Z)$ to $M^{out} \setminus \partial M^{out}$ is positive.}
\end{equation}

\begin{example} \label{th:model-function}
Suppose that we make the following ansatz:
\begin{equation}
\left\{ \begin{aligned}
& H = c^{in} \text{ on $M^{in} \setminus \kappa^{in}([1-\delta,1] \times \partial M^{in})$}, \\
& H(\kappa^{in}(r,y)) = \chi^{in}(r) \text{ for $(r,y) \in [1-\delta,1+\delta] \times \partial M^{in}$}, \\
& H = c^{out} \text{ on $M^{out} \setminus (\kappa^{in}([1,1+\delta] \times \partial M^{in})
\cup \kappa([1-\delta,1] \times \partial M))$}, \\
& H(\kappa(r,y)) = \chi^{out}(r) \text{ for $(r,y) \in [1- \delta,1] \times \partial M$}.
\end{aligned} \right.
\end{equation}
Here $\delta>0$ is small; $c^{in},c^{out}>0$ are constants; $\chi^{in} \in \smooth(\R,(0,\infty))$ is a function satisfying $\chi^{in}(r) = c^{in}$ for $r \leq 1-\delta$, $\chi^{in}(r) = c^{out}$ for $r \geq 1+\delta$, and such that $\chi^{in}(r)-r$ vanishes to infinite order at $r = 1$; similarly, $\chi^{out}$ should satisfy $\chi^{out}(r) = c^{out}$ for $r \leq 1-\delta$, and $\chi^{out}(r)$ agrees with $r$ to infinite order when  $r = 1$ (see Figure \ref{fig:hamiltonian} for a schematic picture). Then
\begin{equation}
H- dH(Z) = \left\{ \begin{aligned}
& \chi^{in} - r \cdot d\chi^{in}/dr \text{ at points $\kappa^{in}(r,y)$, $r \in (1,1+\delta]$,} \\
& c^{out} \text{ on $M^{out} \setminus (\kappa^{in}([1,1+\delta] \times \partial M^{in})
\cup \kappa([1-\delta,1] \times \partial M))$}, \\
& \chi^{out} - r \cdot d\chi^{out}/dr \text{ at points $\kappa(r,y)$, $r \in [1-\delta,1)$.}
\end{aligned} \right.
\end{equation}
It is easy to choose our functions and constants so that this is $>0$ (note the elementary calculus fact that $\chi- r \cdot d\chi/dr$ is where the tangent line to the graph of $\chi$ at $r$ hits the $y$-axis).
\end{example}

From this point onwards, we will assume that \eqref{eq:no-integer-chords} holds for both $\partial L \subset \partial M$ and $\partial L^{in} \subset \partial M^{in}$. Since by assumption, the Hamiltonian vector field $X$ of $H$ restricts to the Reeb fields $R \in \smooth(T(\partial M))$ and $R^{in}\in \smooth(T(\partial M^{in}))$, it follows that these two hypersurfaces contain no integer $X$-chords. We can then can split the set $\scrX_w$ of $X$-chords of length $w$ into two parts
\begin{equation} \label{eq:divide-x}
\scrX_w = \scrX_w^{in} \cup \scrX_w^{out}
\end{equation}
containing the integer $X$-chords lying in $M^{in} \setminus \partial M^{in}$ and $M^{out} \setminus \partial M^{out}$, respectively.

\subsection{Rescaling the inner part\label{subsec:rho}}
Somewhat unconventionally, we will write $\psi^\rho$ for the flow of $Z$ at time $\log(\rho)$. By definition, this shrinks all the symplectic data by a factor of $\rho$. Our strategy will be to use this flow, with $\rho \leq 1$, to shrink the inner Liouville domain to
\begin{equation} \label{eq:inside-shrunk}
M^{in,\rho} = \psi^\rho(M^{in}).
\end{equation}
By definition, $\psi^\rho(\kappa^{in}(r,y)) = \kappa^{in}(\rho r,y)$. In particular, the boundary of \eqref{eq:inside-shrunk} is $\partial M^{in,\rho} =\kappa^{in}(\{\rho\} \times \partial M^{in})$.
\begin{figure}[ht]
\begin{centering}
\begin{picture}(0,0)%
\includegraphics{hamiltonian.pstex}%
\end{picture}%
\setlength{\unitlength}{3355sp}%
\begingroup\makeatletter\ifx\SetFigFont\undefined%
\gdef\SetFigFont#1#2#3#4#5{%
  \reset@font\fontsize{#1}{#2pt}%
  \fontfamily{#3}\fontseries{#4}\fontshape{#5}%
  \selectfont}%
\fi\endgroup%
\begin{picture}(4287,6745)(226,-6194)
\put(1501,-6136){\makebox(0,0)[lb]{\smash{{\SetFigFont{10}{12.0}{\rmdefault}{\mddefault}{\updefault}{\color[rgb]{0,0,0}$\rho \leq r \leq 1$}%
}}}}
\put(3901,-2311){\makebox(0,0)[lb]{\smash{{\SetFigFont{10}{12.0}{\rmdefault}{\mddefault}{\updefault}{\color[rgb]{0,0,0}$\kappa(r,y)$,}%
}}}}
\put(1801,-2536){\makebox(0,0)[lb]{\smash{{\SetFigFont{10}{12.0}{\rmdefault}{\mddefault}{\updefault}{\color[rgb]{0,0,0}$1-\delta \leq r \leq 1+\delta$}%
}}}}
\put(1501,-1861){\makebox(0,0)[lb]{\smash{{\SetFigFont{10}{12.0}{\rmdefault}{\mddefault}{\updefault}{\color[rgb]{0,0,0}$M^{in}$}%
}}}}
\put(2851,-1861){\makebox(0,0)[lb]{\smash{{\SetFigFont{10}{12.0}{\rmdefault}{\mddefault}{\updefault}{\color[rgb]{0,0,0}$M^{out}$}%
}}}}
\put(2851,-61){\makebox(0,0)[lb]{\smash{{\SetFigFont{10}{12.0}{\rmdefault}{\mddefault}{\updefault}{\color[rgb]{0,0,0}$graph(H)$}%
}}}}
\put(2851,-3661){\makebox(0,0)[lb]{\smash{{\SetFigFont{10}{12.0}{\rmdefault}{\mddefault}{\updefault}{\color[rgb]{0,0,0}$graph(H^\rho)$}%
}}}}
\put(226,-5236){\makebox(0,0)[lb]{\smash{{\SetFigFont{10}{12.0}{\rmdefault}{\mddefault}{\updefault}{\color[rgb]{0,0,0}$\rho \cdot c^{in}$}%
}}}}
\put(601,-4186){\makebox(0,0)[lb]{\smash{{\SetFigFont{10}{12.0}{\rmdefault}{\mddefault}{\updefault}{\color[rgb]{0,0,0}$1$}%
}}}}
\put(376,-286){\makebox(0,0)[lb]{\smash{{\SetFigFont{10}{12.0}{\rmdefault}{\mddefault}{\updefault}{\color[rgb]{0,0,0}$c^{out}$}%
}}}}
\put(601,-586){\makebox(0,0)[lb]{\smash{{\SetFigFont{10}{12.0}{\rmdefault}{\mddefault}{\updefault}{\color[rgb]{0,0,0}$1$}%
}}}}
\put(451,-886){\makebox(0,0)[lb]{\smash{{\SetFigFont{10}{12.0}{\rmdefault}{\mddefault}{\updefault}{\color[rgb]{0,0,0}$c^{in}$}%
}}}}
\put(3676,-2536){\makebox(0,0)[lb]{\smash{{\SetFigFont{10}{12.0}{\rmdefault}{\mddefault}{\updefault}{\color[rgb]{0,0,0}$1-\delta \leq r \leq 1$}%
}}}}
\put(2101,-2311){\makebox(0,0)[lb]{\smash{{\SetFigFont{10}{12.0}{\rmdefault}{\mddefault}{\updefault}{\color[rgb]{0,0,0}$\kappa^{in}(r,y)$,}%
}}}}
\put(1576,-5911){\makebox(0,0)[lb]{\smash{{\SetFigFont{10}{12.0}{\rmdefault}{\mddefault}{\updefault}{\color[rgb]{0,0,0}$\kappa^{in}(r,y)$,}%
}}}}
\end{picture}%
\caption{\label{fig:hamiltonian}}
\end{centering}
\end{figure}%
Similarly, we define Lagrangian submanifolds $L^{in,\rho} = \psi^\rho(L^{in})$ and complete them to submanifolds $L^\rho \subset M$ as follows:
\begin{equation}
\left\{
\begin{aligned}
& L^\rho \cap M^{in,\rho} = L^{in,\rho}, \\
& L^\rho \cap \kappa^{in}((\rho,1) \times \partial M^{in}) = (\rho,1) \times
\partial L^{in}, \\
& L^\rho \cap M^{out} = L \cap M^{out};
\end{aligned}
\right.
\end{equation}
There is a natural choice of function $h^\rho$ with $dh^\rho = \theta|L^\rho$, namely
\begin{equation} \label{eq:small-h-rho}
h^\rho = \left\{ \begin{aligned}
& \rho \cdot (h \circ \psi^{1/\rho})
\text{ on $M^{in,\rho}$}, \\
& 0 \text{ at points $\kappa^{in}(r,y)$ with $r \in (\rho,1)$,} \\
& h \text{ on $M^{out}$,}
\end{aligned} \right.
\end{equation}
which makes sense because of \eqref{eq:relative-exactness}. Finally, given any Hamiltonian function $H$ as in \eqref{eq:class-of-h-2} and any $\rho \in (0,1]$, we define the {\em rescaled function} $H^\rho$ as follows:
\begin{equation} \label{eq:h-rho}
H^\rho = \left\{ \begin{aligned}
& \rho \cdot (H \circ \psi^{1/\rho})
\text{ on $M^{in,\rho}$}, \\
& r \text{ at points $\kappa^{in}(r,y)$ with $r \in (\rho,1)$,} \\
& H \text{ on $M^{out}$.}
\end{aligned} \right.
\end{equation}

This is smooth since, along $\partial M^{in,\rho}$, the function $\rho \cdot (H \circ \psi^{1/\rho})$ agrees to infinite order with $\rho \cdot (r \circ \psi^{1/\rho}) = \rho \cdot \rho^{-1} \cdot r = r$ (Figure \ref{fig:hamiltonian} shows what happens if we apply this rescaling process to one of the functions from Example \ref{th:model-function}). Correspondingly, the associated Hamiltonian vector field is
\begin{equation} \label{eq:scaled-x}
X^\rho = \left\{ \begin{aligned}
& \psi^\rho_*(X)
\text{ on $M^{in,\rho}$}, \\
& \kappa_*^{in}(0,R^{in}) \text{ at points $\kappa^{in}(r,y)$ with $r \in (\rho,1)$,} \\
& X \text{ on $M^{out}$.}
\end{aligned} \right.
\end{equation}
There is a natural bijection between $X$-chords with boundary in $L$ and $X^\rho$-chords with boundary in $L^\rho$, which we denote by
\begin{equation} \label{eq:shrink-x}
\begin{aligned}
\scrX_w & \stackrel{\iso}{\longrightarrow} \scrX^{\rho}_w, \\
x & \longmapsto x^\rho.
\end{aligned}
\end{equation}
For chords lying in $M^{out}$, the map is a simple equality $x^\rho = x$, while for the ones in $M^{in}$ it is $x^\rho = \psi^\rho(x)$. As usual, we end our discussion by fixing once and for all a Hamiltonian $H$ as in \eqref{eq:class-of-h-2}, and which in addition should satisfy \eqref{eq:nondegenerate-x} and \eqref{eq:nonconcatenate-x}. These conditions are still generic for $\dim(M) \geq 4$, even within the present smaller class of functions (one can see that by inspecting the argument from Section \ref{subsec:chords}).

We also need to relate our choices of almost complex structures on $M$ and $M^{in}$, but here the connection will be made in a somewhat looser sense. Suppose we are given $I \in \scrJ(M)$ and $I^{in} \in \scrJ(M^{in})$. Then, an {\em interpolating family} of almost complex structures is a family $I^\rho \in \scrJ(M)$ parametrized by $\rho \in (0,1]$, such that $I^1 = I$ and:
\begin{equation} \label{eq:interpolate}
\parbox{30em}{
There is a $\delta>0$ such that for $\rho \leq \delta$, the following properties hold: $I^\rho|M^{in,\rho} = \psi^\rho_*I^{in}$, and $(\kappa^{in})^*(I^\rho)$ is of contact type on
some neighbourhood of $\{\rho\} \times \partial M^{in} \subset [\rho,1] \times \partial M^{in}$.
}
\end{equation}
Note that the first part already implies that $I^\rho$ is of contact type along the boundary of $M^{in,\rho}$; the second part then extends that control slightly to a piece of the cone on the outside of that boundary. The same idea applies to infinitesimal deformations. Namely, assume that we have $I$, $I^{in}$, $I^\rho$ as well as infinitesimal deformations $K \in T\scrJ(M)_I$, $K^{in} \in T\scrJ(M^{in})_{I^{in}}$. Then, an interpolating family consists of infinitesimal deformations $K^\rho$ of $I^\rho$, satisfying the analogue of \eqref{eq:interpolate}. Obviously, this leads to the same relationship between the perturbed almost complex structures which we saw in Equation \eqref{eq:perturbed}.

\begin{remark}
Usually, we will be dealing with families of almost complex structures parametrized by points on a Riemann surface, and possibly additional moduli. In that case, we always assume that the constant $\delta$ and the neighbourhood of $\{\rho\} \times \partial M^{in}$ in \eqref{eq:interpolate} should be the same throughout the family.
\end{remark}

\subsection{Parametrized moduli spaces\label{subsec:parametrized}}
Consider the moduli spaces \eqref{eq:moduli-space} for the given $L \subset M$. In addition to various data which concerns only the domains (strip-like ends, popsicle sticks, one-forms), the construction of those spaces involves choices of almost complex structures on $M$. Namely, for each $w$ we have a family $\bfJ_w$ parametrized by $[0,1]$, and for every stable weighted popsicle we have another family $\bfJ_{S,\bfphi,\bfw}$ parametrized by points of $S$. Choose analogous families $\bfJ^{in}_w$ and $\bfJ_{S,\bfphi,\bfw}^{in}$ on $M^{in}$, and introduce the corresponding moduli spaces, which will be denoted by
\begin{equation} \label{eq:in-spaces}
\scrR^{d+1,\bfp,\bfw}(\bfx)^{in} \quad \text{for $\bfx = (x^0,\dots,x^d)$, $x^k \in \scrX_{w^k}^{in}$.}
\end{equation}
To relate the two constructions, we will then need interpolating families $\bfJ^\rho_w$ and $\bfJ^\rho_{S,\bfphi,\bfw}$, in the sense defined above. To ensure that the interpolating families are asymptotically consistent in a suitable sense, one needs to copy in more detail the construction from Section \ref{subsec:j}, by which we mean the following. Recall that $\bfJ_{S,\bfphi,\bfw}$ was obtained by taking an initial choice $\bfI_{S,\bfphi,\bfw}$, which satisfied consistency in a strict sense, and then perturbing it by exponentiating an infinitesimal deformation $\bfK_{S,\bfphi,\bfw}$. The same applies to the corresponding data on $M^{in}$. Similarly, one should first choose interpolating families $\bfI_{S,\bfphi,\bfw}^\rho$ for the initial choices of data, satisfying strict consistency, and then perturb them using a suitable $\bfK_{S,\bfphi,\bfw}^\rho$ which is asymptotically consistent. We omit the details, which are straightforward.

For any $\rho \in (0,1]$, define $\scrR^{d+1,\bfp,\bfw}(\bfx)^\rho$ to be the moduli space of popsicle maps analogous to \eqref{eq:moduli-space}, but where $L$ is replaced by $L^\rho$, $X$ by $X^\rho$, the limits of the map $u$ are the $x^{k,\rho}$, and the almost complex structures belong to the interpolating families we have just constructed. For $\rho = 1$ this just reproduces the original moduli space $\scrR^{d+1,\bfp,\bfw}(\bfx)$. On the other hand, if $\rho$ is sufficiently small and the $x^k$ all lie in $\scrX_{w^k}^{in}$, it contains a copy of $\scrR^{d+1,\bfp,\bfw}(\bfx)^{in}$, obtained by rescaling the map, $u = \psi^\rho \circ u^{in}$. In the converse direction we have the following statement, proved in Lemma \ref{th:contain}:
\begin{equation} \label{eq:sft}
\parbox{30em}{
There is a positive integer $\nu$, such that for any fixed $(d,\bfp,\bfw)$ with $w^0,\dots,w^d \geq \nu$, the following holds. There is a $\delta$ such that  all points in $\scrR^{d+1,\bfp,\bfw}(\bfx)^\rho$ with $x^0 \in \scrX_{w^0}^{in}$ and $\rho \leq \delta$ are of the form $(S,\bfphi,u = \psi^\rho \circ u^{in})$ for some $(S,\bfphi,u^{in}) \in \scrR^{d+1,\bfp,\bfw}(\bfx)^{in}$. In particular, this means that $\scrR^{d+1,\bfp,\bfw}(\bfx)^\rho = \emptyset$ unless all the other $x^k$, $k>0$, lie in $\scrX_{w^k}^{in}$ as well.
}
\end{equation}
Convention \ref{th:waffle} applies here, meaning that while the statement is formulated in terms which make sense for stable popsicles $(S,\bfphi)$, the analogous result holds for $d+|F| = 1$ as well. We will use this property heavily from now on, hence
\begin{equation} \label{eq:minimal_weight}
\parbox{30em}{
we require throughout the subsequent discussion that all weights should be $\geq \nu$.
}
\end{equation}
The first step is to suitably package our moduli spaces for varying values of $\rho$. The most straightforward way is to introduce {\em parametrized moduli spaces}
\begin{equation} \label{eq:parametrized-1}
\scrP^{d+1,\bfp,\bfw}(\bfx) =  \big\{(\rho,S,\bfphi,u) \;:\; \rho \in (0,1], \;
(S,\bfphi,u) \in \scrR^{d+1,\bfp,\bfw}(\bfx)^{\rho}\big\}.
\end{equation}
Adding the parameter $\rho$ modifies some analytic aspects of the theory, but only in a relatively straightforward way. For instance, if one considers the standard compactification of $\scrR^{d+1,\bfp,\bfw}(\bfx)^\rho$ for any $\rho$, the disjoint union of these compactifications provides a partial compactification of the parametrized moduli space; more precisely, this is a properification with respect to the projection map
\begin{equation} \label{eq:rho-projection}
(\rho,S,\bfphi,u) \longmapsto \rho.
\end{equation}
Next, suppose that the original almost complex structures on $M$, as well as their counterparts on $M^{in}$, have been chosen generically, so that both Theorems \ref{th:smoothness-1} and \ref{th:smoothness-2} apply. In that case, the parametrized moduli spaces are regular near their boundary points, which by definition are the points $(\rho,S,\bfphi,u)$ with $\rho = 1$. Moreover, at each such point, the derivative of the projection \eqref{eq:rho-projection} is onto, so $\rho = 1$ is a regular value of that map.

\begin{theorem} \label{th:smoothness-3}
For a generic choice of interpolating families, satisfying the required asymptotic consistency conditions, the following properties will hold. First, all moduli spaces $\scrP^{d+1,\bfp,\bfw}(\bfx)$ are regular, hence smooth manifolds with boundary of the expected dimension. Secondly, for every fixed value of $\rho$, there is at most one point $(\rho,S,\bfphi,u)$ which belongs to a zero-dimensional parametrized moduli space. Finally, if $\rho$ is such that there is one such point, then all moduli spaces $\scrR^{d+1,\bfp,\bfw}(\bfx)^\rho$ for various choices of $(\bfp,\bfw,\bfx)$  are regular everywhere except at that point.
\end{theorem}

This is the counterpart of Theorem \ref{th:smoothness-1} for parametrized moduli spaces, and is proved by essentially the same methods (see Section \ref{subsec:regular}). There is also a stratified version, which means an analogue of Theorem \ref{th:smoothness-2} applying to the isotropy strata of the $\Sym^{\bfp}$-action on the parametrized moduli spaces. We leave its formulation to the reader. From now on, assume that the interpolating families have been chosen in such a way that the conclusions of Theorem \ref{th:smoothness-3}, and of its stratified version, hold. Points of the zero-dimensional parametrized moduli spaces will be called {\em rigid}; they represent behaviour which is non-generic for a fixed value of $\rho$, but which is unavoidable in the parametrized context.

\subsection{Cascade maps (naive version)\label{subsec:naive-cascade}}
Fix $(d,\bfp)$ as well as $(\bfw,\bfx)$, where it is assumed throughout that
\begin{equation} \label{eq:0-in}
x^0 \in \scrX_{w^0}^{in}.
\end{equation}
Let $(T,\bfF)$, $\{\bfw_v\}$, $\{\bfp_v\}$ and $\{\bfx_v\}$ be as in the definition of the Gromov compactification in Section \ref{subsec:gromov}. A stable {\em cascade map} modelled on $(T,\bfF)$ is a collection of points
\begin{equation} \label{eq:cascade-map}
(\rho_v,S_v,\bfphi_v,u_v) \in \scrP^{|v|,\bfp_v,\bfw_v}(\bfx_v)
\end{equation}
satisfying the causal ordering condition \eqref{eq:time-ordering}. We denote the moduli space of cascade maps by $\scrQ^{d+1,\bfp,\bfw}(\bfx)$. Each component of this space, corresponding to a particular choice of $(T,\bfF)$ and $\bfx_v$, is an open subset of a product of parametrized moduli spaces $\scrP^{|v|,\bfp_v,\bfw_v}(\bfx_v)$, hence regular by Theorem \ref{th:smoothness-3}.
There is also a natural partial compactification $\bar\scrQ^{d+1,\bfp,\bfw}(\bfx)$, whose structure is analogous to that for moduli spaces of cascades considered in Section \ref{subsec:cascade}. Namely, in addition to $(T,\bfF)$,  assume that we have subsets $\bar{E} \subset E$ of the set of internal edges of $T$, as in \eqref{eq:5-label}. Points of the resulting stratum in the partial compactification are collections \eqref{eq:cascade-map}, with the causal ordering condition replaced by equality $\rho_{v_+} = \rho_{v_-}$ for those edges contained in $E$.

\begin{theorem} \label{th:compact-cascade}
The map $\bar\scrQ^{d+1,\bfp,\bfw}(\bfx) \rightarrow (0,1]$ which associates to any cascade
$\{(\rho_v,S_v,\bfphi_v,u_v)\}$ the value $\rho_{v_1}$, where $v_1$ is the vertex closest to the root, is proper.
\end{theorem}

This is the counterpart of \eqref{eq:properity}. The strategy of proof is as follows. Fixing $\bfx$ gives an a priori bound on the energy available for all components of the cascade. The next step is to bound the complexity of $(T,\bfF)$. Because the number of semi-infinite edges of $T$ is $d+1$, and $|F| = \sum_v |F_v|$, this comes down to bounding the length of chains of Floer trajectories. However, a standard monotonicity argument, based on the fact that integer $X^\rho$-chords never merge as $\rho$ varies, shows that each such trajectory consumes at least some fixed amount $\epsilon$ of energy, as long as the parameter values stay in some interval $[\delta,1]$. This achieves the desired bound, and parametrized Gromov compactness theory takes care of bubbling processes.

\subsection{The zero-dimensional case}
In the case where the moduli space of cascade maps is zero-dimensional, its partial compactification is regular, which means that all strata of positive codimension are actually empty. To see that, consider a point lying in such a stratum. Since the codimension is $|E|$, there is at least one edge belonging to that subset. Take the two vertices $v_\pm$ connected by that edge, with the usual convention according to which $v_+$ is supposed to lie closer to the root. For simplicity, denote the associated components by $(\rho_\pm = \rho, S_\pm,\bfphi_\pm,u_\pm) \in \scrP^{d_\pm + 1,\bfp_\pm,\bfw_\pm}(\bfx_\pm)$. For index reasons, each component must be a rigid point in its parametrized moduli space. On the other hand, Theorem \ref{th:smoothness-3} tells us that at most one rigid point exists for any fixed value of the parameter $\rho$. Hence, the two components must necessarily agree. Suppose first that these components are both Floer trajectories. We then know that their limits satisfy $x^0_- = x^1_+$ (this is part of our general condition on the $\bfx_{v}$) as well as $x^1_- = x^1_+$ (because the two components agree). But this is impossible; a Floer trajectory has non-negative energy and hence cannot have the same limits at both ends. In the remaining case where $d_\pm + |F_\pm| \geq 2$, the argument is parallel but based on weights instead of the limits. Suppose that $v_-$ is obtained by leaving $v_+$ in direction of its $k$-th flag. Then $w_-^0 = w_+^k$ (by definition of the $\bfw_v$), but also $w_-^k = w_+^k$ (because the two components are the same). But by \eqref{eq:weights} we necessarily have $w_-^0 > w_-^k$ for any $k>0$, which is again contradictory.

\begin{cor} \label{th:0-cascade}
The zero-dimensional moduli spaces $\scrQ^{d+1,\bfp,\bfw}(\bfx)$ are finite sets.
\end{cor}

To see why that is true, take a point of a zero-dimensional moduli space, and look at the component associated to the vertex $v = v_1$ nearest to the root. By \eqref{eq:0-in}, $x_v^0 = x^0 \in \scrX_{w^0}^{in}$. If $\rho_v$ was sufficiently small, \eqref{eq:sft} would ensure that this component comes from a suitable $\scrR^{in}$ moduli space, hence is regular in the unparametrized sense, which is in contradiction to its rigidity. To be a little more precise, one has to consider all the possibilities for $(|v|,\bfp_v,\bfw_v,\bfx_v)$, but since there are only finitely many, the resulting lower bound for $\rho_v$ is uniform. In view of the previous discussion, $\scrQ^{d+1,\bfp,\bfw}(\bfx)$ agrees with its partial compactification, so applying Theorem \ref{th:compact-cascade} yields the desired result.

\begin{addendum}
At this point we tweak the existing definitions in a small, but useful, way. Consider a zero-dimensional space $\scrQ^{d+1,\bfp,\bfw}(\bfx)$ where $d = 1$, $F = \emptyset$, $\bfw = \{w^0 = w, w^1 = w\}$ and $x^0 = x^1 = x \in \scrX_w^{in}$. We formally add a point, denoted by $\emptyset_x$, to this space. One can consider that as the degenerate case of an empty cascade (whose underlying tree would be an infinite line, with no vertices). Of course, doing that does not affect Corollary \ref{th:0-cascade}.
\end{addendum}

\subsection{The one-dimensional case\label{subsec:1d-cascade}}
Now consider the partial compactifications of one-dimensional moduli spaces of cascade maps. Points in such a space have the property that one component belongs to a one-dimensional parametrized moduli space (we will refer to it as the {\em moveable component}), while the remaining ones are all rigid. As before, one can show that the vertices belonging to an edge in $E$ cannot both carry rigid components. In particular, if $|E| = 1$ then equality of parameter values can only happen between the moveable component and a single adjacent rigid component. In that case, another look at Theorem \ref{th:smoothness-3} shows that the moveable component is a regular point in its unparametrized moduli space, and that implies that the codimension one strata of $\bar\scrQ^{d+1,\bfp,\bfw}(\bfx)$ are regular. We now turn to the remaining case $|E| \geq 2$, where the associated strata should be empty. A slight generalization of the previous argument shows that no two rigid components can be connected by a sequence of edges in $E$ in which each step goes further away from the root. The only remaining possibility is the following one: $E$ is a set of edges starting at a fixed vertex, which is the one where the moveable component is located, and going in direction of the leaves towards other vertices; and moreover, all of those other vertices have the same (rigid) component associated to them. In this situation, which is unfortunately not ruled out by Theorem \ref{th:smoothness-3}, regularity fails. We will have to return to this issue later on.

\begin{cor} \label{th:1-cascade}
Suppose that $\scrQ^{d+1,\bfp,\bfw}(\bfx)$ is one-dimensional. Then its partial compactification $\bar\scrQ^{d+1,\bfp,\bfw}(\bfx)$ contains finitely many additional points.
\end{cor}

This is essentially the same argument as in Corollary \ref{th:0-cascade}. The only case requiring a little extra thought is when the vertex closest to the root is the one carrying the moveable component. But then, the rigid components tied to it by edges in $E$ must have the same parameter value, and if that value is small, \eqref{eq:sft} applies, leading to the same contradiction as before.

With that in mind, we can describe the ends of the partially compactified one-dimensional moduli spaces. Take a partition of $d$ into $d_{-,1} + \cdots + d_{-,l}$, and set $d_+ = l$. This should come with a decomposition of $F$ into subsets $F_+$ and $F_{-,j}$ ($j = 1,\dots,l$), and associated flavours $\bfp_+$, $\bfp_{-,j}$, such that the following holds:
\begin{equation} \label{eq:branch-1}
\parbox{30em}{
For each $f \in F$ such that $d_{-,1} + \cdots + d_{-,j-1} < p_f \leq d_{-,1} + \cdots + d_{-,j}$, either $f \in F_+$ and $p_{+,f} = j$, or $f \in F_{-,j}$ and $p_{-,j,f} = p_f - d_{-,1} - \cdots - d_{-,j-1}$.
}
\end{equation}
Choose weights $\bfw_+$, $\bfw_{-,j}$ and collections of chords $\bfx_+$, $\bfx_-$ as follows:
\begin{equation} \label{eq:branch-2} \left\{
\begin{aligned} &
x_{-,j}^k = x^{d_{-,1} + \cdots + d_{-,j-1} + k}, \;
w_{-,j}^k = w^{d_{-,1} + \cdots + d_{-,j-1} + k} \;
 \text{ for $k>0$; } \\ &
x_{-,j}^0 = x_+^j, \; w_{-,j}^0 = w_+^j = w_{-,j}^1 + \cdots + w_{-,j}^{d_{-,j}} + |F_{-,j}|; \\ &
x_+^0 = x^0, \; w_+^0 = w^0.
\end{aligned} \right.
\end{equation}

\begin{cor} \label{th:ends}
Suppose that $\scrQ^{d+1,\bfp,\bfw}(\bfx)$ is one-dimensional. Then, its partial compactification has finitely many ends, and these are modelled on the disjoint union of products
\begin{equation} \label{eq:ends}
(0,\delta] \times \scrR^{d_+ + 1,\bfp_+,\bfw_+}(\bfx_+)^{in} \times
\prod_j \scrQ^{d_{-,j} + 1,\bfp_{-,j},\bfw_{-,j}}(\bfx_{-,j}).
\end{equation}
The combinatorial data is as described in \eqref{eq:branch-1}, \eqref{eq:branch-2}, with the degrees of the $x^+_j$ chosen so that all the moduli space factors in \eqref{eq:ends} are zero-dimensional.
\end{cor}
\begin{figure}[ht]
\begin{centering}
\begin{picture}(0,0)%
\includegraphics{running.pstex}%
\end{picture}%
\setlength{\unitlength}{3355sp}%
\begingroup\makeatletter\ifx\SetFigFont\undefined%
\gdef\SetFigFont#1#2#3#4#5{%
  \reset@font\fontsize{#1}{#2pt}%
  \fontfamily{#3}\fontseries{#4}\fontshape{#5}%
  \selectfont}%
\fi\endgroup%
\begin{picture}(4011,4037)(1126,-1670)
\put(1166,1876){\makebox(0,0)[lb]{\smash{{\SetFigFont{10}{12.0}{\rmdefault}{\mddefault}{\updefault}{\color[rgb]{0,0,0}$x^5$}%
}}}}
\put(2830,-1621){\makebox(0,0)[lb]{\smash{{\SetFigFont{10}{12.0}{\rmdefault}{\mddefault}{\updefault}{\color[rgb]{0,0,0}$x^0$}%
}}}}
\put(2372,-545){\makebox(0,0)[lb]{\smash{{\SetFigFont{10}{12.0}{\rmdefault}{\mddefault}{\updefault}{\color[rgb]{0,0,0}$\scrR^{3,\cdots}(\cdots)^{in}$}%
}}}}
\put(1645,1338){\makebox(0,0)[lb]{\smash{{\SetFigFont{10}{12.0}{\rmdefault}{\mddefault}{\updefault}{\color[rgb]{0,0,0}$\scrQ^{3,\cdots}(\cdots)$}%
}}}}
\put(3441,1364){\makebox(0,0)[lb]{\smash{{\SetFigFont{10}{12.0}{\rmdefault}{\mddefault}{\updefault}{\color[rgb]{0,0,0}$\scrQ^{4,\cdots}(\cdots)$}%
}}}}
\put(4676,1589){\makebox(0,0)[lb]{\smash{{\SetFigFont{10}{12.0}{\rmdefault}{\mddefault}{\updefault}{\color[rgb]{0,0,0}$x^1$}%
}}}}
\put(4401,2055){\makebox(0,0)[lb]{\smash{{\SetFigFont{10}{12.0}{\rmdefault}{\mddefault}{\updefault}{\color[rgb]{0,0,0}$x^2$}%
}}}}
\put(3408,2235){\makebox(0,0)[lb]{\smash{{\SetFigFont{10}{12.0}{\rmdefault}{\mddefault}{\updefault}{\color[rgb]{0,0,0}$x^3$}%
}}}}
\put(2471,2145){\makebox(0,0)[lb]{\smash{{\SetFigFont{10}{12.0}{\rmdefault}{\mddefault}{\updefault}{\color[rgb]{0,0,0}$x^4$}%
}}}}
\end{picture}%
\caption{\label{fig:running}}
\end{centering}
\end{figure}%

As usual, the situation is clearer in graphical terms (see Figure \ref{fig:running}). Each end \eqref{eq:ends} consists of cascades where the vertex nearest to the root carries the moveable component. For that moveable component, the parameter can become arbitrarily small, and then \eqref{eq:sft} applies, identifying it with a point in $\scrR^{in}$. Note that it is crucial to include the additional points $\emptyset_x$ in the $\scrQ$ factors in \eqref{eq:ends}, since they cover the cases where some of the positive ends of the moveable components do not have additional components attached to them (in the simplest case there is only one component, which must be moveable).

\subsection{The linear restriction homomorphism\label{subsec:chain}}
Temporarily, let's concentrate on the moduli spaces of cascade maps with $d = 1$. This gets rid of the transversality issue described in the discussion preceding Theorem \ref{th:1-cascade} (since that issue assumed the existence of a vertex with valency $|v|-1 \geq |E| \geq 2$). Of course, the situation also simplifies in other respects, which makes it possible to shed some notational ballast. The maps $\bfp$ are constantly equal to $1$, so we replace them in the notation everywhere by $F$. The trees $T$ reduce to chains of two-valent vertices, which we number as $\{v_1,\dots,v_l\}$ starting from the one closest to the root.  We also write
\begin{equation}
\begin{aligned}
& F_j = F_{v_j} \text{ for the subsets forming the decomposition $\bfF$}, \\
& \{w_{v_j}^0 = w_{j-1}, w_{v_j}^1 = w_j\}
\text{ where } w_j = w^1 + |F_{j+1}| + \cdots + |F_l|, \\
& \{x_{v_j}^0 = x_{j-1}, x_{v_j}^1 = x_j\},
\text{ where } x_j \in \scrX_{w_j}, \; x_0 = x^0, x_l = x^1.
\end{aligned}
\end{equation}
Finally, we identify each surface $S_{v_j}$ with $Z$. With that in mind, the components of a cascade map in our moduli space are of the form
\begin{equation} \label{eq:linear-cascade}
(\rho_j,Z,\bfphi_j,u_j) \in \scrP^{2,F_j,\{w_{j-1},w_j\}}(x_{j-1},x_j).
\end{equation}
More precisely, the parametrized moduli space is the space of $(\rho_j,Z,\bfphi_j,u_j)$ divided by the diagonal action of $\R$ on the last two factors, reflecting the fact that the identification $S_{v_j} \iso Z$ is not unique.

Let's begin translating the geometry of the moduli spaces into algebra. For this purpose, assume that $L$ satisfies \eqref{eq:standard-floer}, and choose a grading as well as a $Pin$ structure. These structures carry over to $L^{in}$ (by restriction), and to $L^\rho$ (extending them over the obvious isotopy). To any point $\{(\rho_j,Z,\bfphi_j,u_j)\}$ in a zero-dimensional moduli space $\scrQ^{2,F,\bfw}(\bfx)$ one can associate a map
\begin{equation} \label{eq:chain-orientation-map}
|o_{\{(\rho_j,Z,\bfphi_j,u_j)\}}^{red}|_{\K} : |o_{x^1}|_\K \longrightarrow |o_{x^0}|_{\K}
\end{equation}
analogous to \eqref{eq:orientation-map}. The orientation issues underlying this definition are discussed in Section \ref{subsec:homotopy-signs}. In the case of the formal point $\emptyset_x$ of the moduli space, we set \eqref{eq:chain-orientation-map} to be the identity on $|o_x|_{\K}$. Let $q^{1,F,\bfw}(\bfx)$ be the sum of \eqref{eq:chain-orientation-map} over all points in the zero-dimensional moduli space. A symmetry argument parallel to Lemma \ref{th:vanishing} (see again Section \ref{subsec:homotopy-signs} for more details) shows that

\begin{lemma} \label{th:permute-cascade}
If $|F|>1$ then $q^{1,F,\bfw}(\bfx)$ vanishes.
\end{lemma}

Of the two remaining cases, consider first that where $F = \emptyset$. In that situation, we omit $F$ from the notation for moduli spaces and the resulting maps between orientation spaces, following the model of \eqref{eq:moduli-space}. Consider a partially compactified one-dimensional space $\bar\scrQ^{2,\bfw}(\bfx)$, where $\bfw = \{w^0 = w,w^1 =w\}$ and $\bfx = \{x^0,x^1\}$. The ends of such spaces were described in Corollary \ref{th:ends}. In this specific situation, they are modelled on
\begin{equation} \label{eq:0-boundary}
(0,\delta] \times \scrR^{2,\{w,w\}}(x^0,x^{\mathit{new}})^{in} \times \scrQ^{2,\{w,w\}}(x^{\mathit{new}},x^1),
\end{equation}
where $x^{\mathit{new}} \in \scrX^{in}_w$ satisfies $\deg(x^{\mathit{new}}) = \deg(x^0)-1 = \deg(x^1)$. On the other hand, we have three kinds of boundary points. The first kind are those already contained in $\scrQ^{2,\bfw}(\bfx)$ itself, namely those where the last component is moveable, and reaches time $\rho_l = 1$. The set of such boundary points is the union of the spaces
\begin{equation} \label{eq:1-boundary}
\scrQ^{2,\{w,w\}}(x^0,x^{\mathit{new}}) \times \scrR^{2,\{w,w\}}(x^{\mathit{new}},x^1)
\end{equation}
over all $x^{\mathit{new}} \in \scrX_w$ with $\deg(x^{\mathit{new}}) = \deg(x^0) = \deg(x^1)+1$. As in the case of \eqref{eq:ends} and \eqref{eq:0-boundary}, this works because the left hand factor may contain a point $\emptyset_x$, which represents the case where the boundary point is a chain of total length $l = 1$. The remaining kinds of boundary points lie in the partial compactification, where the set $E$ in \eqref{eq:5-label} consists of a single edge, connecting $v_j$ to $v_{j+1}$ for some $j$, and where $\bar{E}$ is either $\emptyset$ or $E$. It turns out that for each collection $\{(\rho_j,Z,\bfphi_j,u_j)\}$ satisfying $\rho_j = \rho_{j+1}$, the contributions coming from the two choices of $\bar{E}$ cancel each other, so that the total contribution is zero. In other words, one could glue together these boundary points in pairs, in a way compatible with the orientations of the one-dimensional moduli spaces induced from isomorphisms $o_{x^0} \iso \R$ and $o_{x^1} \iso \R$. Hence, the remaining relation is the one which counts points of the form \eqref{eq:0-boundary} and \eqref{eq:1-boundary}. Taking into account signs, which will be elaborated on in Section \ref{subsec:homotopy-signs}, one finds that
\begin{equation} \label{eq:q0-relation}
\begin{aligned}
 & \textstyle\sum m^{1,\{w,w\},in}(x^0,x^{\mathit{new}}) \circ q^{1,\{w,w\}}(x^{\mathit{new}},x^1) \\ & -
 \textstyle \sum q^{1,\{w,w\}}(x^0,x^{\mathit{new}}) \circ m^{1,\{w,w\}}(x^{\mathit{new}},x^1) = 0,
\end{aligned}
\end{equation}
where $m$ and $m^{in}$ are the maps defined in Section \ref{subsec:algebraic-relation} on $M$ and $M^{in}$, respectively, and the sum is over all $x^{\mathit{new}}$ with the appropriate Maslov index. The entire discussion carries over to the remaining case $F = \{1\}$, $\bfw = \{w^0 = w+1,w^1 =w\}$ with minor modifications. In considering the analogue of \eqref{eq:0-boundary}, there are two possibilities, depending on whether the component carrying the popsicle is the first one or not. This leads to two types of ends,
\begin{equation} \label{eq:0-boundary-again}
\begin{aligned}
& (0,\delta] \times \scrR^{2,\{w+1,w+1\}}(x^0,x^{\mathit{new}})^{in} \times
\scrQ^{2,\{1\},\{w+1,w\}}(x^{\mathit{new}},x^1), \\
& (0,\delta] \times \scrR^{2,\{1\},\{w+1,w\}}(x^0,x^{\mathit{new}})^{in} \times
\scrQ^{2,\{w,w\}}(x^{\mathit{new}},x^1).
\end{aligned}
\end{equation}
Similarly, there are two possibilities for the boundary points corresponding to \eqref{eq:1-boundary},
\begin{equation} \label{eq:1-boundary-again}
\begin{aligned}
& \scrQ^{2,\{w+1,w+1\}}(x^0,x^{\mathit{new}}) \times \scrR^{2,\{1\},\{w,w+1\}}(x^{\mathit{new}},x^1), \\
& \scrQ^{2,\{1\},\{w+1,w\}}(x^0,x^{\mathit{new}}) \times \scrR^{2,\{w,w\}}(x^{\mathit{new}},x^1).
\end{aligned}
\end{equation}
The upshot are the algebraic relations
\begin{equation} \label{eq:q1-relation}
\begin{aligned}
 & -\textstyle \sum m^{1,\{w+1,w+1\},in}(x^0,x^{\mathit{new}}) \circ q^{1,\{1\},\{w+1,w\}}(x^{\mathit{new}},x^1) \\ & + \textstyle\sum
 m^{1,\{1\},\{w+1,w\},in}(x^0,x^{\mathit{new}}) \circ q^{1,\{w,w\}}(x^{\mathit{new}},x^1) \\
 & - \textstyle\sum
 q^{1,\{w+1,w+1\}}(x^0,x^{\mathit{new}}) \circ m^{1,\{1\},\{w+1,w\}}(x^{\mathit{new}},x^1) \\ & -
 \textstyle \sum q^{1,\{1\},\{w+1,w\}}(x^0,x^{\mathit{new}}) \circ m^{1,\{w,w\}}(x^{\mathit{new}},x^1)
= 0.
\end{aligned}
\end{equation}

For each $w \geq \nu$, let $\gamma: CF^*(L;wH) \rightarrow CF^*(L^{in};wH^{in})$ and $\lambda: CF^*(L;wH) \rightarrow CF^*(L^{in};(w+1)H^{in})[-1]$ be the maps whose matrix coefficients are $q^{1,\{w,w\}}(x^0,x^1)$ and $q^{1,\{1\},\{w+1,w\}}(x^0,x^1)$, respectively. We want to combine them into a single $\partial_q$-linear map between suitable wrapped cochain groups. More precisely, if $C^\nu$ is the subcomplex from Lemma \ref{th:partially-forget}, we define
\begin{equation} \label{eq:f1}
\begin{aligned}
& \scrF^1: C^\nu \longrightarrow CW^*(L^{in};H^{in}), \\
& \scrF^1(a + q b) = \gamma(a) + q \gamma(b) + \lambda(b).
\end{aligned}
\end{equation}
From \eqref{eq:q0-relation}, \eqref{eq:q1-relation} and the definition of the wrapped differential \eqref{eq:wrapped-differential}, it follows that $\scrF^1$ is a chain homomorphism.

\begin{remark} \label{th:corrected-projection}
The decomposition \eqref{eq:divide-x} induces a projection $CW^*(L;H) \rightarrow CW^*(L^{in};H^{in})$, whose restriction to $C^\nu$ is precisely the summand in $\scrF^1$ coming from the formal points $\emptyset_x$ in our moduli spaces. One can therefore think of $\scrF^1$ as the projection map plus ``instanton corrections'' coming from rigid points in the parametrized moduli spaces $\scrP^{2,F,\bfw}(\bfx)$ with $|F| \leq 1$. Of course, the projection by itself would generally not be a chain map.
\end{remark}

\subsection{Almost complex structures on cascades}
In the original definition, moduli spaces of cascade maps were (open subsets of) products of parametrized moduli spaces. There is no intrinsic need for maintaining this product structure, as long as the necessary relations between codimension one boundary faces are preserved. This allows us to enlarge the space of perturbations, which solves the previously encountered transversality issues.

Fix $d,\bfp,\bfw$, and let $\bar\scrQ^{d+1,\bfp,\bfw}$ be the partially compactified moduli space of weighted cascades. As in the case of popsicles, this is just a copy of the space $\bar\scrQ^{d+1,\bfp}$ defined in Section \ref{subsec:cascade}, but where we wish to emphasize that on each $(T,\bfF,E,\bar{E})$ stratum, the vertices of the tree come with weights $\bfw_v$ as well as flavours $\bfp_v$. Associated to any point in this stratum, represented by a collection $\{(\rho_v,S_v,\bfphi_v)\}$, we want to have families of almost complex structures of the following kind:
\begin{align} \label{eq:edge-j}
& \parbox{30em}{
For each edge $e$ of $T$, a family $\bfJ_e$ of almost complex structures. If the edge is finite and connects two vertices $v_\pm$, this family is parametrized by $(\rho,t) \in [\rho_{v_+},\rho_{v_-}] \times [0,1]$. If the edge is semi-infinite and connects $v$ to one of the leaves, $\bfJ_e$ is parametrized by $[\rho_v,1] \times [0,1]$; similarly, if it connects the root to $v$, $\bfJ_e$ is parametrized by $(0,\rho_v] \times [0,1]$.
}\\[0.5em] & \label{eq:vertex-j}
\parbox{30em}{
For each vertex $v$ of the tree, we want to have a family $\bfJ_v$ of almost complex structures parametrized by $S_v$, which is asymptotically compatible with the strip-like ends (in the usual sense: as $s \rightarrow \pm\infty$, $J_{v,\epsilon^k(s,t)}$ converges to a family depending only on $t$, and this convergence is faster than exponential in any $C^r$ norm).
}
\end{align}
It is crucial that $\bfJ_e$ and $\bfJ_v$ may depend on the whole collection $\{(\rho_v,S_v,\bfphi_v)\}$, and not just on the part of the data associated to that particular vertex or edge. This should really be reflected in the notation, but we avoid doing that in order to keep things reasonably brief. It is implicit in our definition that the dependency on the points in the moduli space of cascades is smooth. In fact, as one goes to the boundary strata in connected components \eqref{eq:delta-simplex} where some popsicles $(S_v,\bfphi_v)$ degenerate into broken popsicles, one wants the same convergence behaviour as in Remark \ref{th:concrete-consistency}. Moreover, we require all our choices to be invariant under the action of $\Sym^{\bfp}$.

The additional requirements on \eqref{eq:edge-j}, \eqref{eq:vertex-j} are of four kinds. The {\em large parameter value} restriction says the following. Let $e$ be an edge, and suppose that the weight on the flags associated to it is $w$. Then, the restriction of $\bfJ_e$ to $\rho = 1$, if it is nonempty, agrees with the family $\bfJ_w$ previously chosen as part of the construction of wrapped Floer cohomology of $L$. Similarly, suppose that $\rho_v = 1$ for some vertex $v$. Then $\bfJ_v = \bfJ_{S_v,\bfphi_v,w_v}$ should agree with the almost complex structure on that particular stable popsicle. Finally, note that the vertices where $\rho_v = 1$ necessarily form a finite union of subtrees inside $T$, each of which reaches up into the leaves. If we remove those subtrees, the result is another tree $T'$, which comes with its own data $(\bfp',\bfF',E',\bar{E}')$. We then require that our families of almost complex structures, when restricted to the edges and vertices of $T'$, should depend only on that part of the cascade, and should in fact agree with the ones the cascade carries as part of our universal choice.

Secondly, we have a {\em small parameter value} requirement, whose structure is roughly similar to the previous one. Namely, there should be a $\delta>0$ (depending on $d,\bfp,\bfw$) such that the restriction of $\bfJ_e$ to any $\rho \leq \delta$, if nonempty, has properties analogous to those for an interpolating family. Namely, on $M^{in,\rho}$ it should agree with $\psi_*^{\rho}\bfJ_w^{in}$, as previously chosen when defining the wrapped Floer cohomology of $L^{in}$. Moreover, in some neighbourhood of $\partial M^{in,\rho} \subset \kappa^{in}([\rho,1] \times \partial M^{in})$ it should be of contact type. Similarly, if $\rho_v \leq \delta$ for some vertex $v$, then $\bfJ_v|M^{in,\rho}$ is the pushforward of the family $\bfJ^{in}_{S_v,\bfphi_v,w_v}$ associated to that stable popsicle for $M^{in}$, and we have the same local condition near $\partial M^{in,\rho}$ as before.  Slightly more interestingly, the vertices where $\rho_v \leq \delta$ form a subtree of $T$, reaching down to the root. Removing $T$ yields a finite disjoint union of similar connected trees $T'_k$, and we require that on the part of our cascade corresponding to each of these, the almost complex structures should be those given by considering only the subtree.

Third, there is an {\em unstable consistency} condition. Pick a vertex $v$, a number $k \in \{0,\dots,|v|-1\}$, and let $e$ be the edge belonging to the $k$-th flag adjacent to $v$. Then the limit of $\bfJ_v$ over the $k$-th end must agree with the family $\bfJ_e$ specialized to $\rho = \rho_v$. Finally, we have {\em stable consistency}, which says that all the almost complex structures depend only on $E$, and are independent of $\bar{E}$. In particular, if we consider pairs of isomorphic codimension one boundary faces defined by $|E| = 1$ and $\bar{E} = \emptyset$ or $E$, then the almost complex structures on both copies of such a face must coincide.

One can build families of almost complex structures satisfying these conditions by gradually increasing the complexity, say by induction on $d+|F|$. For any given $(T,\bfF)$, one first considers the large and small parameter value requirements. In those situations, the almost complex structures are completely determined by those for smaller trees, which had been chosen in a previous step. Note that even though both situations can apply simultaneously, the resulting combined condition is non-contradictory. Having done that, one then extends the choices over the entire moduli space, while fulfilling the relatively weak consistency conditions over the boundary strata, and finally averages to obtain $\Sym^{\bfp}$-invariance.

\subsection{Cascade maps (corrected versions)}
Suppose now that \eqref{eq:edge-j}, \eqref{eq:vertex-j} have been chosen, satisfying all the requirements stated above. Fix $(d,\bfp,\bfw)$. Consider pairs $(T,\bfF)$ as before, but where now two-valent vertices with empty subsets $F_v = \emptyset$ are allowed, violating the previous stability condition. A {\em cascade map} modelled on $(T,\bfF)$ is a collection $\{(\rho_v,S_v,\bfphi_v,u_v)\}$, where the maps $u_v$ have to satisfy certain $\bar\partial$-equations, which we will now specify.

Suppose that we erase each unstable vertex of $T$, directly sewing together the two edges which meet there, and correspondingly forget all components $(\rho_v,S_v,\bfphi_v)$ associated to such vertices. The outcome is a stabilized tree, and a stable cascade modelled on it, which defines a point of $\scrQ^{d+1,\bfp,\bfw}$. Take the almost complex structures \eqref{eq:edge-j}, \eqref{eq:vertex-j} associated to that cascade. Consider first a stable vertex $v$ of $T$, which survives to a vertex of the stabilized tree. In that case, we require that $u_v$ should be a solution of \eqref{eq:dbar} for the almost complex structure $\bfJ_v$ belonging to that vertex. On the other hand, an unstable vertex $v$ gets mapped to an edge of the stabilized tree, and in that case $u_v$ should be a Floer trajectory for the almost complex structure $\bfJ_e$ associated to that edge. As usual, two such trajectories which differ by a translation are considered to be the same.

From now on, $\scrQ^{d+1,\bfp,\bfw}(\bfx)$ denotes the moduli space of cascade maps in the modified sense we have just explained. Essentially by construction, the definition and structure of the partial compactification $\bar\scrQ^{d+1,\bfp,\bfw}(\bfx)$ developed in Section \ref{subsec:naive-cascade} carries over, and so does the transversality theory in so far as it yielded positive results. The breakdown in transversality for one-dimensional partially compactified moduli spaces, discussed in Section \ref{subsec:1d-cascade}, is resolved for the following simple reason. Take two vertices $v_{1,-}$ and $v_{2,-}$ of $T$, both of which are connected to the same vertex $v_+$ by edges pointing towards the root. Suppose first that the $v_{k,-}$ are stable ($v_+$ is automatically stable, since it is at least trivalent by assumption). By the previous discussion, we only need to consider the situation when $\rho_{v_+} = \rho_{v_-,1} = \rho_{v_-,2} = \rho$, and the weighted popsicles $(S_{v_-,k},\bfphi_{v_-,k},\bfw_{v_-,k})$ are isomorphic. In the framework of \eqref{eq:vertex-j}, the families of almost complex structures on $S_{v_-,k}$ can usually be chosen independently of each other, which means that generically, rigid solutions of the $\bar\partial$-equation will not appear for the same parameter value, thus removing the obstruction to transversality. There are two exceptions to the statement about independence. If $\rho = 1$, the almost complex structures on both $v_{-,k}$ components are necessarily the same. But in that case, transversality for the original construction of wrapped Floer cohomology of $L$ ensures that the components $\{(\rho_v,S_v,\bfphi_v,u_v)\}$ for $v = v_+,v_{k,-}$ must all be moveable, which can never happen in a one-dimensional moduli space. The other, parallel but slightly more complicated, exception occurs when $\rho \leq \delta$. Supposing that $\delta$ has been chosen sufficiently small, \eqref{eq:sft} will then apply all components of our popsicle map lying between $v_{-,k}$ and the root, which implies that these components must be moveable, leading to the same dimensional contradiction as before (this is not a circular argument: the upper bounds on parameters $\rho$ provided by Lemma \ref{th:contain} depend only on action considerations, hence are independent of the choice of almost complex structures). Finally, one also needs to consider the situation when both $\rho_{v_-,k}$ are unstable, but the argument there is essentially the same.

\begin{remark}
To summarize, we have addressed the transversality issue by enlarging the class of allowed perturbations, while preserving enough consistency to retain the basic relations between the boundary strata of one-dimensional moduli spaces. There are of course other possibilities. For instance, \cite{fooo} uses Kuranishi structures induced by virtual perturbations for the same purpose (the version of Kuranishi space theory described in \cite{joyce07} seems particularly relevant here, since we need to deal with moduli spaces for all $d$ at the same time).
\end{remark}

\subsection{The restriction homomorphism}
Assume now that generic choices of \eqref{eq:edge-j} and \eqref{eq:vertex-j} have been made, subject to all the conditions above, making the spaces $\scrQ^{d+1,\bfp,\bfw}(\bfx)$ and their partial compactifications regular. As usual we consider only the case when $\Sym^{\bfp}$ is trivial, and replace $\bfp$ by $F \subset \{1,\dots,d\}$ in the notation. Write
\begin{equation} \label{eq:q-maps}
q^{d+1,F,\bfw}(\bfx) : |o_{x^d}|_{\K} \otimes \cdots \otimes |o_{x^1}|_{\K} \longrightarrow
|o_{x^0}|_{\K}
\end{equation}
for the map obtained by counting the contributions from a zero-dimensional moduli space $\scrQ^{d+1,F,\bfw}(\bfx)$ (as announced at the beginning of the section, we won't properly discuss the signs that enter into this definition). For $d = 1$ this reproduces the definition from Section \ref{subsec:chain}, in a slightly generalized version since the choices of almost complex structures considered there form a strict subset of the ones allowed here. We assemble the \eqref{eq:q-maps} for all $\bfx$ into a multilinear map
\begin{equation} \label{eq:ff}
CF^*(L;w^d H)[q] \otimes \cdots \otimes CF^*(L;w^1 H)[q]
\xrightarrow{\scrF^{d,F,\bfw}} CF^*(L^{in};w^0 H^{in})[q]
\end{equation}
of degree $1-d$. More precisely, the first summand of \eqref{eq:ff}, which takes values in $CF^*(L^{in};w^0 H)$, has coefficients \eqref{eq:q-maps} with suitable signs, and we define the second summand to make $\scrF^{d,F,\bfw}$ a $\partial_q$-linear map in the same sense as in \eqref{eq:commutes-with-dq}. In the case where $d = 1$ and $F = \emptyset$, we also add a trivial term which is projection $CF^*(L;w H)[q] \rightarrow CF^*(L^{in};w H^{in})[q]$, compare Remark \ref{th:corrected-projection}.

Recall from \eqref{eq:minimal_weight} that all weights that we are considering are larger than a fixed weight $\nu >0$.  In particular, we extend each \eqref{eq:ff} by zero to a map defined on the whole of $(C^\nu)^{\otimes d} \rightarrow CW^*(L;H)^{\otimes d}$, and taking values in $CW^*(L^{in};H^{in})$. Define $\scrF^d$ to be the sum of all these maps for fixed $d$ but varying $(F,\bfw)$.

\begin{thm}
The maps $\cF^{d}$ satisfy the $A_\infty$-homomorphism equations
\begin{equation}
\begin{aligned}
& \sum_{d_1 + \cdots + d_k = d} \mu^{in,k}(\scrF^{d_k}(c_d,\dots,c_{d-d_k+1}), \dots,
\scrF^{d_1}(c_{d_1},\dots,c_1)) \\ & =
\sum_k (-1)^{\bardeg(c^1) + \cdots + \bardeg(c^{i-1})}
\scrF^{d_+}(c^d,\dots,c^{i+d_-}, \\[-1em]
& \qquad \qquad \qquad \qquad \qquad \mu^{d_-}(c^{i+d_--1},\dots,c^i),c^{i-1},\dots,c^1).
\end{aligned}
\end{equation}
Moreover, they can be extended to an $A_{\infty}$-homomorphism $CW^{*}(L) \to  CW^{*}(L^{in})$, in a way which is unique up to homotopy.
\end{thm}

\begin{proof}
The underlying geometry is very simple: the left hand side of the equation collects contributions from the ends of one-dimensional spaces $\bar\scrQ^{d+1,\bfp,\bfw}(\bfx)$, where the component for the vertex next the root is the moveable one, and \eqref{eq:sft} applies to it; the right hand side takes into account those boundary points where, for one of the vertices nearest to the leaves, the parameter value becomes equal to $1$. All other boundary strata cancel in pairs, just as in the corresponding argument in Section \ref{subsec:chain}. 

We know from Lemma \ref{th:partially-forget} that the inclusion $C^\nu \rightarrow CW^*(L;H)$ is a quasi-isomorphism of $A_\infty$-algebras. It is a general algebraic fact that there is always a quasi-isomorphism in inverse direction, $CW^*(L;H) \rightarrow C^\nu$, whose restriction to $C^\nu$ is the identity, and this is unique up to homotopy (see for instance \cite[Section 3.7]{keller99} and the references there). By composing the given $\scrF$ with such an inverse, one can get an $A_\infty$-homomorphism defined on the whole of $CW^*(L;H)$.
\end{proof}

\begin{remark} \label{th:modify-theta}
The reason why $\scrF$ is initially only defined on $C^\nu$ comes down to Lemma \ref{th:positive-action}. There is a simpler alternative, in the spirit of Lemmas \ref{th:bend-boundary} and \ref{th:not-relatively-exact}. Namely, one extends $h$ to a function on $M$, which vanishes near $\partial M \cup \partial M^{in}$, and then makes a replacement $\theta \mapsto \theta - dh$. One effect of this is that the $h$ terms in \eqref{eq:action} become absorbed into the $\int \theta$, which removes the previous obstruction. Of course, it simultaneously changes $Z$, hence also our family of rescaled Lagrangian submanifolds $L^\rho$.
\end{remark}

\section{Complements\label{sec:many-objects}}
Our first, and rather straightforward, task in this section is to extend the previous construction to wrapped Fukaya categories. After that, we return to the two-dimensional case, which could not be treated previously due to the failure of \eqref{eq:nonconcatenate-x} to be generic, and outline a simple workaround, which is to increase dimensions by taking the product with a standard factor (other solutions are possible; see Remark \ref{th:other}).

\subsection{Wrapped Fukaya categories}
Let $\bfL = \{L_i\}_{i \in I}$ be a collection of Lagrangian submanifolds of $M$, which is at most countable. Suppose that each of them satisfies \eqref{eq:asymptotically-parallel} as well as \eqref{eq:standard-floer}, and choose gradings as well as $Pin$ structures. An integer Reeb chord from $\partial L_{i_0}$ to $\partial L_{i_1}$ is a map $x: [0,1] \rightarrow \partial M$ satisfying $dx/dt = wR$ for some positive integer $w$, such that $x(0) \in \partial L_{i_0}$, $x(1) \in \partial L_{i_1}$. Integer $X$-chords going from $L_{i_0}$ to $L_{i_1}$ are defined in the same way. Generalizing \eqref{eq:no-integer-chords}, we assume that there are no integer Reeb chords between the boundaries of any two Lagrangian submanifolds in our collection. Having that, we choose some Hamiltonian $H$ as in \eqref{eq:class-of-h}, such that the appropriate generalizations of \eqref{eq:nondegenerate-x} and \eqref{eq:nonconcatenate-x} hold. As before, these are generic assumptions (true after a generic rescaling of $\theta$, and small isotopies of the $L_i$) provided that $\dim(M) \geq 4$. Let $CF^*(L_{i_0},L_{i_1};wH)$ be the Floer cochain complex generated by $|o_x|_{\K}$, where $x$ runs over all $X$-chords of length $w$ going from $L_{i_0}$ to $L_{i_1}$. Generalizing \eqref{eq:wrapped-cochain}, we define the wrapped complex by
\begin{equation} \label{eq:wrapped-cochain-01}
 CW^*(L_{i_0},L_{i_1}; H) = \bigoplus_{w = 1}^\infty CF^*(L_{i_0},L_{i_1}; wH)[q]
\end{equation}
The wrapped Fukaya category $\scrW(\bfL)$ has the $L_i$ as objects, \eqref{eq:wrapped-cochain-01} as morphism spaces, and an $A_\infty$-structure defined by extending the construction from Section \ref{sec:wrapped}, in an essentially straightforward way to multiple Lagrangians.

Now suppose that we have a Weinstein subdomain $M^{in} \subset M$, and that each $L_i$ satisfies \eqref{eq:relative-exactness} with a uniform lower bound on the size of the neighbourhood of the boundary of $L_i$ where the restriction of $\theta$ vanishes.  In this situation, one can define an $A_\infty$-functor $\scrW(\bfL) \rightarrow \scrW(\bfL^{in})$, where $\bfL^{in} = \{L_i^{in} = L_i \cap M^{in}\}$ by following the same procedure as before. To be slightly more precise, for every pair of elements in $I$ there will be a constant $\nu_{i_0,i_1}$ such that the analogue of Lemma \ref{th:positive-action} holds with $w \geq \nu_{i_0,i_1}$. By inspecting the proof of Lemma \ref{th:positive-action}, we find that we may choose $ \nu_{i_0, i_1} $ so that for any collection $(L_{i_0}, \cdots, L_{i_d})$, we have
\begin{equation} \nu_{i_{0}, i_d} \leq \sum_{k} \nu_{i_{k}, i_{k+1}} .\end{equation}
In particular, the collection of morphism spaces which are sums \eqref{eq:wrapped-cochain-01} restricted to $w \geq \nu_{i_0,i_1}$  forms a subcategory on which we can define a restriction functor. However, this subcategory is quasi-isomorphic to the whole thing by the same argument as before, and one can use that to formally extend the functor.

\subsection{Stabilization\label{subsec:stabilization}}
Consider $D^*S^1 = [-1,1] \times S^1$. We use standard coordinates $(p,q)$, so that the tautological one-form is $p\,dq$. Let $M$ be a Liouville domain, and $k: M \times S^1 \rightarrow \R$ a function which vanishes in a neighbourhood of $\partial M \times S^1$. One can then consider the product $M \times D^*S^1$, with the one-form
\begin{equation} \label{eq:product-form}
\tilde\theta = \theta + C\, p\,dq - dk
\end{equation}
for some constant $C$ (pullbacks by projection to various factors are implicit in this definition). The associated Liouville field is $\tilde{Z} = Z + C \, p \partial_p + Y_q + (\partial k/\partial q)  \partial_p$, where $Y_q$ is the Hamiltonian vector field of $k(\cdot,q)$ on $M$. This will point strictly outwards along both boundary faces, provided that $C \gg 0$ is big enough. Hence, after rounding off the corners suitably, we get a Liouville domain $\tilde{M} \subset M \times D^*S^1$.

Suppose that $L \subset M$ is a Lagrangian submanifold satisfying \eqref{eq:relative-exactness}. We take the product of $L$ with some cotangent fibre $D^*_aS^1 = [-1,1] \times \{a\}$. Assume that the function $k$ appearing above has been chosen in such a way that
\begin{equation} \label{eq:k-function}
dk(\cdot,a)|L = \theta|L.
\end{equation}
In that case, the Lagrangian submanifold $\tilde{L} = (L \times D^*_aS^1) \cap \tilde{M}$ satisfies $\tilde\theta|\tilde{L} = 0$, hence our previous construction of the wrapped $A_\infty$-algebra, with a suitably chosen Hamiltonian function $\tilde{H}$, applies to it (gradings and $Pin$ structures for $L$ induce the same kind of structures on $\tilde{L}$). Recall that generators of $CF^*(\tilde{L};w\tilde{H})$ correspond naturally to $\tilde{X}$-chords of length $w$, with boundary on $\tilde{L}$. The image of any such chord $\tilde{x}$ under projection $\tilde{M} \rightarrow S^1$ is a closed loop, whose winding number we denote by $\alpha(\tilde{x})$. Obviously, the moduli spaces $\scrR^{d+1,\bfp,\bfw}(\boldsymbol{\tilde{x}})$ can be nonempty only if $\alpha(\tilde{x}^0) = \alpha(\tilde{x}^1) + \cdots + \alpha(\tilde{x}^d)$. Hence, the subcomplex of $CW^*(\tilde{L};\tilde{H})$ consisting only of those generators with zero winding number is in fact an $A_\infty$-subalgebra. For $2n = \mathit{dim}(M) > 2$, one can prove that it is quasi-isomorphic to the previously defined $CW^*(L;H)$. We will not prove this here, but with this motivation in mind, one can use the subalgebra as a replacement for the wrapped complex in the two-dimensional case. It is easy to see that the Viterbo restriction maps also preserve these subalgebras.

Slightly more generally, given a finite collection $\bfL = \{L_i\}_{i \in I}$ satisfying \eqref{eq:relative-exactness}, one can take $\tilde{L}_i = (L_i \times D^*_{a_i}S^1) \cap \tilde{M}$ for pairwise distinct $a_i$, and choose $k$ in such a way that the analogue of \eqref{eq:k-function} is satisfied for all $i$. In addition, fix preimages of all the $a_i$ under the projection $\R \rightarrow \R/\Z = S^1$. For any chord $\tilde{x}$ which connects $\tilde{L}_{i_0}$ to $\tilde{L}_{i_1}$, one then has an integer $\alpha(\tilde{x})$ measuring the difference in winding numbers between the projection of $\tilde{x}$ and that of any path in $\R$ connecting the preimages of $a_{i_0},a_{i_1}$. The subcomplexes of the $CW^*(\tilde{L}_{i_0},\tilde{L}_{i_1}; \tilde{H})$ generated by chords with $\alpha(\tilde{x}) = 0$ form an $A_\infty$-subcategory of $\scrW(\tilde{\bfL})$. As before, this can be used as a replacement for $\scrW(\bfL)$ in the two-dimensional case.

\begin{remark} \label{th:other}
The stabilization construction is actually a version of the approach outlined in Remark \ref{th:modify-theta}, and correspondingly comes with some limitations. Most substantially, it only works for finitely many Lagrangian submanifolds at once. In principle, one can overcome that particular restriction by using the universal cover $D^*\R = [-1,1] \times \R$ as a stabilizing factor, but that no longer strictly falls into the previously considered framework, hence requires some additional work (in particular, in the choice of Hamiltonian functions, to avoid problems arising from the noncompactness of the $\R$ factor). Among the other possible approaches which avoid stabilization entirely, the most obvious one, since we are dealing with a transversality problem, is to use Kuranishi structures. Alternatively, along more classical lines, one could try to replace the multiples $wH$ by a more general sequence of Hamiltonian functions, or make the Hamiltonians time-dependent (which has the disadvantage of making action computations less straightforward; we have not considered this in detail).
\end{remark}

\section{Popsicles and stable maps\label{sec:lollipop}}

This section revisits the compactified moduli spaces $\bar\scrR^{d+1,\bfp}$, filling in some gaps left in the original discussion (Section \ref{sec:popsicle}). Instead of starting with the combinatorial description of these spaces as union of strata, as in \eqref{eq:compactification}, we construct them from a complex geometry perspective, based on the existing theory of stable maps (this follows the strategy used in \cite{fukaya-oh98} for Stasheff polyhedra).

\subsection{A complex analogue}
Fix $d \geq 1$ and $\bfp = \{p_f\}$, $f \in F$. Consider the projective line $\C P^1 = \C \cup \{\infty\}$, and let $Y = (\C P^1)^F$ be the product of such lines indexed by $F$, with coordinates $y_f$. For $k \in \{0,\dots,d\}$, define subsets $H^k \subset Y$ by
\begin{equation}
H^k = \begin{cases} \text{the point where all $y_f = 0$} & \text{if $k = 0$,} \\
\text{the subset where $y_f = \infty$ for all $f$ with $p_f = k$} & \text{if $k \in \{1,\dots,d\}$.}
\end{cases}
\end{equation}
A $\bfp$-flavoured {\em lollipop} is a closed genus zero Riemann surface $C$, equipped with a collection of $d+1$ distinct ordered marked points $\bfzeta = \{\zeta^0,\dots,\zeta^d\}$, together with a holomorphic map $\bfpsi: C \rightarrow Y$ of degree $(1,\dots,1)$ such that
\begin{equation} \label{eq:incidence}
 \bfpsi(\zeta^k) \in H^k \;\; \text{ for all $k$.}
\end{equation}
Equivalently, the components $\psi_f: C \rightarrow \C P^1$ are isomorphisms satisfying $\psi_f(\zeta^0) = 0$, $\psi_f(\zeta^{p_f}) = \infty$. More generally, define a {\em broken lollipop} to be a genus zero surface $C$ with nodes, equipped with a collection $\bfzeta$ of smooth points, and with a map satisfying the same incidence condition \eqref{eq:incidence} as before. The stability condition says that on each irreducible component of $C$ which carries less than three special points (nodes or marked points), $\bfpsi$ may not be constant; in our case, stability and the nature of \eqref{eq:incidence} rule out having just one special point on a component, so there are always at least two. For $d+|F| \geq 2$, we then have the moduli space of lollipops $\scrM^{d+1,\bfp}$ and its compactification $\bar\scrM^{d+1,\bfp}$, the space of stable broken lollipops, considered as a subspace of the space of stable maps with target $Y$. The compactification has a decomposition into strata $\scrM^{T,\bfF}$ indexed by data similar to that in \eqref{eq:compactification}: the difference is that instead of a ribbon structure, $T$ only comes with a labeling of its semi-infinite edges by $\{0,\dots,d\}$. The pair $(T,\bfF)$ fixes the topological type of the broken lollipop (the structure of the nodal curve, including the choice of irreducible components on which the marked points lie; and the degrees of the map $\bfpsi$ on each irreducible component). There is a natural action of $(\C^*)^F$ on $Y$ by rescaling coordinates, and therefore an induced action on our moduli spaces. Similarly, we have an action of $\Sym^{\bfp}$, which comes from permuting the factors in the target space $Y$. In both cases, the important point is that each $H^k$ remains invariant.

\begin{lemma} \label{th:stable-maps}
$\bar\scrM^{d+1,\bfp}$ is a smooth compact complex manifold.
\end{lemma}

\proof The moduli space of all genus zero stable maps with target $Y$ and degree $(1,\dots,1)$ is known to be smooth and compact: compactness is a general property, and smoothness follows from the convexity of $Y$ \cite[Definition 2.4.2]{kontsevich-manin94}. By definition, $\bar\scrM^{d+1,\bfp}$ is a closed subspace, hence itself compact. Smoothness does not follow directly, but the proof proceeds along well-established lines. First of all, the infinitesimal deformation space ${\mathcal T}$ of any stable broken lollipop sits in a long exact sequence
\begin{multline} \label{eq:deformation-sequence}
 0 \longrightarrow Ext^0(\Omega^1_C(\zeta^0+\cdots+\zeta^d),{\mathcal O}_C) \longrightarrow
 H^0(C,{\mathcal E}) \longrightarrow {\mathcal T} \rightarrow \\ \longrightarrow Ext^1(\Omega^1_C(\zeta^0+\cdots+\zeta^d),{\mathcal O}_C)
 \longrightarrow H^1(C,{\mathcal E}),
\end{multline}
where ${\mathcal E} \subset \bfpsi^*TY$ is the subsheaf of vector fields which at each marked point $\zeta^k$ are tangent to $H^k$. More specifically, ${\mathcal E}$ is a direct sum of line bundles ${\mathcal E}_f \subset \psi_f^*T\C P^1$. Each $\psi_f^*T\C P^1$ itself is of degree $2$ on exactly one irreducible component of $C$, and trivial on all others. Moreover, ${\mathcal E}_f \subset \psi_f^*T\C P^1$ is defined by imposing vanishing conditions at the points $\zeta^0,\zeta^{p_f}$. From this and the fact that our curve has genus zero, it is easy to see that $H^0(C,{\mathcal E}_f) \iso \C$, $H^1(C,{\mathcal E}_f) = 0$. Hence the obstruction term, which is the last one in \eqref{eq:deformation-sequence}, vanishes, ensuring smoothness. \qed

This computation has another useful consequence. Fix a stable broken lollipop, modelled on a tree $T$. Nodes of $C$ then correspond canonically to interior edges $e$ of $T$. Choose, for each node, local holomorphic coordinates $\{z_{e,\pm}\}$ near its preimages in the normalization. Then there is a deformation
\begin{equation} \label{eq:smoothing}
\scrC \longrightarrow \mathcal{U}
\end{equation}
of the underlying surface with marked points $(C,\{\zeta^0,\dots,\zeta^d\})$, parametrized by some neighbourhood of the origin $\mathcal U \subset \prod_e \C$. The local model for the deformation is the standard smoothing
\begin{equation} \label{eq:local-smoothing}
z_{e,-}z_{e,+} = \delta_e,
\end{equation}
applied to each node. Possibly after shrinking $\mathcal U$, the map $\bfpsi$ extends to $\Psi: \scrC \rightarrow Y$, in such a way that \eqref{eq:incidence} is satisfied on each fibre. This is a general fact about deformations of lollipops, which follows directly from the Kodaira-Spencer deformation theory of holomorphic maps: in the notation from Lemma \ref{th:stable-maps}, the obstruction to such an extension lies in $H^1(C,{\mathcal E})$, which we know to be zero. Of course, the extension is not unique, but the ambiguity is easy to understand: the group of maps $\mathcal U \rightarrow (\C^*)^F$ which equal $(1,\dots,1)$ at the origin acts simply transitively on the set of all possible choices. Once we have chosen $\Psi$, the resulting family is classified by a holomorphic map $\mathcal U \rightarrow \bar\scrM^{d+1,\bfp}$ (this uses the universality property of the moduli space, which in turn is based on the fact that stable broken lollipops have no nontrivial automorphisms).

Instead of starting with a single stable broken lollipop, one can also take a family, in particular that corresponding to one of the strata $\scrM^{T,\bfF}$. Suppose that for each surface in the family, we have chosen local coordinates near the nodes, as in the previous discussion. Given that, carry out the smoothing process, and find suitable maps $\Psi$. As before, this yields a classifying map, which now takes a form analogous to \eqref{eq:gluing-map}:
\begin{equation} \label{eq:complex-gluing-map}
\scrM^{T,\bfF} \times \prod_e \C \supset {\mathcal U}^{T,\bfF} \longrightarrow
\bar\scrM^{d+1,\bfp}.
\end{equation}

\begin{lemma} \label{th:tubular}
Possibly after shrinking ${\mathcal U}^{T,\bfF}$, \eqref{eq:complex-gluing-map} is a diffeomorphism onto a neighbourhood of $\scrM^{T,\bfF} \subset \bar\scrM^{d+1,\bfp}$.
\end{lemma}

\proof Take a stable broken lollipop $(C,\bfpsi)$ representing a point of $\scrM^{T,\bfF}$. Let $\tilde{C} \rightarrow C$ be the normalization, and $\{\nu_{e,\pm}\}$ the preimages of the nodes in $\tilde{C}$. Concerning the tangent space of our stratum at this point, we have an exact sequence
\begin{multline}
 0 \longrightarrow Ext^0(\Omega^1_C(\zeta^0+\cdots+\zeta^d),{\mathcal O}_C) \longrightarrow
 H^0(C,{\mathcal E}) \longrightarrow T\scrM^{T,\bfF} \rightarrow \\ \longrightarrow Ext^1(\Omega^1_{\tilde{C}}(\textstyle \sum_k \zeta^k + \sum_e \nu_{e,+} + \sum_e \nu_{e,-}),{\mathcal O}_{\tilde{C}})
 \longrightarrow 0.
\end{multline}
By comparing this with \eqref{eq:deformation-sequence}, one sees that $T\scrM^{d+1,\bfp}/T\scrM^{T,\bfF}$ can be identified with the quotient $Ext^1(\Omega^1_C(\sum_k \zeta^k),{\mathcal O}_C)/Ext^1(\Omega^1_{\tilde{C}}(\sum_k \zeta^k + \sum_e \nu_{e,+} + \sum_e \nu_{e,-}),{\mathcal O}_{\tilde{C}})$. On the other hand, there is a well-known short exact sequence \cite[p.\ 100]{harris-morrison} which identifies that quotient with
\begin{equation} \label{eq:quotient-deformations}
\bigoplus_e\, T\tilde{C}_{\nu_{e,+}} \!\otimes T\tilde{C}_{\nu_{e,-}}.
\end{equation}
Geometrically, projection from $Ext^1(\Omega^1_C(\sum_k \zeta^k),{\mathcal O}_C)$ to \eqref{eq:quotient-deformations} measures the amount of smoothing that is applied to the nodes (the reduction to first order of the parameter $\delta_e$ from \eqref{eq:local-smoothing}, essentially). This means that the differential of \eqref{eq:complex-gluing-map} induces an isomorphism from the normal bundle of $\scrM^{T,\bfF}$ to \eqref{eq:quotient-deformations}.
\qed

\begin{addendum} \label{th:differentiable}
In the construction of \eqref{eq:complex-gluing-map}, we have implicitly assumed that the local coordinates used in the gluing process are not only holomorphic on each single surface, but also depend holomorphically on the moduli. In our intended applications, the dependance is in fact only differentiable, which makes it difficult to apply the universality property directly. One way to get around this problem is to restrict attention to {\em rational} local coordinates, which are ones that extend to isomorphisms from $\C P^1$ to the relevant component of the broken lollipop. There is an enhanced moduli space $\hat\scrM^{T,\bfF}$ of broken lollipops equipped with such coordinates for all the nodes, and this comes with a canonical holomorphic gluing map
\begin{equation} \label{eq:universal-complex-gluing}
\hat\scrM^{T,\bfF} \times \prod_e \C \supset \hat{\mathcal U}^{T,\bfF} \longrightarrow
\bar\scrM^{d+1,\bfp}.
\end{equation}
On the other hand, the forgetful map $\hat\scrM^{T,\bfF} \rightarrow \scrM^{T,\bfF}$ is a holomorphic principal bundle. A choice of rational local coordinates, depending differentiably on moduli, is simply a smooth section of that bundle. One then defines the generalization of \eqref{eq:complex-gluing-map} by composing this section and \eqref{eq:universal-complex-gluing}. This is obviously differentiable, and it is easy to see that Lemma \ref{th:tubular} continues to hold.
\end{addendum}

\subsection{Real structures\label{subsec:real-gluing}}
The complex conjugate of a lollipop is the complex conjugate $\bar{C}$ of the underlying curve, equipped with the map $\bar\bfpsi$. A real lollipop is one which is invariant under this operation. This means that it has an anti-holomorphic involution $\iota: C \rightarrow C$ preserving each $\zeta^k$, and such that $\bfpsi \circ \iota = \bar\bfpsi$. The same applies to broken lollipops. Note that in the stable case, a real structure (if it exists) is necessarily unique. Complex conjugation induces an anti-holomorphic involution of the compactified moduli space, whose fixed point set $\bar\scrM^{d+1,\bfp}_\R$ parametrizes real broken lollipops.

We now revisit \eqref{eq:smoothing} in the presence of real structures. Start with $(C,\bfpsi)$ which is real. For each node of $C$, choose local holomorphic coordinates around its preimages in the normalization, and assume that these coordinates are compatible with the anti-holomorphic involution. In that case, the family $\scrC$ inherits a real structure as well, and the map $\Psi$ on it can be chosen to be compatible with that structure. This is still not unique, but we can make a more specific choice based on additional geometric data, as follows. Suppose that $(C,\bfpsi)$ is modelled on $(T,\bfF)$. This means in particular that the irreducible components $C_v$ are indexed by vertices of $T$, and that for each $f \in F_v$, the restriction $\psi_f|C_v: C_v \rightarrow \C P^1 = \C \cup \{\infty\}$ is an isomorphism. Choose, for each $f$, an additional marked point $q_f$ lying in $\psi_f^{-1}(S^1) \subset C_v$, which should not be a node. Since the deformation $\scrC$ is defined by a local process near the nodes, these points will survive into the nearby fibres, yielding sections $\kappa_f: \mathcal U \rightarrow \scrC$ of \eqref{eq:smoothing}. For any choice of $\Psi$ compatible with the real structure, we get an associated collection of real analytic functions
\begin{equation} \label{eq:real-psi}
|\Psi_f \circ \kappa_f|: {\mathcal U}_\R \longrightarrow \R^*,
\end{equation}
each of which sends the origin to $1$. If we change $\Psi_f$ by a a real analytic map ${\mathcal U} \rightarrow \C^*$, the effect is to multiply \eqref{eq:real-psi} with the restriction of the same map to ${\mathcal U}_\R$. Hence, possibly after shrinking ${\mathcal U}_\R$, there is a unique choice of $\Psi$ such that all the \eqref{eq:real-psi} become constant equal to $1$. We call the resulting $(\scrC, \Psi)$ a {\em geometric deformation} of the real broken lollipop $(C,\bfpsi)$. Such deformations depend on the choice of local holomorphic coordinates near the nodes, and on the points $q_f$.

Take a stratum $\scrM^{T,\bfF}_\R$ in the real part of the lollipop moduli space. Suppose that along this stratum, we choose rational local coordinates around the nodes, as well as additional marked points as described above. For the same reason as in Addendum \ref{th:differentiable}, it is sufficient to have these coordinates and points vary differentiably over the moduli space. The geometric deformation construction then yields a gluing map, whose real part is a diffeomorphism onto a neighbourhood of our stratum:
\begin{equation} \label{eq:real-gluing}
\scrM^{T,\bfF}_\R \times \prod_e \R \supset {\mathcal U}^{T,\bfF}_\R \longrightarrow
\bar\scrM^{d+1,\bfp}_\R.
\end{equation}

\subsection{Popsicles vs lollipops\label{subsec:real-vs-complex}}
Take, as before, a tree $T$ with $d+1$ semi-infinite edges, labeled by $\{0,\dots,d\}$. We say that $T$ is of ribbon type if it can be embedded into the plane in such a way that the given numbering of the semi-infinite edges is compatible with the counterclockwise cyclic ordering. If such an embedding exists, it is essentially unique, so we do indeed get a distinguished ribbon structure in the ordinary sense. Suppose from now on that $T$ is of this type; and let $(C,\bfpsi)$ be a stable broken lollipop modelled on $(T,\bfF)$, which carries a real structure. From the ribbon structure of $T$, we get a canonical ordering of the flags adjacent to any given vertex $v$ (compare Section \ref{subsec:define-popsicle}), hence of the special points (nodes and marked points) on each irreducible component $C_v$. We say that the real structure is admissible if there is some orientation of the real parts $C_{v,\R} \iso S^1$ which is compatible with the ordering of the special points, and such that for each $f \in F_v$, the diffeomorphism $\psi_{f,\R}: C_{v,\R} \rightarrow \R \cup \{\infty\} = Y_{f,\R}$ becomes orientation-preserving. Again, there can be at most one such orientation. Supposing that one exists, one can construct a popsicle modelled on $(T,\bfF)$ as follows. For each $v$, take $\bar{S}_v$ to be the disc in $C_v$ whose boundary is $C_{v,\R}$ with the positive orientation, and then remove all special points to get a pointed disc $S_v$. The popsicle maps are defined by
\begin{equation} \label{eq:logarithm}
\phi_f = \textstyle\frac{1}{\pi} \log(\psi_f) \,|\, S_v.
\end{equation}
Conversely, the broken popsicle with components $S_v$ and maps $\phi_f$, $f \in F_v$, determines the original $(C,\psi)$. Now suppose that we are given popsicle sticks on each $S_v$, which allow one to convert the information contained in the popsicle maps $\phi_f$ into sprinkles $q_f \in S_v \subset C_v$. Because of \eqref{eq:logarithm}, these points necessarily satisfy $\psi_f(q_f) \in S^1$. Hence, the geometric deformation introduced above actually reproduces the gluing process from Section \ref{subsec:sticks}. More precisely, geometric deformation with non-positive gluing parameter corresponds to gluing of popsicles, whereas a positive gluing parameter results in a non-admissible real structure.

To conclude, $\bar\scrR^{d+1,\bfp}$ can be identified with the subset of $\bar\scrM^{d+1,\bfp}_\R$ consisting of those broken lollipops which have admissible real structures. Consider some boundary stratum $\scrR^{T,\bfF}$, and suppose that for each family $\scrS^{|v|,\bfp_v} \rightarrow \scrR^{|v|,\bfp_v}$, we have made choices of rational strip-like ends as well as popsicle sticks. The associated gluing map \eqref{eq:gluing-with-sprinkles} is a restriction of \eqref{eq:real-gluing}, which means that

\begin{cor} \label{th:glue-well}
Possibly after making $\mathcal U^{T,\bfF}$ smaller, \eqref{eq:gluing-with-sprinkles} is a diffeomorphism onto a neighbourhood of $\scrR^{T,\bfF}$ inside $\bar\scrR^{d+1,\bfp} \subset \bar\scrM^{d+1,\bfp}_\R$.
\end{cor}

\begin{cor} \label{th:real-and-negative}
$\bar\scrR^{d+1,\bfp} \subset \bar\scrM^{d+1,\bfp}_\R$ is a smooth submanifold with corners.
\end{cor}

\subsection{Compactifying the universal family\label{subsec:singularity}}
We now discuss, in slightly less detail, the construction of the compactified total space $\bar\scrS^{d+1,\bfp}$. Consider stable maps $\bfpsi: C \rightarrow Y$ of the same genus and degree as before, but where $C$ now carries $d+2$ marked non-special points $\{\zeta^0,\dots,\zeta^d,z\}$. The first $d+1$ of these should satisfy the usual incidence condition \eqref{eq:incidence}, and the last one can be mapped freely to $Y$. On one hand, the computation from Lemma \ref{th:stable-maps} shows that such maps form a smooth moduli space $\scrC^{d+1,\bfp}$. On the other hand, it is part of the general theory of stable maps (the prototypical case of stable curves was treated in \cite{knudsen83}) that there is a natural holomorphic forgetful map
\begin{equation} \label{eq:forget-z}
\bar\scrC^{d+1,\bfp} \longrightarrow \bar\scrM^{d+1,\bfp},
\end{equation}
and that this can in fact be identified with the universal family of stable curves over the moduli space.

Consider the subset $\bar\scrS^{d+1,\bfp} \subset \bar\scrC^{d+1,\bfp}$ consisting of those points $(C,\bfzeta \cup \{z\},\bfpsi)$ which satisfy the following two conditions. First, the image under \eqref{eq:forget-z} should lie in $\bar\scrR^{d+1,\bfp}$. Secondly, if it is the case that $(C,\bfzeta,\bfpsi)$ is stable (note that we are forgetting $z$ here), we require additionally that $z$ should lie in the subset of $C$ given by the associated compactified broken popsicle. One checks easily that (as a set) this can be identified with the compactification of the universal family \eqref{eq:compactified-total-space}. The case where forgetting $z$ makes the map unstable deserves special mention. This occurs whenever $\bfpsi$ is constant on an irreducible component of $C$ which contains $z$ as well as
\begin{equation}
\parbox{30em}{
\begin{itemize}
\item[(i)] exactly one of the points $\zeta^k$, and exactly one nodal point; or \item[(ii)] none of the $\zeta^k$, and exactly two nodal points.
\end{itemize}
}
\end{equation}
In terms of \eqref{eq:compactified-total-space}, this corresponds to having the compactification $\bar{S}$ of an (ordinary or broken) popsicle, together with a point $z \in \bar{S}$ which either (i) agrees with one of the points at infinity, or (ii) is a singular point. The correspondence is shown schematically in Figure \ref{fig:unstable}.
\begin{figure}
\begin{centering}
\begin{picture}(0,0)%
\includegraphics{unstable.pstex}%
\end{picture}%
\setlength{\unitlength}{3947sp}%
\begingroup\makeatletter\ifx\SetFigFont\undefined%
\gdef\SetFigFont#1#2#3#4#5{%
  \reset@font\fontsize{#1}{#2pt}%
  \fontfamily{#3}\fontseries{#4}\fontshape{#5}%
  \selectfont}%
\fi\endgroup%
\begin{picture}(4958,5698)(-299,-2726)
\put(-299,764){\makebox(0,0)[lb]{\smash{{\SetFigFont{10}{12.0}{\rmdefault}{\mddefault}{\updefault}{\color[rgb]{0,0,0}(ii)}%
}}}}
\put(526,524){\makebox(0,0)[lb]{\smash{{\SetFigFont{10}{12.0}{\rmdefault}{\mddefault}{\updefault}{\color[rgb]{0,0,0}$\zeta^1$}%
}}}}
\put(601,-111){\makebox(0,0)[lb]{\smash{{\SetFigFont{10}{12.0}{\rmdefault}{\mddefault}{\updefault}{\color[rgb]{0,0,0}$\zeta^0$}%
}}}}
\put(3001,524){\makebox(0,0)[lb]{\smash{{\SetFigFont{10}{12.0}{\rmdefault}{\mddefault}{\updefault}{\color[rgb]{0,0,0}$\zeta^2$}%
}}}}
\put(3001,-111){\makebox(0,0)[lb]{\smash{{\SetFigFont{10}{12.0}{\rmdefault}{\mddefault}{\updefault}{\color[rgb]{0,0,0}$\zeta^3$}%
}}}}
\put(1701,934){\makebox(0,0)[lb]{\smash{{\SetFigFont{10}{12.0}{\rmdefault}{\mddefault}{\updefault}{\color[rgb]{0,0,0}$z$}%
}}}}
\put(2701,-1411){\makebox(0,0)[lb]{\smash{{\SetFigFont{10}{12.0}{\rmdefault}{\mddefault}{\updefault}{\color[rgb]{0,0,0}$\zeta^3$}%
}}}}
\put(2776,-2686){\makebox(0,0)[lb]{\smash{{\SetFigFont{10}{12.0}{\rmdefault}{\mddefault}{\updefault}{\color[rgb]{0,0,0}$\zeta^2$}%
}}}}
\put(676,-2686){\makebox(0,0)[lb]{\smash{{\SetFigFont{10}{12.0}{\rmdefault}{\mddefault}{\updefault}{\color[rgb]{0,0,0}$\zeta^1$}%
}}}}
\put(751,-1411){\makebox(0,0)[lb]{\smash{{\SetFigFont{10}{12.0}{\rmdefault}{\mddefault}{\updefault}{\color[rgb]{0,0,0}$\zeta^0$}%
}}}}
\put(1876,-2086){\makebox(0,0)[lb]{\smash{{\SetFigFont{10}{12.0}{\rmdefault}{\mddefault}{\updefault}{\color[rgb]{0,0,0}$z$}%
}}}}
\put(546,2869){\makebox(0,0)[lb]{\smash{{\SetFigFont{10}{12.0}{\rmdefault}{\mddefault}{\updefault}{\color[rgb]{0,0,0}$z$}%
}}}}
\put(-224,2189){\makebox(0,0)[lb]{\smash{{\SetFigFont{10}{12.0}{\rmdefault}{\mddefault}{\updefault}{\color[rgb]{0,0,0}$\zeta^0$}%
}}}}
\put(1501,1814){\makebox(0,0)[lb]{\smash{{\SetFigFont{10}{12.0}{\rmdefault}{\mddefault}{\updefault}{\color[rgb]{0,0,0}$\zeta^1$}%
}}}}
\put(2026,1889){\makebox(0,0)[lb]{\smash{{\SetFigFont{10}{12.0}{\rmdefault}{\mddefault}{\updefault}{\color[rgb]{0,0,0}$\zeta^2$}%
}}}}
\put(1726,2464){\makebox(0,0)[lb]{\smash{{\SetFigFont{10}{12.0}{\rmdefault}{\mddefault}{\updefault}{\color[rgb]{0,0,0}$\zeta^3$}%
}}}}
\put(4126,2864){\makebox(0,0)[lb]{\smash{{\SetFigFont{10}{12.0}{\rmdefault}{\mddefault}{\updefault}{\color[rgb]{0,0,0}$\zeta^3$}%
}}}}
\put(3526,2189){\makebox(0,0)[lb]{\smash{{\SetFigFont{10}{12.0}{\rmdefault}{\mddefault}{\updefault}{\color[rgb]{0,0,0}$z = \zeta^0$}%
}}}}
\put(3526,1514){\makebox(0,0)[lb]{\smash{{\SetFigFont{10}{12.0}{\rmdefault}{\mddefault}{\updefault}{\color[rgb]{0,0,0}$\zeta^1$}%
}}}}
\put(4426,1514){\makebox(0,0)[lb]{\smash{{\SetFigFont{10}{12.0}{\rmdefault}{\mddefault}{\updefault}{\color[rgb]{0,0,0}$\zeta^2$}%
}}}}
\put(-299,2714){\makebox(0,0)[lb]{\smash{{\SetFigFont{10}{12.0}{\rmdefault}{\mddefault}{\updefault}{\color[rgb]{0,0,0}(i)}%
}}}}
\end{picture}%
\caption{\label{fig:unstable}}
\end{centering}
\end{figure}%

It is now easy to construct local models for $\bar\scrS^{d+1,\bfp}$. Near any point where $(C,\bfzeta,\psi)$ is stable, this is a manifold with corners, the codimension of the corner being given as usual by the number of double points. In the remaining case (i), if we consider the node lying on the same component as $z$, the relevant gluing parameter does not have to be real, but (with suitable conventions) it should lie in the upper half-plane; otherwise, the outcome of the gluing process does not lie in the subspace $\bar\scrS^{d+1,\bfp}$. Hence, the local model is a corner whose codimension is the number of double points of the image of our point under \eqref{eq:forget-z}. Finally, in case (ii) one has two gluing parameters $\delta_\pm$ associated to the nodes which lie on the same component as $z$. Both of them must lie in the upper half-plane, and moreover, their product needs to be real and positive, so we get a singularity of type \eqref{eq:singularity}; this is particularly easy to see if one thinks of $\arg(\delta_\pm)$ as an angle of rotation in Figure \ref{fig:unstable}.

As we have briefly mentioned before, even though $\bar\scrS^{d+1,\bfp}$ is singular, one can define the notion of smooth function on that space, either by being smooth with respect to the local coordinates mentioned above, or (equivalently) as the restriction of a smooth function on $\bar\scrC^{d+1,\bfp}$. For our purposes, only a very special class of such functions will be relevant. Namely, suppose that $h: \scrS^{d+1,\bfp} \rightarrow \R$ is a smooth function with the following property: it admits a smooth extension $\bar{h}$ to the compactification, which vanishes to arbitrary order at special points (points at infinity and double points). From the more concrete viewpoint of gluing parameters, this property has the following implications. First, consider a single popsicle $(S,\bfphi)$, identified with some fibre of the universal family $\scrS^{d+1,\bfp}$, and let $h_S$ be the restriction of $h$ to that fibre. Take a strip-like end, written for simplicity as an inclusion $Z_\pm \hookrightarrow S$. Then, for every $C>0$ there is a constant $D>0$ such that on the strip-like end,
\begin{equation} \label{eq:superexponential-1}
|h(s,t)|  < D e^{\pm Cs}.
\end{equation}
In words, $h$ decays faster than exponentially near the points at infinity. To get \eqref{eq:superexponential-1} from the given vanishing condition, one simply notices that $\exp(-\pi(\pm s \pm it))$ is part of the local coordinate system on $\bar\scrS^{d+1,\bfp}$ centered at a point at infinity; therefore, $h$ grows smaller than any of its powers.

Next, suppose that we have a one-parameter family of popsicles $(S_l,\bfphi_l)$, obtained by gluing together two fixed popsicles $(S_\pm,\bfphi_\pm)$ with gluing length $l \gg 0$. By definition, $S_l$ contains a finite strip $[0,l] \times [0,1]$. On that strip, the functions $h_{S_l} = h|S_l$ must have the following behaviour. For any $C>0$ there is a constant $D>0$ (independent of $l$) such that
\begin{equation} \label{eq:superexponential-2}
|h_{S_l}(s,t)| \leq \min\{ D e^{-C s}, D e^{C(s-l)} \}.
\end{equation}
This means that $h$ decays faster than exponentially (uniformly in $l$), whether one counts from the left end or from the right end of the strip. The reason is essentially the same as in \eqref{eq:superexponential-1}: in terms of \eqref{eq:singularity}, the point $(s,t)$ on the neck of the surface $S_l$ has local coordinates $\delta_- = -\exp(-\pi(s+it))$, $\delta_+ = -\exp(-\pi(l-s-it))$. We know from the assumption that locally near $\delta_- = \delta_+ = 0$, $|h|$ is less than any $|\delta_-|^C$, $|\delta_+|^C$. This yields \eqref{eq:superexponential-2} on a piece $[L,l-L] \times [0,1]$, where the constant $L$ is independent of $l$. On the other hand, if one restricts $h_{S_l}(s,t)$ to $[0,L] \times [0,1]$, it converges uniformly as $l \rightarrow \infty$, and the same is true for $h_{S_l}(s-l,t)$ on $[-L,0] \times [0,1]$. Using that fact, one easily extends the given bound to the whole of $[0,l] \times [0,1]$.

Bounds similar to \eqref{eq:superexponential-1}, \eqref{eq:superexponential-2} hold for the derivatives of $h_{S_l}(s,t)$. There are also analogues where one considers families of popsicles, as well as more complicated gluing processes.

\begin{remark}
We have used $\bar\scrS^{d+1,\bfp}$ because it is the most straightforward way of compactifying the universal family. There is at least one other possibility, which would be to consider triples $(C,\bfzeta \cup \{z,\bar{z}\},\phi)$, admitting a real involution which exchanges the two additional marked points. The compactification obtained in this way is somewhat larger, but has the advantage of being a genuine manifold with corners. If we restrict to a single popsicle $(S,\bfphi)$ representing a fibre of $\scrS^{d+1,\bfp}$, what this alternative compactification does is to add an {\em interval at infinity} to each end, rather than a single point.
\end{remark}

\section{A priori estimates\label{sec:compactness}}

This section collects various elementary analytic arguments.  Action and energy considerations encode what, in more standard terms, is known as the technique of monotone homotopies, going back at least to \cite{floer-hofer94}. We also need an elementary SFT-type convexity result, which constrains pseudo-holomorphic curves to certain parts of the target manifold; this replaces the maximum principle and Monotonicity Lemma arguments used for the same purpose in \cite{viterbo97a}.

\subsection{Preliminaries\label{subsec:prelims}}
Let's start by summarizing the geometric data that enters into the formulation of our inhomogeneous $\bar\partial$-equations (the following considerations will be largely independent of how that data is chosen in the specific context of wrapped Fukaya categories). On the target space side, let $M$ be a Weinstein domain, $L \subset M$ a Lagrangian submanifold satisfying \eqref{eq:asymptotically-parallel} as well as \eqref{eq:no-integer-chords}, and $H$ a function as in \eqref{eq:class-of-h}, \eqref{eq:nondegenerate-x}. On the source side, we have a pointed disc $S$, which comes with strip-like ends $\bfepsilon$ as well as weights $\bfw$. Besides that, $S$ carries a one-form $\gamma$ satisfying \eqref{eq:gamma-properties}. By applying Stokes one sees that
\begin{equation} \label{eq:stokes}
\textstyle{-\int_S}\,d\gamma = w^0 - w^1 - \cdots - w^d.
\end{equation}
In particular, $d\gamma = 0$ if and only if $w^0$ is the sum of the other weights. Finally, we have a family $\bfJ = \{J_z\}$ of almost complex structures parametrized by points of $S$, each of which is of contact type at the boundary. Over the strip-like ends $\epsilon^k(s,t)$, this family should converge to one which depends only on $k$ and $t$, and the convergence is assumed to be faster than exponential in any $C^r$ topology. Given that, we consider solutions of \eqref{eq:dbar}, which (to save the reader the trouble of looking it up) means maps
\begin{equation} \label{eq:restate-dbar}
\left\{
\begin{aligned}
 & u: S \longrightarrow M, \\
 & u(\partial S) \subset L, \\
 & \textstyle\lim_{s \rightarrow \pm\infty} u(\epsilon^k(s,\cdot)) = x^k \in \scrX_{w^k}, \\
 & (du-X \otimes \gamma)^{0,1} = 0.
\end{aligned}
\right.
\end{equation}
If we choose local holomorphic coordinates $z = s+it$ around some point of $S$, and write accordingly $\gamma = p(s,t) ds + q(s,t) dt$, the $\bar\partial$-equation in \eqref{eq:restate-dbar} becomes
\begin{equation} \label{eq:local-floer}
\partial_t u - q(s,t)X = J_{s,t}(\partial_s u - p(s,t)X).
\end{equation}
This class of equations is compatible with certain changes of variables, given by the action of the flow $\phi$ of $X$. More precisely, if we make the ansatz
\begin{equation} \label{eq:ansatz}
u(s,t) = \phi^{l(s,t)}(\tilde{u}(s,t)),
\end{equation}
then $\tilde{u}$ satisfies the same equation with modified data
\begin{equation} \label{eq:gauge-data}
\begin{aligned}
 & \tilde{J}_{s,t} = \big(\phi^{l(s,t)}\big)^*J_{s,t}, \\
 & \tilde{p}(s,t) = p(s,t) - \partial_s l(s,t), \\
 & \tilde{q}(s,t) = q(s,t) - \partial_t l(s,t).
\end{aligned}
\end{equation}
In a less coordinate-bound way, what this means is that one has replaced $\gamma$ with $\tilde\gamma = \gamma-dl$. Suppose that the point we are considering is an interior point of $S$, and that $d\gamma = 0$ in a neighbourhood of that point. If we then take $l$ such that $dl = \gamma$, the associated transformation turns \eqref{eq:local-floer} into a standard pseudo-holomorphic map equation (still with almost complex structure which depends on $z$). Similarly, near boundary points (where $d\gamma$ always vanishes, by assumption) one has $\gamma = dl$, where $l$ is constant along the boundary, and that transforms our equation into a standard pseudo-holomorphic map equation with Lagrangian boundary conditions.

In the general case where $d\gamma \neq 0$, one can still get rid of the inhomogeneous term by applying the traditional Gromov trick. Namely, consider $S \times M$ with the symplectic form $\omega - d(H\gamma) + \text{\it (a large positive two-form on $S$)}$. There is a unique compatible almost complex structure such that the equation for pseudo-holomorphic sections $S \rightarrow S \times M$ becomes equivalent to $(du - X \otimes \gamma)^{0,1} = 0$. One consequence of this reduction, and of standard pseudo-holomorphic curve theory, is the following:

\begin{lemma} \label{th:intersection}
Two different solutions of \eqref{eq:restate-dbar} can only agree on a discrete set of points.
\end{lemma}

\subsection{Energy and action\label{subsec:energy}}
The geometric and topological energies of a solution of \eqref{eq:restate-dbar} are defined by
\begin{equation} \label{eq:energy}
\begin{aligned}
 E^{geom}(u) & = \int_S \half \|du - X \otimes \gamma\|^2 = \displaystyle\int_S u^*\omega - u^*dH \wedge \gamma, \\
 E^{top}(u) & = \int_S u^*\omega - d(u^*H \cdot \gamma) = E^{geom}(u) - \int_S u^*H \cdot d\gamma.
\end{aligned}
\end{equation}
Since $H > 0$ and $d\gamma \leq 0$, it follows that
\begin{equation} \label{eq:energy-inequality}
0 \leq E^{geom}(u) \leq E^{top}(u).
\end{equation}
Moreover, the first inequality is an equality if and only if $du = X \otimes \gamma$ (see Section \ref{subsec:trivial-solutions} for further discussion of these special solutions), while the second one is an equality if and only if $d\gamma = 0$. On the other hand, the topological energy can be computed in terms of action differences by a direct application of Stokes's theorem. Namely, define the action of $x \in \scrX_w$ to be
\begin{equation} \label{eq:action}
A_{wH}(x) = \int_0^1 -x^*\theta + wH(x(t))\,dt + h(x(1)) - h(x(0)),
\end{equation}
where $h$ is as in \eqref{eq:asymptotically-parallel}. Then, for any solution of \eqref{eq:restate-dbar} with limits $x^k$,
\begin{equation} \label{eq:a-priori}
 E^{top}(u) = A_{w^0H}(x^0) - \sum_{k=1}^d A_{w^kH}(x^k).
\end{equation}

\subsection{A convexity argument\label{subsec:sft}}
Suppose that we have a manifold with boundary $M$, together with a one-form $\theta$ such that $d\theta = \omega$ is symplectic, but where now the Liouville field $Z$ points strictly inwards along $\partial M$. This is sometimes called a {\em concave} contact type boundary, and the natural completion of such a manifold is to attach the small part of the symplectization $(0,1] \times \partial M$ to the boundary. Let $L \subset M$ be a Lagrangian submanifold satisfying \eqref{eq:asymptotically-parallel}, where we assume additionally that
\begin{equation} \label{eq:f-vanishes-along-the-boundary}
h|\partial L = 0.
\end{equation}
Suppose also that we have a function $H$ satisfying the analogue of \eqref{eq:class-of-h}. Finally, there is a natural notion of almost complex structure which is of contact type on $\partial M$, analogous to the one in Section \ref{subsec:j}.

Take a connected compact Riemann surface $S$ with corners. The corners divide $\partial S$ into pieces which are closed intervals or circles. In our case, each piece should carry an $n$ or $l$ label, with the property that two pieces ending at the same corner must be labeled differently. We denote by $\partial S = \partial_n S\cup \partial_l S$ the resulting decomposition, and assume that $\partial_n S \neq \emptyset$. Suppose also that our surface comes with some $\gamma \in \Omega^1(S)$, satisfying $\gamma|\partial_l S = 0$ and $d\gamma \leq 0$. Additionally, we want to have a family $\bfJ$ of almost complex structures, parametrized by points of $S$, which are of contact type on $\partial M$. Given this, consider the equation
\begin{equation} \label{eq:j-concave}
\left\{
\begin{aligned}
 & u: S \longrightarrow M, \\
 & u(\partial_l S) \subset L, \\
 & u(\partial_n S) \subset \partial M, \\
 & (du-X \otimes \gamma)^{0,1} = 0.
\end{aligned}
\right.
\end{equation}

\begin{lemma} \label{th:penetrate}
There are no solutions of \eqref{eq:j-concave} unless $d\gamma = 0$. Even then, all solutions must satisfy $du = X \otimes \gamma$ and $u(S) \subset \partial M$.
\end{lemma}

\proof We apply Stokes to the integral defining the topological energy \eqref{eq:energy}:
\begin{align} \label{eq:integrate-out2}
E^{top}(u)  & = \int_{\partial S} u^*\theta - u^*H \cdot \gamma.
\end{align}
By assumption, we have $H = \theta(X) = 1$ along $\partial M$, and moreover $\theta|L = dh$, where $h$ vanishes along $\partial L$. Using this, one can rewrite \eqref{eq:integrate-out2} as
\begin{equation}
\begin{aligned}
E^{top}(u) & = \int_{\partial_n S} \theta \circ (du-X \otimes \gamma) \\
& = \int_{\partial_n S} (\theta \circ J_z) \circ (du-X \otimes \gamma) \circ (-j).
\end{aligned}
\end{equation}
We know that along $\partial M$, $\theta \circ J_z$ is the radial one-form (written as $dr$ in terms of the coordinates on the cone). In particular, since $X|\partial M$ is the Reeb flow, $(\theta \circ J_z)(X) = 0$. If $\xi$ is a vector tangent to $\partial_nS$, and positive with respect to the boundary orientation, then $j\xi$ points inwards, hence $du(j\xi)$ does not point outwards along $\partial M$, whence $(\theta \circ J_z \circ du)( j\xi) \geq 0$. Integrating that, we find that $E^{top}(u) \leq 0$. This contradicts \eqref{eq:energy-inequality} and the nonnegativity of the geometric energy, unless $E^{geom} = 0$ and $E^{top} = E^{geom}$, which means that $du = X \otimes \gamma$ and $d\gamma = 0$. Since $du$ is always a multiple of $X$, it follows that the image of $u$ lies inside a single orbit of that vector field, which in turn must be contained in $\partial M$ since $\partial_nS \neq \emptyset$.
\qed

Note that Lemma \ref{th:penetrate} still holds if we replace $H$ by $cH$ for some constant $c>0$, since that can be compensated by rescaling $\gamma$. This will be occasionally useful.

\subsection{Solutions escaping to infinity\label{subsec:escape}}
We now return to the original setup of equation \eqref{eq:restate-dbar}, and modify it by completing all the target space structures. This means that we consider
\begin{equation} \label{eq:completed-target}
\left\{
\begin{aligned}
& u: S \longrightarrow \hat{M}, \\
& u(\partial S) \subset \hat{L}, \\
& \textstyle\lim_{s \rightarrow \pm\infty} u(\epsilon^k(s,\cdot)) = x^k \in \scrX_{w^k}, \\
& (du-\hat{X} \otimes \gamma)^{0,1} = 0.
\end{aligned}
\right.
\end{equation}
Here, $\hat{X}$ is the Hamiltonian vector field of the natural extension $\hat{H}$ from \eqref{eq:class-of-h}. Similarly, we use the extended almost complex structures $\hat{J}_z$ which are of contact type on the infinite cone. By assumption \eqref{eq:no-integer-chords}, all the possible limits $x^k$ still lie in $M \setminus \partial M \subset \hat{M}$.

The infinite cone $[1,\infty) \times \partial M$ is a manifold with a concave contact type boundary. Suppose first that $u$ intersects $\{1\} \times \partial M$ transversally. Then, one can apply Lemma \ref{th:penetrate} to the part of $u$ that gets mapped to the infinite cone. Note that the additional condition \eqref{eq:f-vanishes-along-the-boundary} holds automatically here, since $\hat\theta$ vanishes along $[1,\infty) \times \partial L$. The outcome is that $u$ does not actually penetrate into the interior of the cone, hence stays in $M \subset \hat{M}$ after all. To remove the technical transversality assumption, take some $\epsilon \geq 0$ such that $u$ intersects $\{1+\epsilon\} \times \partial M$ transversally. Lemma \ref{th:penetrate} (or rather, the slight generalization mentioned at the end of Section \ref{subsec:sft}) applies to $[1+\epsilon,\infty) \times \partial M$, showing that the image of $u$ remains inside $M \cup ([1,1+\epsilon] \times \partial M)$. Since $\epsilon$ can be chosen arbitrarily small, we then come to the same conclusion as before.

To see the specific implications of this, suppose that we have a one-dimensional regular moduli space $\scrR^{d+1,\bfp,\bfw}(\bfx)$ and a point in that space, represented by a triple $(S_0,\bfphi_0,u_0)$, such that $u_0^{-1}(\partial M) \neq \emptyset$. Because of regularity and the dimension assumption, this point sits in a family $(S_r,\bfphi_r,u_r)$, where a priori the $u_r$ could leave $M$ and penetrate into $(1,\infty) \times \partial M \subset \hat{M}$. However, the argument above shows that this does not happen (the family of surfaces $u_r(S_r)$ becomes tangent to $\partial M$ at $r = 0$, but does not actually cross the boundary). Hence, $(S_0,\bfphi_0,u_0)$ is not a boundary point of the original moduli space.

\subsection{Solutions escaping from a subdomain}
On the target space side, let's suppose that we have $M$ together with a Liouville subdomain $M^{in}$ and a Lagrangian submanifold $L$, as in Section \ref{subsec:geometry-of-viterbo}, and assuming that \eqref{eq:no-integer-chords} holds for both $L \subset M$ and $L^{in} \subset M^{in}$. Take a Hamiltonian $H$ of the form \eqref{eq:class-of-h-2}. In particular, for each $w$ we have the splitting of the set of integer $X$-chords of length $w$ into ones lying in the interior and exterior of the dividing hypersurface, as in \eqref{eq:divide-x}.

\begin{lemma} \label{th:positive-action}
There is a constant $\nu$ such that for all $w \geq \nu$ and all $x \in \scrX_w^{out}$, the action $A_{wH}(x)$ is positive.
\end{lemma}

\proof By assumption \eqref{eq:relative-exactness}, we can extend $h|L^{out}$ to a function defined on the whole of $M^{out}$, which vanishes in a neighbourhood of $\partial M^{out} = \partial M^{in} \cup \partial M$. Note that on both hypersurfaces we have $H = 1$, with nontrivial derivative in transverse direction. Hence, we can choose our neighbourhood to be bounded by level sets of $H$ on either side, which means that it will be invariant under the flow of $X$.

Take $x \in \scrX_w^{out}$, and write the action as
\begin{equation} \label{eq:action2}
A_{wH}(x) = \int_0^1 w(H(x(t)) - dH(Z)_{x(t)}) \, dt + h(x(1)) - h(x(0)).
\end{equation}
If $x$ lies in either of the two neighbourhoods defined above, the $h$ term in \eqref{eq:action2} vanishes, and the integral is over a positive function by \eqref{eq:class-of-h-2}, which yields the desired result for any $w$. On the other hand, outside those neighbourhoods $H - dH(Z)$ is bounded below by some positive constant, and of course $h$ is itself bounded, which means that \eqref{eq:action2} will necessarily be positive if $w$ is sufficiently large.
\qed

Consider the rescaled family $H^\rho$ as in \eqref{eq:h-rho}, whose Hamiltonian vector field $X^\rho$ satisfies \eqref{eq:scaled-x}. The crucial point about this rescaling process is that it affects the action values on the inside and outside differently. Namely, for $x^\rho$ as in \eqref{eq:shrink-x},
\begin{equation} \label{eq:action-shrinks}
A_{wH^\rho}(x^\rho) = \begin{cases}
A_{wH}(x) & \textrm{if } x \in \scrX_w^{out}, \\
\rho \cdot A_{wH}(x) & \textrm{if } x \in \scrX_w^{in}.
\end{cases}
\end{equation}
Take $(S,\bfepsilon,\bfw,\gamma)$ as in Section \ref{subsec:prelims}. For the almost complex structures, we want to have the following: a family $\bfJ = (J_z)$ in $\scrJ(M)$ parametrized by $S$, which asymptotically over the ends depends only on $t$; a similar family $\bfJ^{in} = (J_z^{in})$ in $\scrJ(M^{in})$ of the same kind; and a third family interpolating between the two, in the sense of Section \ref{subsec:rho}, again with the same kind of asymptotic behaviour. For any $\rho \in (0,1]$, we then consider solutions to the parametrized equation
\begin{equation} \label{eq:restate-parametrized-dbar}
\left\{
\begin{aligned}
 & u: S \longrightarrow M, \\
 & u(\partial S) \subset L^{\rho}, \\
 & \textstyle\lim_{s \rightarrow \pm\infty} u(\epsilon^k(s,\cdot)) = x^{k,\rho}, \\
 & (du_z - X_{u(z)}^\rho \otimes \gamma_z) \circ j + J_{z,u(z)}^\rho \circ (du_z - X_{u(z)}^\rho \otimes \gamma_z) = 0.
\end{aligned}
\right.
\end{equation}
Clearly, for $\rho = 1$ this reduces to \eqref{eq:restate-dbar}. Moreover, if $u^{in}$ is a solution of the corresponding equation on $M^{in}$ (using $H^{in} = H|M^{in}$, and the almost complex structures $\bfJ^{in}$), then for each sufficiently small $\rho$, the rescaled map $\psi^{\rho} \circ u^{in}$ is a solution of \eqref{eq:restate-parametrized-dbar}. The last observation has a partial converse:

\begin{lemma} \label{th:contain}
Suppose that all weights are $\geq \nu$, $\rho$ is sufficiently small, and $x^0 \in \scrX_{w^0}^{in}$. Then every solution of \eqref{eq:restate-parametrized-dbar} is of the form $\psi^{\rho} \circ u^{in}$. In particular, its other asymptotics are $x^k \in \scrX_{w^k}^{in}$, and the image $u(S)$ is entirely contained inside $M^{in,\rho}$.
\end{lemma}

\proof Combining \eqref{eq:energy-inequality} with \eqref{eq:a-priori} and \eqref{eq:action-shrinks}, one finds that
\begin{equation} \label{eq:rho-or-not}
0 \leq E^{top}(u) = \rho \cdot A_{wH}(x^0) -
\rho \sum_{\substack{k>0 \\ x^k \in \scrX_{w^k}^{in}}} A_{wH}(x^k)
- \sum_{\substack{k>0 \\ x^k \in \scrX_{w^k}^{out}}} A_{Hw}(x^k).
\end{equation}
In the last term on the right hand side, all the actions are positive by Lemma \ref{th:positive-action}, and in fact bounded below by a positive constant which is independent of $\rho$. The other terms are bounded by $\rho$ times a constant, which depends only on the weights (the maximum of $|A_{wH}(x)|$ among the finite sets $\scrX_{w^0},\dots,\scrX_{w^d}$). This leads to a contradiction for $\rho \ll 1$ unless we assume that the last sum in \eqref{eq:rho-or-not} is empty, which means that $x^k \in \scrX_{w^k}^{in}$ for all $k$.

Having excluded the possibility that some input lies in $M^{out}$, we now prove that all solutions to \eqref{eq:restate-parametrized-dbar} lie in $ M^{in, \rho}$.  Consider $M^{out} \cup \kappa^{in}([\rho,1) \times \partial M^{in})$. This has $\partial M^{in,\rho}$ as a concave contact type boundary (of course, there is another boundary component, namely $\partial M$, but that won't be relevant). Moreover, for all $\rho \ll 1$, the almost complex structures $J^\rho_z$ are of contact type in a neighbourhood of $\partial M^{in,\rho} \subset \kappa^{in}([\rho,1) \times \partial M^{in})$. From the previous argument, we know that all limits $x^{k,\rho}$ lie in the interior of $M^{in,\rho}$. Suppose first that $u$ intersects $\partial M^{in,\rho}$ transversally, so that the preimage $u^{-1}(M^{out} \cup \kappa^{in}([\rho,1) \times \partial M^{in}))$ is a compact surface with corners. Applying Lemma \ref{th:penetrate} then shows that $u(S) \subset M^{in,\rho}$. To remove the transversality assumption, one argues as before by moving the hypersurface slightly, then taking the limit. \qed

\section{Transversality}

We now revisit in more detail the transversality claims made throughout the main part of the paper. The proofs are basically routine, in spite of the presence of a finite symmetry group (see Remark \ref{th:equi} for an informal explanation of why this group does not cause any problems).

\subsection{Nondegeneracy of chords\label{subsec:chords}}

Let $M$ be a compact symplectic manifold possibly with boundary, $L \subset M$ a Lagrangian submanifold, and $H \in \smooth(M,\R)$ a function, with its associated Hamiltonian vector fields $X$ and flow $\phi$. We will be interested in the generic behaviour of integer $X$-chords. Here, generic is understood in the sense of Baire, with respect to perturbations of $H$. Actually, we will sometimes find it useful to perturb $L$ as well, in a Hamiltonian (exact Lagrangian) way. However, as far as $X$-chords are concerned, passing from $L$ to $\psi(L)$ is equivalent to going from $H$ to $H \circ \psi^{-1}$ while keeping $L$ fixed. Hence, all the results ultimately translate into statements about perturbing $H$ alone, for arbitrary fixed $L$.

As a first step, after a perturbation of $H$, we may assume that: first, at any stationary point $x$, the derivative $D\phi^1_x$ has no eigenvalues which are roots of unity; secondly, if $x$ is a periodic point whose minimal period is some integer $w>0$, then $D\phi^w_x$ has $X$ as the only eigenvector for the eigenvalue $1$, and no other eigenvalues which are roots of unity. As a consequence of these properties, if $x$ belongs to any (possibly stationary of multiply covered) periodic orbit of integer period $w$, the image of
\begin{equation} \label{eq:1-minus-dphi}
Id - D\phi^w_x: TM_x \longrightarrow TM_x,
\end{equation}
is precisely the symplectic orthogonal complement of $X_x$. In addition, we can assume that the stationary points are disjoint from $L$, and that the non-stationary integer periodic orbits intersect $L$ transversally. Fix an $H$ with these properties. We want to know if, by a small generic Hamiltonian isotopy of $L$, one can make the intersections $\phi^w(L) \cap L$ transverse. By a standard Sard-theory argument, this is possible if for every $x \in \phi^w(L) \cap L$, the map
\begin{equation}
\begin{aligned} \label{eq:difference}
\mathit{Def}(L) & \longrightarrow TM_x/(TL_x + T(\phi^w(L))_x), \\
\xi & \longmapsto \xi_x - D\phi^w(\xi_{\phi^{-w}(x)})
\end{aligned}
\end{equation}
is surjective. Here, $\mathit{Def}(L)$ is the space of infinitesimal deformations of $L$, seen as normal vector fields. If $x \neq \phi^{-w}(x)$, the map \eqref{eq:difference} subtracts values of $\xi$ at two different points, so surjectivity is automatic. In the remaining case, we have $X_x \notin TL_x$ by assumption, which implies that $TL_x$ is not contained in the symplectic orthogonal complement of $X_x$. From this and the previous remark concerning \eqref{eq:1-minus-dphi}, it follows that \eqref{eq:difference} is still onto. Hence:

\begin{lemma} \label{th:nondegenerate-g}
For generic $H$, all integer $X$-chords are nondegenerate. \qed
\end{lemma}

Next, consider triple intersections $L \cap \phi^w(L) \cap \phi^v(L)$, for integers $0 < w < v$. In that case, the map corresponding to \eqref{eq:difference} is
\begin{equation}
\begin{aligned}
\mathit{Def}(L) \oplus TM_x &
\longrightarrow
TM_x/TL_x \oplus TM_x/T(\phi^w(L))_x \oplus TM_x/T(\phi^v(L))_x, \\
(\xi,\delta) & \longmapsto (\xi_x+\delta,D\phi^w(\xi_{\phi^{-w}(x)}) +\delta,
D\phi^v(\xi_{\phi^{-v}(x)}) + \delta).
\end{aligned}
\end{equation}
This is obviously surjective if the points $(x,\phi^w(x),\phi^v(x))$ are all different. If exactly two of them coincide, it is surjective by the previous argument. In the remaining situation, where all three points coincide, surjectivity fails. This happens whenever $L$ intersects an integer periodic orbit of $X$, which doesn't happen in dimensions $4$ and higher due to our previous transverse intersection assumption. The upshot is:

\begin{lemma} \label{th:nonconcatenate-g}
Assume that $\mathit{dim}\,M \geq 4$. For generic $H$, starting points of integer $X$-chords are never endpoints of integral $X$-chords. \qed
\end{lemma}

%The same arguments apply to the slightly more complicated setup where $M$ is a Weinstein domain, %$L$ a Lagrangian submanifold satisfying \eqref{eq:z-parallel} as well as \eqref{eq:no-integer-chords}, %and one only allows Hamiltonian functions constructed as in \eqref{eq:class-of-h}.

\subsection{Trivial solutions\label{subsec:trivial-solutions}}
Consider an equation \eqref{eq:restate-dbar}, in the setup described there. A {\em trivial solution} of that equation is a map $u$ such that
\begin{equation} \label{eq:du-equal-x}
du = X \otimes \gamma.
\end{equation}
These are precisely the solutions whose geometric energy is zero, or more intuitively, those maps which satisfy \eqref{eq:restate-dbar} for all possible $J$ at the same time. In local coordinates as in \eqref{eq:local-floer}, this means that
\begin{equation} \label{eq:local-trivial}
\partial_s u - p(s,t)X = 0, \quad
\partial_t u - q(s,t)X = 0.
\end{equation}
The simplest examples are constant maps $u(z) = x$, where $X_x = 0$.

\begin{lemma}
Suppose that $d\gamma$ is not identically zero. The only trivial solutions are the constant ones taking values at stationary points of $X$.
\end{lemma}

\proof This follows from the integrability condition for \eqref{eq:local-trivial}. Namely, taking $\nabla$ to be any torsion-free connection on $TM$, we have
\begin{equation}
\begin{aligned}
0 & = \nabla_{\partial_s u} \partial_t u - \nabla_{\partial_t u} \partial_s u \\
 & = q \cdot \nabla_{\partial_s u} X + \partial_sq \cdot X
 - p \cdot \nabla_{\partial_t u}X - \partial_tp \cdot X \\
 & = q \cdot \nabla_{pX} X - p \cdot \nabla_{qX}X + (\partial_s q- \partial_t p) \cdot X
 = d\gamma(\partial_s,\partial_t) \cdot X.
\end{aligned}
\end{equation}
Since $d\gamma$ is not identically zero by assumption, there is some $z \in S$ such that $u(z)$ is a stationary point of $X$. But then, \eqref{eq:du-equal-x} and elementary ODE theory show that $u$ must be constant.
\qed

\begin{lemma} \label{th:z-only}
Suppose that $d\gamma = 0$, and that $H$ satisfies \eqref{eq:nonconcatenate-x}. The only trivial solutions are the ones with domain $S \iso Z$.
\end{lemma}

\proof Since $S$ is simply-connected we can write $\gamma = dl$, and then integrate \eqref{eq:du-equal-x}, which yields
\begin{equation} \label{eq:integrate-out}
u(z) = \phi^{l(z)}(x)
\end{equation}
for some $x$ which is independent of $z$. Equivalently, one could view this as the result of applying a transformation \eqref{eq:ansatz} globally, and getting a $\tilde{u}$ which is a holomorphic curve with energy zero. The function $l$ is locally constant on $\partial S$, and over the strip-like ends it satisfies
\begin{equation}
l(\epsilon^k(s,t)) = \begin{cases} w^0 t & k = 0, \\
w^1 + \cdots + w^{k-1} + w^k t & k > 0.
\end{cases}
\end{equation}
For the limits $\bfx = \{x^k\}$ this means that
\begin{equation}
\begin{aligned}
 & x^1(0) = x^0(0), \\
 & x^2(0) = x^1(1) = \phi^{w_1}(x^1(0)), \\
 & \cdots \\
 & x^d(0) = x^{d-1}(1) = \phi^{w_{d-1}}(x^{d-1}(0)), \\
 & x^1(1) = x^d(1).
\end{aligned}
\end{equation}
In particular, if $d > 1$ then $x^2(0)$ would be both an endpoint and a starting point of an integral $X$-chord, which contradicts \eqref{eq:nonconcatenate-x}. \qed

\begin{lemma} \label{th:discrete-1}
Suppose that $d\gamma = 0$, and let $u$ be a solution of \eqref{eq:restate-dbar} which is not a trivial solution. Then the set of points where $du = X \otimes \gamma$ holds is a discrete subset of $S$.
\end{lemma}

\proof Locally near each point of $S$ we can apply a transformation \eqref{eq:ansatz}, which takes \eqref{eq:restate-dbar} into an ordinary pseudo-holomorphic curve equation, and transforms trivial solution into constant maps. Using standard pseudo-holomorphic curve techniques (the Carleman similarity principle; see \cite[Corollary 2.3]{floer-hofer-salamon94} for interior points, and \cite[Corollary 2.2]{oh96c} for boundary points), one sees that the subset of $U$ consisting of points at which $du-X \otimes \gamma$ vanishes to infinite order is both open and closed. The rest is straightforward. \qed

\begin{lemma} \label{th:discrete-2}
Suppose that $d\gamma \neq 0$, and that $H$ satisfies \eqref{eq:nonconcatenate-x}. Let $U \subset S$ be a connected open subset, containing $\partial S$ as well as neighbourhoods of all the points at infinity, and such that $d\gamma|U = 0$ (such a subset always exists, because of the assumptions on $\gamma$). Then, for any solution $u$ of \eqref{eq:restate-dbar}, the set of points in $U$ where $du = X \otimes \gamma$ holds is a discrete subset of $U$.
\end{lemma}

\proof Suppose that the contrary is true. Applying the same argument as in Lemma \ref{th:discrete-1}, one then sees that $du = X \otimes \gamma$ holds on all of $U$. As in Lemma \ref{th:z-only}, it follows that the limits $\bfx$ of our map satisfy $x^1(0) = x^0(0)$. In view of \eqref{eq:stokes}, the assumption $d\gamma \neq 0$ implies that $w^0>w^1$. Therefore $x^1(1) = \phi^{w^1}(x^1(0)) = \phi^{w^1-w^0}(x^0(1))$ is both starting point and endpoint of an integral $X$-chord, contradicting the other assumption. \qed

\subsection{Linearization\label{subsec:linear}}
We now consider the linearization of \eqref{eq:restate-dbar} with respect to infinitesimal changes of $\bfJ$ as well as as $u$. This is given by an operator
\begin{equation} \label{eq:d2}
\begin{aligned}
 \scrY \oplus W^{1,r}(u^*TM) & \longrightarrow L^r(Hom^{0,1}(TS,u^*TM)), \\
 (\bfY,X) & \longmapsto \half Y_{z,u(z)} \circ (du_z-X_{u(z)} \otimes \gamma_z) \circ j_z + D_uX.
\end{aligned}
\end{equation}
Here,  $\scrY$ is the space of infinitesimal deformations of our given family of almost complex structure, which have superexponential decay over the ends of $S$. The second component is the standard linearized $\bar\partial$-operator $D_u$, taking the Sobolev completion for some $r>2$.

\begin{lemma} \label{th:surjective}
Suppose that $u$ satisfies \eqref{eq:restate-dbar}, and is not a trivial solution of that equation.
Then the map \eqref{eq:d2} is surjective.
\end{lemma}

\proof This is a standard argument. Since $D_u$ itself is Fredholm, we only need to show that the image of \eqref{eq:d2} has zero orthogonal complement. Suppose that $T \neq 0$ lies in the kernel of the adjoint $D_u^*$. Then it is automatically smooth, and vanishes only on a discrete subset of points in $S$. If $d\gamma = 0$, one can use Lemma \ref{th:discrete-1} to find a point $z \in S$ where both $T$ and $du-X \otimes \gamma$ are nonzero, and a $\bfY$ such that
\begin{equation} \label{eq:nonzero-pairing}
\leftsc Y_{z,u(z)} \circ (du_z-X_{u(z)} \otimes \gamma_z) \circ j_z, T_z \rightsc \neq 0.
\end{equation}
Strictly speaking, our almost complex structures are constrained to be of contact type along the boundary. Hence, to ensure that \eqref{eq:nonzero-pairing} can be achieved, we need $u(z)$ to lie in the interior of $M$. But since we are assuming that \eqref{eq:no-integer-chords} holds, all limits $x$ lie in the interior, hence so does $u(z)$ outside a compact subset of $S$. Having found \eqref{eq:nonzero-pairing}, one multiplies $\bfY$ with a function supported very close to $z$, and gets a contradiction to the assumption that $T$ lies in the orthogonal complement to the image of \eqref{eq:d2}. In the case where $d\gamma \neq 0$, the same argument works except that one has to take $z$ in a subset $U \subset S$ of the kind considered in Lemma \ref{th:discrete-2}. \qed

Now consider two solutions $u,v$ of \eqref{eq:restate-dbar}, with possibly different limits at infinity. We can then consider the linearization of the equation at both solutions, with respect to a common variation of the family of almost complex structures. In parallel with \eqref{eq:d2}, this is given by an operator
\begin{equation} \label{eq:d-double}
\begin{aligned}
 & \scrY \oplus W^{1,r}(u^*TM) \oplus W^{1,r}(v^*TM) \longrightarrow \\
 & \qquad \rightarrow L^r(Hom^{0,1}(TS,u^*TM)) \oplus L^r(Hom^{0,1}(TS,v^*TM)). \\
\end{aligned}
\end{equation}
One can now argue as before, varying the almost complex structure near a point $z \in S$ such that $(du - X \otimes \gamma)_z \neq 0$, $(dv - X \otimes \gamma)_z \neq 0$, and $u(z) \neq v(z)$.
With the obvious exceptions, the existence of sufficiently many points with that property is ensured by Lemma \ref{th:intersection}. The outcome is:

\begin{lemma} \label{th:double-surjective}
Suppose that neither $u$ nor $v$ are trivial solutions, and that $u \neq v$. Then \eqref{eq:d-double} is surjective.
\end{lemma}

\begin{addendum} \label{th:partial-restriction}
Inspection of the proofs above shows that one can restrict to the subspace of those $\bfY$ which vanish outside some open subset of $M$, as long as that subset contains all integer $X$-chords.
\end{addendum}

\subsection{Application\label{subsec:regular}}
We will now explain how these arguments apply to the moduli spaces $\scrR^{d+1,\bfp,\bfw}(\bfx)$. It is a classical fact that a generic choice of $\bfI_w$ ensures that the moduli spaces of solutions to Floer's equation \eqref{eq:old-floer} are regular. With that done, let's turn to the case of stable weighted popsicles $(S,\bfphi,\bfw)$. Recall that in \eqref{eq:dbar} we use families of almost complex structures obtained by taking a fixed universal choice $\bfI = \bfI_{S,\bfphi,\bfw}$ and then changing it in the way indicated by an infinitesimal deformation $\bfK = \bfK_{S,\bfphi,\bfw}$. The latter is restricted by $\Sym^\bfp$-invariance as well as the asymptotic consistency requirement. Take a dense Banach subspace $\scrK$ of the space of all possible choices of infinitesimal deformations, and consider the universal moduli space, which is the space of solutions of \eqref{eq:dbar} with $\bfK$ considered as an additional variable. The linearization of this equation is a map
\begin{equation} \label{eq:d3}
\scrK \times T\scrR^{d+1,\bfp,\bfw} \times
W^{1,r}(u^*TM) \longrightarrow L^r(Hom^{0,1}(TS,u^*TM)),
\end{equation}
whose first and third components are as in \eqref{eq:d2}, and where the tangent space to $\scrR^{d+1,\bfp,\bfw}$ is taken at the point corresponding to $(S,\bfphi)$. The implications of asymptotic consistency were spelled out in Remark \ref{th:concrete-consistency}. Inspection of that discussion shows that for any given $S$ the choice of $\bfK$ is free, subject only to the conditions of superexponential decay over the ends. It therefore follows from Lemma \ref{th:surjective} (with $\bfY$ restricted to a dense subspace of $\scrY$, which does not affect the result) that \eqref{eq:d3} is onto. By a standard Sard-Smale argument, a generic choice of $\bfK \in \scrK$ will ensure that for every $(S,\bfphi,u) \in \scrR^{d+1,\bfp,\bfw}(\bfx)$ the restriction of \eqref{eq:d3} to fixed almost complex structure, which is an operator
\begin{equation} \label{eq:parametrized-d}
T\scrR^{d+1,\bfp,\bfw} \times
W^{1,r}(u^*TM) \longrightarrow L^r(Hom^{0,1}(TS,u^*TM)),
\end{equation}
is again onto. This is precisely the statement of Theorem \ref{th:smoothness-1}, and the kernel of \eqref{eq:parametrized-d} is the tangent space of our moduli space. The same argument, applied to the stratification of $\scrR^{d+1,\bfp,\bfw}$ by isotropy groups of the $\Sym^\bfp$-action, proves Theorem \ref{th:smoothness-2}.

\begin{remark} \label{th:equi}
In view of the generally problematic nature of equivariant transversality, it behoves us to explain why in this specific instance, the action of $\Sym^\bfp$ does not pose any problems. As a toy model, suppose we have a manifold $B$ and a vector bundle $E$ over it, both of which come with compatible actions of a finite group $\Gamma$. Crucially, suppose that {\em for every point $b \in B$, the isotropy group $\Gamma^b$ acts trivially on $E_b$}. In that case, for every given $\eta \in E_b$ one can find an equivariant section $s$ such that $s(b) = \eta$, simply by averaging. In other words, the evaluation map on the space of equivariant sections,
\begin{equation}
B \times \smooth(B,E)^\Gamma \longrightarrow E
\end{equation}
is a submersion. Therefore, a generic $s \in \smooth(B,E)^\Gamma$ is transverse to the zero-section. Returning to our real-world case, the assumption on isotropy groups corresponds to the fact that if $(S,\bfphi)$ lies in the subset of $\scrR^{d+1,\bfp}$ fixed by some subgroup of $\Sym^\bfp$, then that subgroup acts trivially on the range of \eqref{eq:parametrized-d} for any map $u: S \rightarrow M$ satisfying \eqref{eq:restate-dbar}.
\end{remark}

To conclude, we consider parametrized moduli spaces as defined in \eqref{eq:parametrized-1}. As before, we can allow variations of the associated interpolating family of almost complex structures (taking suitable technical precautions, so that the variations form a Banach space), leading to an infinite-dimensional universal moduli space, which will be a smooth Banach manifold with boundary. Moreover, projection to the parameter is a submersion from the universal moduli space to $(0,1]$. This follows from Lemma \ref{th:surjective} and Addendum \ref{th:partial-restriction}, where the latter makes up for the restrictions imposed on the behaviour of interpolating families for small $\rho$. As before, this implies regularity of the actual parameter moduli space for generic choices of the interpolating family, which is the first part of Theorem \ref{th:smoothness-3}.

Now consider the product of two parametrized moduli spaces, and its universal version (where the almost complex structures on both factors are varied in the same way). Lemma \ref{th:double-surjective} shows that the universal space is smooth, and that the projection map to $(0,1]^2$ is a submersion away from the diagonal. This implies the second part of Theorem \ref{th:smoothness-3}, and the final part is similar. Of course, these arguments strictly speaking only apply in the range $d+|F| \geq 2$ where the underlying popsicles are stable, but the corresponding results in ordinary Floer theory are well-known.

\section{Signs\label{sec:signs}}

The issues involved in constructing orientations of moduli spaces in Floer type theories, and the role played by $Pin$ structures in the open string case, are well-understood \cite{floer-hofer93, fooo, desilva00}. Therefore, we only give an outline the general theory, and then concentrate on two phenomena where contributions coming from different points in the moduli space cancel. The first of these comes from the action of $\Sym^\bfp$ on popsicle moduli spaces. The second one occurs in the context of cascades, where the same boundary component appears several times, see Example \ref{th:two-boundary}. We will also discuss the origin of the signs appearing in various of our formulae. This mainly comes down to clarifying the conventions, and the choices of orientations for $\scrR^{d+1,\bfp}$.

\subsection{Degrees and orientation spaces\label{subsec:grading-and-stuff}}
Let $L \subset M$ be a Lagrangian submanifold satisfying \eqref{eq:standard-floer}. One can then associate to any $x \in \scrX_w$ a Maslov index $\deg(x) \in \Z$. This is easiest to explain if one assumes that $[x] \in \pi_1(M,L)$ is trivial, which is when $x$ can be filled in by a half-disc. More precisely, we want to think of this half-disc as a map $w: H \rightarrow M$, where $H = \R \times [0,\infty) \subset \C$ is the upper half-plane, satisfying $w(\partial H) \subset L$ and $\lim_{s \rightarrow - \infty} w(se^{-i\pi t}) = x(t)$. Such a $w$ comes with a linear $\bar\partial$-operator $D_w$, and one defines
\begin{equation} \label{eq:degree}
\deg(x) = \mathrm{index}\,D_w.
\end{equation}
Vanishing of $2c_1(M,L)$ ensures that this integer is independent of the choice of $w$. In fact, one can allow higher genus surfaces as well and thereby extend the treatment to the case where only the homology class $[x] \in H_1(M,L)$ is trivial. However, the answer in full generality is not as canonical and requires a slightly different approach. Fix a compatible almost complex structure on $M$, and let
\begin{equation}
\scrK^{-2} = (\Lambda^{top}_{\C} TM)^{\otimes 2}
\end{equation}
be the bundle with first Chern class $2c_1(M)$. Along $L$ this bundle has a preferred trivialization, given by taking an orthonormal basis $\{\xi_k\}$ of $TL$ at a point, and considering the nonzero element $(\xi_1 \wedge \cdots \wedge \xi_n)^{\otimes 2}$ in the fibre of $\scrK^{-2}$ at that point, which is independent of the choice of basis. The vanishing of $2c_1(M,L)$ allows us to equip $L$ with a {\em grading}, which is an extension of the given trivialization to the whole of $M$. The choice of grading in turn singles out, for each $x \in \scrX_w$, a class of $\bar\partial$-operators on $H$ which are abstract generalizations of the $D_w$ introduced above. Those operators, called {\em orientation operators}, are then used to extend \eqref{eq:degree} (one reference is \cite[Section 11]{seidel04}, but the idea goes back to \cite{floer88}). Changing the homotopy class of the grading by some $c \in H^1(M,L;\Z)$ will affect the degree of each $x$ by $\leftsc c, [x] \rightsc$, which explains why for trivial $[x]$ one can get away without choosing a grading.

Next, we want to associate to each $x \in \scrX_w$ a one-dimensional real vector space $o_x$, the so-called {\em orientation space}. If one considers only those $x$ that can be filled in by a half-disc $w$, and assumes additionally that $H^1(M;\Z/2) \rightarrow H^1(L;\Z/2)$ is onto, one can simply define
\begin{equation} \label{eq:o-x}
o_x = det\, D_w
\end{equation}
to be the determinant line of the operator $D_w$. The additional assumption we just made ensures that if we have another choice of half-disc $\tilde{w}$, the loop formed by the boundaries of $w$ and $\tilde{w}$ is contractible in $L$. A choice of contraction of that loop yields an isomorphism $det\,D_w \iso det\,D_{\tilde{w}}$; and since $w_2(L) = 0$, that isomorphism is independent of the specific contraction, up to multiplication with a positive constant. The general procedure is slightly different: one chooses a $Pin$ structure on $L$, and then replaces $D_w$ by abstract orientation operators which are compatible with that structure in a suitable way. Twisting the $Pin$ structure by a real line bundle $\lambda$ has the effect of changing each $o_x$ to $o_x \otimes \lambda_{x(0)} \otimes \lambda_{x(1)}$.

We assume from now that a grading and $Pin$ structure have been chosen, so that $\deg(x)$ and $o_x$ are unambiguously defined for all $x$. Additionally, we need the following notation. Let $\xi$ be a one-dimensional real vector space, and denote by $\sigma^{\pm}$ its two possible orientations. For a fixed coefficient field $\K$, the {\em $\K$-normalization} of $\xi$ is the one-dimensional $\K$-vector space
\begin{equation}
|\xi|_{\K} = \K \sigma^+ \oplus \K \sigma^- / \K( \sigma^+ + \sigma^-).
\end{equation}
By construction, $|\xi|_{\K}$ is identified with $\K$ in a way which is canonical up to sign. A choice of identification is precisely the same as an orientation of $\xi$. Similarly, any isomorphism between one-dimensional real vector spaces induces an isomorphism of their $\K$-normalizations, which is $\pm 1$ if one chooses identifications of those spaces with $\K$.

\subsection{Zero-dimensional moduli spaces\label{subsec:sign-zero}}
Consider a point $(S,\bfphi,u)$ in a moduli space $\scrR^{d+1,\bfp,\bfw}(\bfx)$. If we take the linearized operator $D_u$ associated to $u$, and glue it at the positive ends to orientation operators for $x^1,\dots,x^d$, the outcome is an operator on $H$ which is homotopic to an orientation operator for $x^0$. The gluing relation for determinant lines of $\bar\partial$-operators therefore yields an isomorphism
\begin{equation} \label{eq:composition-law}
o_{x^o} \iso det\,D_u \otimes o_{x^d} \otimes \cdots \otimes o_{x^1}.
\end{equation}
To be more precise, one should say that by choosing the homotopy of operators to be compatible with the given $Pin$ structures, one gets an isomorphism of one-dimensional vector spaces which is canonical up to multiplication by a positive constant (the same interpretation should be applied to similar statements made later on). Assume from now on that all moduli spaces are regular, as in Theorem \ref{th:smoothness-1}, and also that $d+|F| \geq 2$, which means that the underlying popsicle $(S,\bfphi)$ is stable. Then, because the tangent space to the moduli space is the kernel of an operator of the form \eqref{eq:parametrized-d}, we have
\begin{equation} \label{eq:decompose-t}
\Lambda^{top}(T\scrR^{d+1,\bfp,\bfw}(\bfx))_{S,\bfphi,u} \iso
\Lambda^{top}(T\scrR^{d+1,\bfp})_{S,\bfphi} \otimes det\,D_u,
\end{equation}
Plugging those two relations into each other, we get an isomorphism
\begin{equation} \label{eq:plug}
o_{S,\bfphi,u} : \Lambda^{top}(T\scrR^{d+1,\bfp,\bfw}(\bfx))_{S,\bfphi,u} \otimes o_{x^d} \otimes \cdots \otimes o_{x^1} \iso \Lambda^{top}(T\scrR^{d+1,\bfp})_{S,\bfphi} \otimes o_{x^0}.
\end{equation}
At this point, let's choose an orientation of each moduli space $\scrR^{d+1,\bfp}$ (which is possible since they are contractible). This allows us to trivialize the $T\scrR$ term in \eqref{eq:plug}. Assume additionally that the moduli space of popsicle maps is zero-dimensional, so that $\Lambda^{top}(T\scrR^{d+1,\bfp,\bfw}(\bfx))$ is canonically trivial. We then get an isomorphism
\begin{equation} \label{eq:reduced-o}
o_{S,\bfphi,u}^{red}: o_{x^d} \otimes \cdots \otimes o_{x^1} \iso o_{x^0}.
\end{equation}
Since this is unique up to multiplication by a positive constant, the induced map on $\K$-normalizations, denoted by $|o_{S,\bfphi,u}^{red}|_\K$, is canonically well-defined. By definition, this is the contribution \eqref{eq:orientation-map} of $(S,\bfphi,u)$ to the algebraic count of points in our moduli space.

\begin{lemma} \label{th:cancellation}
Take $\sigma \in \Sym^{\bfp} \subset \Sym^{|F|}$, acting on the zero-dimensional moduli space $\scrR^{d+1,\bfp,\bfw}(\bfx)$ by permuting popsicle maps, $(S,\bfphi,u) \mapsto (S,\sigma(\bfphi),u)$. Then, up to a positive constant,
\begin{equation}
o_{S,\sigma(\bfphi),u}^{red} = \mathrm{sign}(\sigma)\, o_{S,\bfphi,u}^{red}.
\end{equation}
\end{lemma}

This is clear from the construction above. Since $\scrR^{d+1,\bfp}$ projects to $\scrR^{d+1}$ with fibers $\bR^F$, and  $\Sym^{\bfp} \subset \Sym^{F}$ acts fibrewise by the standard representation of permutations on $|F|$ letters,  $\sigma$ preserves orientations if and only if $\sigma$ is an even permutation. Hence, if we compare $(S,\sigma(\bfphi),u)$ with $(S,\sigma,u)$, the associated maps \eqref{eq:composition-law} agree, but the trivializations $det\, D_u \iso \R$ obtained from \eqref{eq:decompose-t} differ by $\mathrm{sign}(\sigma)$, and so do the maps \eqref{eq:reduced-o}.

\subsection{One-dimensional moduli spaces\label{subsec:1-dim}}

Consider a point in the boundary of a one-dimensional moduli space $\bar\scrR^{d+1,\bfp,\bfw}(\bfx)$, represented by a broken map with two components $(S_\pm,\bfphi_\pm,u_\pm)$. As before, we assume stability of the underlying popsicles, which means that $(S_\pm,\bfphi_\pm,u_\pm) \in \scrR^{d_\pm + 1, \bfp_\pm, \bfw_\pm}(\bfx_\pm)$ with $d_\pm + |F_\pm| \geq 2$. Gluing them together, one gets a nearby point $(S,\bfphi,u)$ in the interior $\scrR^{d+1,\bfp,\bfw}(\bfx)$. The gluing parameter yields a preferred local orientation for the moduli space (pointing towards the boundary), hence an isomorphism $\Lambda^{top}(T\scrR^{d+1,\bfp,\bfw}(\bfx)) \iso \R$ at $(S,\bfphi,u)$. With that in mind, and an additional choice an orientation of $\scrR^{d+1,\bfp}$, one can turn $o_{S,\bfphi,u}$ into a map $o_{S,\bfphi,u}^{red}$ of the same form as in \eqref{eq:reduced-o}. The elementary fact that boundary points come in pairs implies that the sum of the normalizations $|o_{S,\bfphi,u}^{red}|_\K$ over all $(S_\pm,\bfphi_\pm,u_\pm)$ vanishes. Of course, this is not strictly speaking true, since we are ignoring the case where a Floer trajectory bubbles off, but we will make up for that omission later (Section \ref{subsec:unstable}).

On the other hand, we shall presently see that $o_{S,\bfphi,u}^{red}$ can be decomposed as follows:
\begin{equation} \label{eq:in-steps}
\begin{aligned}
 & o_{x^d} \otimes \cdots \otimes o_{x^1} \\
 & \iso o_{x^d} \otimes \cdots \otimes o_{x^{i+d_-}} \otimes \Lambda^{top}(T\scrR^{d_- + 1,\bfp_-})_{S_-,\bfphi_-} \otimes o_{x^{\mathit{new}}}
 \otimes o_{x^{i-1}} \otimes \cdots \otimes o_{x^1} \\
 & \iso \Lambda^{top}(T\scrR^{d_- + 1,\bfp_-})_{S_-,\bfphi_-} \otimes o_{x^d} \otimes \cdots \otimes o_{x^{\mathit{new}}} \otimes \cdots \otimes o_{x^1} \\
 & \iso \Lambda^{top}(T\scrR^{d_+ + 1, \bfp_+})_{S_+,\bfphi_+} \otimes \Lambda^{top}(T\scrR^{d_- + 1,\bfp_-})_{S_-,\bfphi_-} \otimes
 o_{x^0} \\
 & \iso \Lambda^{top}(T\scrR^{d+1,\bfp})_{S,\bfphi} \otimes o_{x^0} \\
 & \iso o_{x^0}.
\end{aligned}
\end{equation}
The first step consists of applying $o_{S_-,\bfphi_-,u_-}$, keeping in mind that $\scrR^{d_-+1,\bfp_-,\bfw_-}(\bfx_-)$ is necessarily of dimension zero; here, $x^{\mathit{new}} = x^0_- = x^i_+$ is the $X$-chord at the ends where our two components are glued together. The second step is just a permutation of tensor factors (with suitable Koszul signs given by dimensions and degrees; this is natural, since each factor encodes orientations of a vector space). Thirdly, we apply $o_{S_+,\bfphi_+,u_+}$, and exchange the first two factors of the resulting expression. Next, recall that $(S,\bfphi)$ is a point of $\scrR^{d+1,\bfp}$ close to the boundary stratum to which the broken popsicle $(S_\pm,\bfphi_\pm)$ belongs. Using any collar neighbourhood of that stratum, we get an isomorphism
\begin{equation} \label{eq:boundary-1}
\Lambda^{top}(T\scrR^{d+1,\bfp})_{S,\bfphi} \iso \Lambda^{top}(T\scrR^{d_++1,\bfp_+})_{S_+,\bfphi_+} \otimes \Lambda^{top}(T\scrR^{d_-+1,\bfp_-})_{S_-,\bfphi_-},
\end{equation}
which is used in the fourth step in \eqref{eq:in-steps}. The fifth and last one is just to put in the chosen orientation of $\scrR^{d+1,\bfp}$. After unwinding the definition of \eqref{eq:plug}, the fact that the sequence of operations in Equation \eqref{eq:in-steps} agrees with the original definition of $o_{S,\bfphi,u}^{red}$ reduces to applying the gluing isomorphism $det\,D_u \iso det\,D_{u_+} \otimes det\,D_{u_-}$ and using its associativity properties.

From now on, we always assume that $\Sym^\bfp$ is trivial, and accordingly think of $F \subset \{1,\dots,d\}$ as in Section \ref{subsec:algebraic-relation}. The same property then holds for $\bfp_-$, but not necessarily for $\bfp_+$, compare Example \ref{th:ex-popsicle}(ii). Given $\sigma \in \Sym^{\bfp_+}$, let $(\tilde{S},\tilde{\bfphi},\tilde{u})$ be the element of $\scrR^{d+1,\bfp,\bfw}(\bfx)$ obtained by gluing together $(S_+,\sigma(\bfphi_+),u_+)$ and $(S_-,\bfphi_-,u_-)$. We can see the effect of $\sigma$ on \eqref{eq:in-steps} by looking at the diagram of isomorphisms
\begin{equation}
\xymatrix{
 & \hspace{-2em} o_{x^d} \otimes \cdots \otimes o_{x^1} \ar[dl] \ar[dr] \hspace{-2em}
 %\ar@/_10pc/[ddl]
 %\ar[ddr]
 & \\
 {
 \begin{matrix}
 \Lambda^{top}(T\scrR^{d_+ + 1, \bfp_+})_{(S_+,\bfphi_+)} \\
 \otimes \, \Lambda^{top}(T\scrR^{d_- + 1,\bfp_-})_{(S_-,\bfphi_-)} \\
 \otimes \, o_{x^0}
 \end{matrix}
 }
 \ar[rr]^-{\Lambda^{top}(D\sigma) \otimes id \otimes id}
 \ar[d] && \ar[d]
 {
 \begin{matrix}
 \Lambda^{top}(T\scrR^{d_+ + 1,\bfp_+})_{(S_+,\sigma(\bfphi_+))} \\
 \otimes \, \Lambda^{top}(T\scrR^{d_- + 1,\bfp_-})_{(S_-,\bfphi_-)} \\
 \otimes \, o_{x^0}
 \end{matrix}
 }
 \\
 \Lambda^{top}(T\scrR^{d+1,\bfp})_{(S,\bfphi)} \otimes o_{x^0} \ar[dr] &&
 \Lambda^{top}(T\scrR^{d+1,\bfp})_{(\tilde{S},\tilde{\bfphi})} \otimes o_{x^0} \ar[dl]
 \\
 & o_{x^0}
}
\end{equation}
Going down the left or right sides yields the contributions of $(S,\bfphi,u)$ and $(\tilde{S},\tilde{\bfphi},\tilde{u})$ to the total count of boundary points, respectively. The upper triangle of the diagram commutes. The lower pentagon depends only on the orientations of the various popsicle moduli spaces, and can be described as follows. On the left side, we compare some orientation of the boundary stratum \eqref{eq:boundary-1} with the chosen orientation of the interior $\scrR^{d+1,\bfp}$. On the right hand side we do the same, but where the orientation of the boundary stratum is changed by the action of $\sigma$, so the difference is $\mathrm{sign}(\sigma)$. In particular, whenever $\mathrm{sign}(\sigma) = -1$, the contributions of $(S,\bfphi,u)$ and $(\tilde{S},\tilde{\bfphi},\tilde{u})$ cancel out.

With that out of the way, we turn to the case where $\Sym^{\bfp_\pm}$ are both trivial, whose combinatorics is captured by the notion of admissible cut (Definition \ref{th:cut}). In that case, the outcome of our discussion can be summarized as follows. Given a choice of orientation of all the popsicle moduli spaces, let $(-1)^\triangle$ be the sign measuring the difference between the product orientation of \eqref{eq:boundary-1} and that of $\scrR^{d+1,\bfp}$. Then what \eqref{eq:in-steps} says is that
\begin{equation} \label{eq:simplified-steps}
\begin{aligned}
& o_{S,\bfphi,u}^{red} = (-1)^{\S}\, o_{S_+,\bfphi_+,u_+}^{red} \circ (id^{\otimes d_+ - i} \otimes
o_{S_-,\bfphi_-,u_-}^{red} \otimes id^{\otimes i-1}), \\
& \S = \triangle + (\dim\,\scrR^{d_- + 1,F_-})(\dim\,\scrR^{d_+ + 1,F_+}) \\ & \qquad \qquad +
(\dim\,\scrR^{d_- + 1,F_-})(\deg\,x^{i+d_-} + \cdots + \deg\, x^d).
\end{aligned}
\end{equation}

\subsection{Explicit orientations\label{subsec:choose-orientation}}
So far, the choice of orientations of the popsicle moduli spaces has been arbitrary. For those cases where $\Sym^{\bfp}$ is nontrivial, there is no need to restrict that freedom, since their contributions cancel out anyway. In contrast, for the cases with trivial symmetry, we do want to make an explicit choice, so as to get a concrete expression for $\triangle$ in \eqref{eq:simplified-steps}.

The standard way of choosing an orientation of $\scrR^{d+1}$ is as follows. For every pointed disc $S$ there is a unique isomorphism $\bar{S} \iso D$ which takes $\zeta^0,\zeta^1,\zeta^2$ to $-1,-i,+1 \in \partial S$, and hence all other points at infinity $\zeta^k$ to points in the upper semicircle. This yields an identification
\begin{equation} \label{eq:orient-stasheff}
\scrR^{d+1} \iso \{0 < t_1 < \dots <t_{d-2} < \pi\} \subset \R^{d-2},
\end{equation}
and one pulls back the standard orientation from the right hand side. Next, consider the moduli spaces $\scrR^{d+1,F}$ for $d \geq 2$ and $F \subset \{1,\dots,d\}$. These come with diffeomorphisms to $\scrR^{d+1} \times \R^F$, obtained by splitting the forgetful map \eqref{eq:forgetful}. The splitting is not canonical, but the induced orientation is well-defined (to make things precise, we should say that the identification of the fibre of the forgetful map with $\R^F$ is such that increasing the coordinates has the effect of moving the sprinkles away from $\zeta^0 = -\infty$; moreover, the orientation of $\R^F$ is the one obtained by ordering the elements of $F$ in the obvious way). We extend this to the case where $d = 1$, $F = \{1\}$ by taking the {\em opposite} of the tautological orientation of $\scrR^{1,\{1\}} = point$. An elementary computation shows that with these conventions, the orientation difference associated to a given cut is
\begin{equation}
\begin{aligned}
\triangle = &\;
d_-d_+ + d_-i + i + 1 \\ & +|F_+| d_- + |\{(k_+,k_-) \in F_+ \times F_- \;:\; \iota_+(k_+) >
\iota_-(k_-)\}|.
%d_-d_+ + (d_- - 1)(i-1) \\ & +|F_+| d_- + |\{(k_+,k_-) \in F_+ \times F_- \;:\; \iota_+(k_+) >
%\iota_-(k_-)\}|.
\end{aligned}
\end{equation}
After substituting this, as well as the obvious formulae for the dimensions of the moduli spaces, into \eqref{eq:simplified-steps}, one gets signs as originally stated in \eqref{eq:orientation-signs}. In principle, any other choice of orientation is equally viable; the only difference is that the sign in \eqref{eq:ad-hoc-sign} would have to be adjusted accordingly, so as to recover the $A_\infty$-associativity relations in the desired form \eqref{eq:as}.

\subsection{The unstable case\label{subsec:unstable}}
We will now describe a slightly different version of the formalism above, for $d = 1$, which is more easily extended to $F = \emptyset$. Take $(S,\bfphi,u) \in \scrR^{2,\bfp,\bfw}(\bfx)$. As in Section \ref{subsec:chain} we identify $S = Z$, so that $(\bfphi,u)$ are maps defined on $Z$, determined up to translation. In that case, the tangent space of the moduli space can be identified with a kernel of an operator
\begin{equation} \label{eq:symmetric-operator}
\R^F \times
W^{1,r}(u^*TM) \longrightarrow \R \times L^r(u^*TM).
\end{equation}
Here, increasing any of the entries in $\bfsigma = \{\sigma_f\} \in \R^F$ has the effect of moving the corresponding sprinkle away from $-\infty$, while leaving the map $u$ constant. The second component of \eqref{eq:symmetric-operator} is the linearization of \eqref{eq:dbar}, with the standard trivialization of $TZ$ taken into account. The first component is the $L^2$ pairing with the infinitesimal generator of the translational $\R$-action, concretely
\begin{equation}
\begin{aligned}
& \R^F \times W^{1,r}(u^*TM) \longrightarrow \R, \\
& (\{\sigma_f\},X) \longmapsto -\textstyle \sum_f \sigma_f + \leftsc \partial_s u, X \rightsc.
\end{aligned}
\end{equation}
From \eqref{eq:symmetric-operator} we get an isomorphism
\begin{equation} \label{eq:r-r}
\Lambda^{top}(T\scrR^{2,\bfp,\bfw}(\bfx)) \otimes o_{x^1} \iso \R^\vee \otimes \Lambda^{top}(\R^F)
\otimes o_{x^0}.
\end{equation}
If we choose an ordering of $F$, hence an orientation of $\R^F$, the right hand side of \eqref{eq:r-r} can be identified with $o_{x^0}$. Using that, one then associates to any point in a zero-dimensional moduli space a map $o_{(Z,\bfphi,u)}^{red}: o_{x^1} \rightarrow o_{x^0}$. In the special case where $|F| = 1$, \eqref{eq:symmetric-operator} takes on the form
\begin{equation}
\begin{pmatrix} -id_\R & \leftsc \partial_s u, \cdot \rightsc \\ D_u(\partial_s u) & D_u \end{pmatrix},
\end{equation}
and one can deform it through invertible operators to $-id_\R \oplus D_u$. In this way, one sees that the sign contribution here agrees with the convention in Section \ref{subsec:choose-orientation}.

Now consider a point on the boundary of a one-dimensional space $\scrR^{2,\bfp,\bfw}(\bfx)$. In parallel with Section \ref{subsec:1-dim} we denote the components of that point by $(Z,\bfphi_\pm,u_\pm)$, and a nearby point in the interior by $(Z,\bfphi,u)$. The preferred orientation at the latter point, going towards the boundary, determines a map $o_{(Z,\bfphi,u)}^{red}: o_{x^1} \rightarrow o_{x^0}$. A gluing argument as in \eqref{eq:in-steps} shows that
\begin{equation}
\begin{aligned}
& o_{(Z,\bfphi,u)}^{red} = (-1)^\# o_{(Z,\bfphi_+,u_+)}^{red} \circ o_{(Z,\bfphi_-,u_-)}^{red}, \\
& \# = |F_+| + 1 + \{(f_+,f_-) \in F_+ \times F_- \;:\; f_+ > f_-\}.
\end{aligned}
\end{equation}
In this formula, the first summand comes from swapping the $\Lambda^{top}(\R^{F_+})$ factor with an adjacent $\R^\vee$, the second one arises because the standard generator of $\R^\vee$ moves the popsicle to the right, whereas in the gluing process, going to the boundary moves popsicles in opposite ways, and the last one comes from our chosen ordering of $F$ and the induced orderings of the subsets $F_\pm$. Of course, in practice only the case $|F| = 1$ is important, in which case the sign is just $(-1)^{|F_-|}$. This explains why in the two terms in Example \ref{th:algebraic-relation}(i) appear with opposite signs, thereby justifying the relevant case of \eqref{eq:orientation-signs}. A similar argument explains the other instances of \eqref{eq:algebraic-relation} where one of the two operations involved is $m^{1,\emptyset,\bfw}$.

\subsection{Parametrized moduli spaces\label{subsec:homotopy-signs}}
We now consider the parametrized moduli spaces and moduli spaces of cascades maps, always assuming that $d = 1$. The notation will generally follow that from Section \ref{subsec:chain}. In the parametrized case, variations of the parameter $\rho$ have to be included into the operator whose kernel describes the tangent space to the moduli space. Hence, the analogue of \eqref{eq:r-r} for a point $(\rho,Z,\bfphi,u) \in \scrP^{2,F,\bfw}(\bfx)$ is a canonical isomorphism
\begin{equation} \label{eq:parametrized-o}
\begin{aligned}
& \Lambda^{top}(T\scrP^{2,F,\bfw}(\bfx)) \otimes o_{x^1} \\
% & \iso \R \otimes \R^\vee \otimes \Lambda^{top}(\R^F) \otimes det\,D_u \otimes o_{x^1} \\
& \iso \R \otimes \R^\vee \otimes \Lambda^{top}(\R^F) \otimes o_{x^0},
\end{aligned}
\end{equation}
where the $\R$ factor is the tangent space to the parameter space $(0,1]$.

Now suppose that we have chosen an ordering of $F$, as well as identifications
$o_{x^1} \iso \R$, $o_{x^0} \iso \R$. After trivializing the $\R$ and $\R^\vee$
factors in \eqref{eq:parametrized-o} in the standard way, we then simply get
an orientation of the parametrized moduli space. Consider a codimension one boundary stratum in the partial compactification $\bar\scrP^{2,F,\bfw}(\bfx)$, which
is a fibre product
\begin{equation} \label{eq:fibre-product}
\scrP^{2,F_+,\bfw_+}(x^0,x^{\mathit{new}}) \times_{(0,1]} \scrP^{2,F_-,\bfw_-}(x^{\mathit{new}},x^1).
\end{equation}
After taking the orderings of $F_\pm$ induced from the one on $F$, and choosing
some identification $o_{x^{\mathit{new}}} \iso \R$, one gets orientations of both factors
in \eqref{eq:fibre-product}, hence an orientation of the fibre product. Note that
since $o_{x^{\mathit{new}}}$ appears in either factor, its trivialization does not really
affect the fibre product orientation.

\begin{remark} \label{th:fibre-product}
To explain the convention here, consider Euclidean spaces $\R^{n_\pm}$ with
coordinates $t_\pm$. Take projections to the first coordinate, $t_{\pm,1}: \R^{n_\pm} \rightarrow \R$. Then, the diffeomorphism
\begin{equation}
\begin{aligned}
& \R^{n_+} \times_\R \R^{n_-} \longrightarrow \R^{n_+ + n_- - 1}, \\
& (t_+,t_-) \longmapsto (t_{+,2},\dots,t_{+,n_+},t_-)
\end{aligned}
\end{equation}
takes the fibre product orientation to the standard orientation of $\R^{n_+ + n_- - 1}$.
In particular, if we consider the fibre product as the boundary of the open set $\{(t_+,t_-) \;:\; t_{+,1} < t_{-,1}\} \subset \R^{n_+} \times \R^{n_-}$, the induced boundary orientation agrees with the fibre product orientation.
\end{remark}

Again arguing in the style of \eqref{eq:in-steps}, one finds the following:

\begin{lemma} \label{th:compare-fibre-product}
The boundary orientation of the stratum \eqref{eq:fibre-product}, which is induced from that of $\scrP^{2,F,\bfw}(\bfx)$, differs from the fibre product orientation by
$(-1)^\triangle$, where
\begin{equation}
\begin{aligned}
\triangle = & \; |\{ (f_+,f_-) \in F_+ \times F_- \subset F \times F \;:\; f_+ > f_-\}| \\ & +
|F_-|(\deg(x^{\mathit{new}}) - \deg(x^0)).
\end{aligned}
\end{equation}
\end{lemma}

Now we turn to the moduli spaces of cascade maps $\scrQ^{2,F,\bfw}(\bfx)$ from Section \ref{subsec:chain}. Recall that a point in such a space consists of a finite chain of $l \geq 1$ components, each of which is of the form \eqref{eq:linear-cascade}. As before, we choose an ordering of $F$, equip the subsets $F_j$ with the induced orderings, and consider the associated isomorphisms
\begin{equation} \label{eq:component-o}
\Lambda^{top}(T\scrP^{2,F_j,(w_{j-1},w_j})(x_{j-1},x_j) \otimes o_{x_j} \longrightarrow
o_{x_{j-1}}.
\end{equation}
Now $\scrQ^{2,F,\bfw}(\bfx)$ is locally equal to the product of the $\scrP^{2,F_j,(w_{j-1},w_j)}(x_{j-1},x_j)$. Hence, by composing all the \eqref{eq:component-o},
one gets an isomorphism $\Lambda^{top}(T\scrQ^{2,F,\bfw}(\bfx)) \otimes o_{x^1} \rightarrow o_{x^0}$. We twist this isomorphism by $(-1)^\star$, where
\begin{equation} \label{eq:star-sign}
\begin{aligned}
\star = & \; |\{ (f_+,f_-) \in F_j \times F_k \subset F \times F \;:\; j>k, \; f_+ > f_-\}| \\
& + \sum_j |F_j|\, (\deg(x_{j-1}) - \deg(x_0))\; +\; l,
\end{aligned}
\end{equation}
and denote the outcome by $o_{\{(\rho_j,Z,\bfphi_j,u_j)\}}^{red}$. Suppose that we change the ordering of $F$ by swapping two elements. If both of these lie in the same subset $F_j$, the associated isomorphism \eqref{eq:component-o} changes sign by definition, while in the other case, each component still gives rise to the same map, but $\star$ changes. In both cases, the result is that the sign of $o_{\{(\rho_j,Z,\bfphi_j,u_j)\}}^{red}$ gets flipped, which explains the cancellations which result in Lemma \ref{th:permute-cascade}.

We now concentrate on the case where $|F| \leq 1$, when the choice of ordering of $F$, hence also the first line in \eqref{eq:star-sign}, become trivial. For simplicity, assume as before that we have trivialized $o_{x^0}$ and $o_{x^1}$, so that the moduli spaces of cascade maps connecting these two chords are canonically oriented. Suppose that $\scrQ^{2,F,\bfw}(\bfx)$ is one-dimensional, and consider a boundary point in its partial compactification, where two components $\rho_{j-1} = \rho_j$ occur for the same parameter value. As in the model case of Example \ref{th:two-boundary}, such points can arise in two different ways: either $\bar{E} = \emptyset$, which means that equality of parameter values is achieved simply as a degenerate limiting case of the usual inequality $\rho_{j-1} < \rho_j$; or else $\bar{E} = E$, which means that our boundary point is the limit of a sequence in which one component is divided into two parts by (the standard Floer-theoretic) bubbling. In the first case, the orientation of the boundary point is just the fibre product orientation of $\scrQ^{2,F_+,\bfw_+}(\bfx_+) \times_{(0,1]} \times \scrQ^{2,F_-,\bfw_-}(\bfx_-)$, where
\begin{equation} \label{eq:blah}
\parbox{30em}{
$F_+ = F_1 \cup \cdots \cup \cdots F_{j-1}$, $F_- = F_j \cup \cdots \cup F_l$, $\bfw_+ = (w^0,w^{\mathit{new}})$, $\bfw_- = (w^{\mathit{new}},w^1)$, $\bfx_+ = (x^0,x^{\mathit{new}})$, $\bfx_- = (x^{\mathit{new}},x^1)$ with $w^{\mathit{new}} = w^0 - |F_+| = w^1 + |F_-|$, $x^{\mathit{new}} = x_j$.
}
\end{equation}
Let's temporarily forget the signs \eqref{eq:star-sign} and use the untwisted orientations of the moduli spaces of cascade maps. Then, Lemma \ref{th:compare-fibre-product} would tell us that the boundary orientation in the first and second cases differ from each other by $|F_-|(\deg(x^{\mathit{new}} - \deg(x^0))$. Turning the twist by \eqref{eq:star-sign} back on cancels out that difference, up to an overall $-1$ which comes from the difference in $l$ (obvious, since bubbling raises the number of components by one). Hence, we have shown that with the conventions used in defining $o_{\{(\rho_j,Z,\bfphi_j,u_j)\}}^{red}$, the contributions from the two boundary points to the standard counting argument will vanish, as claimed in Section \ref{subsec:chain}.

The remaining boundary points of $\bar\scrQ^{2,F,\bfw}(\bfx)$ are those where $\rho_l = 1$. In we take $j = l$ and use notation as in \eqref{eq:blah}, the boundary orientation differs from that of the product $\scrQ^{2,F_+,\bfw_+}(\bfx_+) \times \scrR^{2,F_-,\bfw_-}(\bfx_-)$ by $(-1)^\times$, where the sign comes from \eqref{eq:star-sign} and is
\begin{equation}
\times = |F_-|(\deg(x^{\mathit{new}}) - \deg(x^0)) + 1.
\end{equation}
The other relevant geometric phenomenon are the ends of the partial compactification, which are modelled on $(0,\delta] \times \scrR^{2,F_+,\bfw_+}(\bfx)^{in} \times \scrQ^{2,F_-,\bfw_-}(\bfx_-)$ with data as in \eqref{eq:blah} for $j = 1$. In this case, pointing towards infinity along such an end corresponds to decreasing the parameter. Hence, if we cut off the end, the boundary orientation differs from the product orientation by $(-1)^{\times+1}$. Specializing to $F = \emptyset$ and $F = \{1\}$, and taking into account the fact that the component of each boundary point belong to zero-dimensional moduli spaces, then yields \eqref{eq:q0-relation} and \eqref{eq:q1-relation}, respectively.

%\bibliographystyle{plain}
%\bibliography{../../../bib/all,../../../bib/new}

\end{document}